\documentstyle{amsppt}

\topskip=-8mm
\hoffset=3.2mm



\hfuzz=5pt
\mathsurround 2pt
\overfullrule=0pt

\define\A{{\Bbb A}}
\define\C{{\Bbb C}}

\define\R{{\Bbb R}}
\define\Q{{\Bbb Q}}

\define\Z{{\Bbb Z}}



\define\p{\frak p}

\define\a{\alpha}
\redefine\b{\beta}

\define\e{\varepsilon}   
\redefine\l{\lambda}
\redefine\o{\omega}
\define\ph{\varphi}

\define\s{\sigma}
\redefine\P{\Phi}
\predefine\Sec{\S}
\redefine\L{\Lambda}
\redefine\S{{\Cal S}}


\define\sh{\sharp}


\define\back{\backslash}

\define\lra{\longrightarrow}
\redefine\tt{\otimes}
\define\scr{\scriptstyle}
\define\liminv#1{\underset{\underset{#1}\to\leftarrow}\to\lim}
\define\limdir#1{\underset{\underset{#1}\to\rightarrow}\to\lim}

\define\isoarrow{\ {\overset{\sim}\to{\longrightarrow}}\ }


\define\nass{\noalign{\smallskip}}

\define\End{\text{\rm End}}
\define\Sym{\text{\rm Sym}}
\define\Hom{\text{\rm Hom}}
\define\Spec{\text{\rm Spec}\,}

\font\cute=cmitt10 at 12pt
\font\smallcute=cmitt10 at 9pt
\define\kay{{\text{\cute k}}}
\define\smallkay{{\text{\smallcute k}}}
\define\OK{\Cal O_{\smallkay}}

\define\sig{\text{\rm sig}}
\redefine\ord{\text{\rm ord}}

\define\degh{\widehat{\text{\rm deg}}}  


\magnification=\magstep1
\baselineskip=15pt
\parindent=0pt
\parskip=13pt

\redefine\ord{\text{\rm ord}}

\predefine\oldvol{\vol}
\redefine\vol{\text{\rm vol}}

\define\Diff{\text{\rm Diff}}

\define\OC{{\Cal O_C}}
\define\F{\Bbb F}
\define\fps{{\Bbb F_{p^2}}}
\define\fp{{\Bbb F_p}}
\redefine\L{\Cal L}
\define\K{\Cal K}

\define\CZ{\Cal Z}
\define\CB{\Cal B}              
\define\CT{\Cal T}              

\define\tr{\text{\rm tr}}
\define\diag{\text{diag}}

\define\pr{\text{\rm pr}}       

\redefine\j{{\text{\bf j}}}
\redefine\x{{\text{\bf x}}}

\define\mss{\Cal M^{\text{ss}}}
\define\red{\text{\rm red}}
\define\nn{\frak n}


\define\LL{{\Lambda}}            
\define\zps{{\Z_{p^2}}}          
\define\zpsx{{\Z_{p^2}^\times}}  
\define\qps{{\Q_{p^2}}}          

\define\bob{{\b\!\!\!\!\b}}      

\define\phn{\ph_p'}         


\define\eichler{{\bf 1}}
\define\vandergeer{{\bf 2}}
\define\grosskeating{{\bf 3}}
\define\hlr{{\bf 4}}
\define\hirzebruch{{\bf 5}}
\define\kaiser{{\bf 6}}
\define\kitaokatwo{{\bf 7}}
\define\kitaoka{{\bf 8}}   
\define\kottwitz{{\bf 9}} 
\define\duke{{\bf 10}}     
\define\annals{{\bf 11}}         
\define\krsiegel{{\bf 12}}
\define\krdrin{{\bf 13}}       
\define\tiny{{\bf 14}}
\define\langlands{{\bf 15}}
\define\milne{{\bf 16}}
\define\raynaud{{\bf 17}}
\define\stamm{{\bf 18}}
\define\yang{{\bf 19}}


\define\GSpin{\text{\rm GSpin}}
\define\teta{\tilde{\eta}}  
\define\sgn{\text{\rm sgn}}

\define\mud{\ ; \ }   

\vsize=7.7in
\voffset=-.4in

\centerline{\bf Arithmetic Hirzebruch Zagier cycles}
\smallskip
\centerline{by}
\smallskip
\centerline{\bf Stephen S. Kudla\footnote{Partially supported by NSF Grant
DMS-9622987}}
\smallskip
\centerline{and}
\smallskip
\centerline{\bf Michael Rapoport}

\vskip 25pt

\subheading{Introduction}

In this paper and its companion \cite{\krsiegel}, we establish a relation
between part of the height pairing of special cycles on a Shimura variety
associated to an orthogonal group of signature $(n-1,2)$ and special
values of derivatives of Fourier coefficients of certain Siegel Eisenstein series
of genus $n$ in two new cases. Our results provide some evidence in favor of the
general program set forth in \cite{\annals}, where the case of
Shimura curves was analyzed. In \cite{\krsiegel}, we considered the case
of (twisted) Siegel
threefolds ($n=4$) while, in the present paper, we are concerned
with the case of (twisted) Hilbert-Blumenthal surfaces ($n=3$).

For these surfaces, the special cycles include the modular and Shimura curves
 studied extensively by Hirzebruch and Zagier (comp.\
\cite{\hirzebruch}, \cite{\vandergeer}).
 Among other things,
they showed that a generating function for the intersection numbers of
such curves is an elliptic modular form of weight $2$.
The Hirzebruch-Zagier curves also account
for all Tate classes on such surfaces rational over
abelian extensions of $\Q$, \cite{\hlr}.
Thus it is of interest to investigate the arithmetic analogues of these cycles
on integral models of such surfaces, which are arithmetic threefolds.

The canonical model $M$ over $\Q$ of a Hilbert-Blumenthal surface
is defined as a moduli space of abelian varieties with real
multiplications, and the Hirzebruch-Zagier curves on $M$ can be
defined as the loci where the abelian variety admits an extra
endomorphism of a particular type, a special endomorphism. For a
prime $p$ of good reduction, a model $\Cal M$ of $M$ over
$\Z_{(p)}$ can be defined by a moduli problem. We then consider a
modular extension $\Cal Z$ to $\Cal M$ of a Hirzebruch-Zagier curve
$Z$ on $M$, again defined by imposing a special endomorphism. These
are the arithmetic Hirzebruch-Zagier cycles of the title. More
generally, one can consider the loci where the abelian variety
carries several special endomorphims. If two independent
endomorphisms are imposed one obtains points on the generic fiber
$M$ and curves on $\Cal M$. If three independent endomorphisms are
imposed, then the associated cycle $\Cal Z$ is supported in the
special fiber, but {\it need not have dimension $0$}! In fact, it
is sometimes possible to impose $4$ independent special
endomorphisms and still have a nonempty locus, which may have
dimension $0$ or $1$.

It turns out that the case where $p$ splits in the real quadratic field differs
radically from the case where $p$ is inert. The first case is
similar to the case of modular curves, and in this case a special cycle $\Cal Z$
cut out by three independent special endomorphisms
consists of a finite number of points. The more interesting
case is when $p$ is inert.
In this case, we
determine the precise conditions which ensure that $\Cal Z$ consists of a finite
number of points.

For either type of $p$, when $\Cal Z$ is an
Artin scheme, we show that its length
coincides with the derivative at $s=0$, the center of symmetry, of a
Fourier coefficient of a
Siegel Eisenstein series of genus $3$ (Theorem~7.3 and Theorem~11.5).
This result, which generalizes
Theorem~3 of \cite{\tiny}, is connected with local height pairings
as follows. Suppose that $\Cal Z_1$ (resp.
$\Cal Z_2$) is an arithmetic cycle defined by imposing $1$ (resp. $2$) special
endomorphisms. The intersection of $\Cal Z_1$ and $\Cal Z_2$ is then the locus
where a triple of endomorphisms is imposed, although the relative position
of the first with respect to the second two is no longer specified. The
locus where the triple is independent is then a finite union
of special cycles $\Cal Z_3$. The sum of the lengths of the corresponding
local rings then gives the local contribution to the height
pairing at $p$ for the cycles $Z_1$ and $Z_2$
on the generic fiber, i.e., for a curve and a point on a
Hilbert-Blumenthal surface. In fact, for this interpretation one must
assume a conjecture explained in \cite{\krsiegel} concerning the vanishing
of Tor terms in the intersection multiplicity. The case of a height pairing of
a triple of curves can be handled in the same way.
These results generalize Theorem~14.11 of \cite{\annals}.

We now give a more precise discussion of our results.

The twisted Hilbert-Blumenthal surfaces can be viewed as the Shimura varieties
associated to certain rational quadratic forms.
Let $(V,Q)$ be a quadratic space of signature $(2,2)$
over $\Q$. The even part $C^+(V)$ of the Clifford algebra $C(V)$ is
the base change of an indefinite quaternion algebra $B_0$ over $\Q$
to a real quadratic extension $\kay$ of $\Q$. We will constantly
use the exceptional isomorphism
$$\GSpin(V) =G(\Q)=\{ g\in C^+(V)^{\times};\
\nu(g)\in\Q^{\times}\}\ \ .\tag I.1$$
Here (I.1) is the group of $\Q$-rational points of an algebraic
group $G$ over $\Q$ to which there is associated a Shimura variety
with complex points given as follows:
$$M(\C)= Sh(G,{\Cal D})_K= G(\Q)\setminus [{\Cal D}\times G({\Bbb
A}_f)/K]\ \ .\tag I.2$$
Here $K$ is a compact open subgroup of the finite adele group $G(\A_f)$
and $\Cal D$ is the space of oriented negative $2$-planes in $V(\R)$.
This is a twisted version of a
Hilbert-Blumenthal surface. In the degenerate case
$\kay=\Q\oplus\Q$, $M$ is a product of two modular curves; when
$\kay$ is a field and $B_0=M_2(\Q)$, $M$ is the usual
Hilbert-Blumenthal surface associated to $\kay$. The Shimura
variety (I.2) is the moduli space of principally polarized abelian
varieties of dimension 8 with level structure which are equipped
with an action of $C(V)\otimes\kay$ satisfying certain
compatibilities. Let us fix a prime number $p$ for which there
exists a self-dual $\Z_{(p)}$-lattice $\Lambda$ in $V(\Q)$ (good
reduction condition). We will always take $K$ to be of the form
$K=K^p.K_p$ where $K_p\subset G(\Q_p)$ is the stabilizer of
$\Lambda$ and where $K^p\subset G({\Bbb A}_f^p)$ is sufficiently
small. The modular interpretation of the Shimura variety then allows us
to construct a smooth model ${\Cal M}$ of (I.2) over
${\roman{Spec}}\, \Z_{(p)}$.

Hirzebruch and Zagier have defined special curves on
Hil\-bert-Blu\-men\-thal surfaces. Expressed in adelic language, a
prototype of such a curve is given by the inclusion of Shimura
varieties
$$Sh(B_0^{\times}, {\Cal D}_0)_{K_0}\hookrightarrow Sh(G, {\Cal
D})_K\ \ .\tag I.3$$ Here $B_0$ is as before an indefinite
quaternion algebra over $\Q$, $K_0$ is a compact open subgroup of
the finite adeles $B^\times_{0,\A_f}$ and $\Cal D_0$ is the space
of oriented negative $2$-planes in $V_0(\R)$, where $V_0$ is the
space of trace zero elements in $B_0$. In \cite{\hlr} the
Hirzebruch-Zagier cycles are defined as the images of (I.3) under
the Hecke correspondences defined by elements in $G({\Bbb A}_f)$.
Note, however, that Hirzebruch and Zagier do not use the adelic
language and that their definition of these cycles is different.

In the present paper, we use yet
another version of these cycles. As indicated above,  we define the special
cycles in a modular way by imposing {\it special endomorphisms} on
the abelian varieties parametrized by ${\Cal M}$. Here an
endomorphism is special if it is self-adjoint for the Rosati
involution and if it Galois commutes with the action of
$C(V)\otimes\kay$,
$$\iota(c\otimes a)\circ j=j\circ\iota(c\otimes a^{\sigma})\ \ ,\ \
c\in C(V),\ a\in\kay;\ <\sigma>\ = \ {\roman{Gal}}(\kay /\Q)\ \
.\tag I.4$$
The space of special endomorphisms is a finitely
generated free $\Z_{(p)}$-module\break
equipped with the quadratic form $Q$ given by
$$j^2=Q(j)\cdot 1.$$
For $T\in{\roman{Sym}}_n(\Z_{(p)})_{\geq 0}$ and a $K$-invariant
open compact subset $\omega\subset V({\Bbb A}_f)^n$ (with a certain
integrality property at $p$), we define as in \cite{\krsiegel} a
special cycle $\CZ(T,\o)$ by imposing an $n$-tuple $\j$ of special
endomorphisms with $Q(\j)=T$ and satisfying a compatibility
condition with $\o$. If ${\roman{det}}(T)\neq 0$, the generic fibre
of $\CZ(T,\o)$ is a cycle of codimension $n$ on the
Hilbert-Blumenthal surface. For $n=1$ we essentially obtain the
classical Hirzebruch-Zagier curves. As in \cite{\krsiegel} our aim
is to study the structure of these cycles and their intersection
behaviour.

Let us fix positive integers $n_1,\ldots, n_r$ with $n_1+\ldots
+n_r=3$. For each $i=1,\ldots, r$ we choose $T_i\in
{\roman{Sym}}_{n_i}(\Z_{(p)})_{>0}$ and $\omega_i\subset V({\Bbb
A}_f)^{n_i}$ as above. We form the fibre product of the
corresponding special cycles,
$$\CZ=\CZ(T_1,\omega_1)\times_{\Cal M}\ldots\times_{\Cal
M}\CZ(T_r, \o_r)\ \ .\tag I.5$$ As in \cite{\krsiegel} we have a
disjoint sum decomposition according to the value of the {\it
fundamental matrix}
$$\CZ=\mathop{\coprod}\limits_{\matrix
\scriptstyle T\in{\roman{Sym}}_3(\Z_{(p)})_{\ge 0}\\
\scriptstyle {\roman{diag}}(T)= (T_1,\ldots, T_r)\endmatrix
}
\CZ(T,\o)\tag I.6$$
Here $\omega=\o_1\times\ldots\times \o_r$. Concerning this
decomposition we have the following results.
Let $T\in {\roman{Sym}}_3(\Z_{(p)})_{\ge
0}$.

\proclaim{Theorem 1} If $\det T\ne 0$, then $\CZ(T,\o)$ lies
purely in characteristic $p$ and consists only of supersingular
points.
\endproclaim

\proclaim{Theorem 2} Suppose that ${\roman{det}}\, T\ne 0$.
If $p$ is split in $\kay$, then $\CZ(T,\o)$ is a union of isolated
points. If $p$ is inert in $\kay$, then $\CZ(T,\o)$ is a union of
isolated points
 if and only if $T\not\equiv 0\, {\roman{mod}}\,
p$.
\endproclaim

Suppose that ${\roman{det}}\, T\ne 0$ and that $T\not\equiv 0\,
{\roman{mod}}\, p$ when $p$ is inert in $\kay$. Define
$$<\CZ(T,\o)>_p=\sum\limits_{\xi\in\CZ(T,\o)(\overline{\Bbb F}_p)}
lg({\Cal O}_{\CZ(T,\o),\xi})\ \ .$$
In sections 7 and 11, we define a certain Siegel-Eisenstein series
$E(h,s,\P)$ on $Sp_6(\A)$  depending on $\o$ which vanishes
at $s=0$, the center of the critical strip. Our next result
concerns the Fourier coefficients $E'_T(h,0,\P)$of the derivative $E'(h,0,\P)$
at this point.
\proclaim{Theorem 3} For $T$ as above and for $h\in Sp_6(\R)$,
$$E'_T(h,0,\P) = -\frac12\cdot C\cdot \log(p)\cdot <\CZ(T,\o)>_p
\cdot W_T^2(h),$$
where $C$ is a volume factor and $W_T^2(h)$
is a standard archimedean Whittaker function.
\endproclaim

As explained above, via the decomposition (I.6),
Theorems 1--3  have a direct bearing on
the problem of intersection multiplicities of special cycles. The
proofs here are quite similar to those in \cite{\krsiegel}
and in fact easier since we are here in one dimension less.
Ultimately, as in loc.\ cit., we rely for the proof of Theorem 3 on the
results of Gross and Keating, \cite{\grosskeating},
on deformations of homomorphisms of
one-dimensional formal groups and of Kitaoka, \cite{\kitaoka}, on local
representation densities of quadratic forms.

These results for varying $p$ can be combined. In fact we define
a moduli problem represented by a scheme over $\Spec\  \Z[\frac{1}{N}]$,
for a suitable integer $N$, and special cycles $\CZ(T,\o_N)$
on it, whose base change to $\Spec\ \Z_{(p)}$ are the cycles discussed
above. Theorem 3 then shows that part of the Fourier expansion of
$E'(h,0,\P)$ is a generating series for the degrees of certain
of these cycles. In light of the result of \cite{\tiny}, one might
hope that the whole Fourier series has an interpretation of this kind.

The previous three theorems focus on the supersingular locus of the
Hirzebruch-Zagier cycles. However, we also obtain a more global
view of their geometry. Using the well known list of isogeny types
in the special fiber, we are able to enumerate those that meet the
image of a special cycle. Specifically, suppose that
$T\in\Sym_n(\Z_{(p)})$ with $\det(T)\ne 0$. If $n=3$, the
corresponding special cycle $\CZ(T,\o)$ lies entirely in the
supersingular locus, while if $n=2$, the special cycle can meet at
most one additional isogeny class, uniquely determined by $T$.
Furthermore, a non-supersingular point of ${\Cal M}$ in
characteristic $p$ is ordinary and we may use the Serre-Tate
deformation theory to investigate the local structure of special
cycles along the ordinary locus. Again it turns out that the
divisibility of $T$ by $p$ is decisive.
\proclaim{Theorem 4} Let $\overline{\CZ}(T,\o)^{\ord}$ be the ordinary
part of the special fiber of the cycle $\CZ(T,\o)$. \hfill\break
(i) Suppose that $n=1$, i.e., that $T\in \Z_{(p),>0}$. Then
locally at each point, $\overline{\CZ}(T,\o)^{\ord}$
is the $p^{\ord_p(T)}$-fold multiple of a smooth divisor. \hfill\break
(ii) Suppose that $n=2$, i.e., that $T\in \Sym_2(\Z_{(p)})_{>0}$.
Then $\overline{\CZ}(T,\o)^{\ord}$ is an Artin scheme of length
$p^{\ord_p(\det(T))}$ at each of its finitely many points.
\endproclaim
Note that by Theorem 1, if $T\in \Sym_n(\Z_{(p)})_{>0}$ with
$n\ge 3$, then $\overline{\CZ}(T,\o)^{\ord}$ is empty.

Now suppose that $p$ is inert in $\kay$ and that $p|T$. Then
by Theorem 2, the special cycle $\CZ(T,\o)$ cut out
by three independent
special endomorphisms has dimension $1$ (case of
excess intersection).
The supersingular locus of ${\Cal M}$ is a
union of projective lines, \cite{\stamm}, and
$\CZ(T,\o)$ is a union of
certain of these. We give a complete enumeration of the lines in
$\CZ(T,\o)$ in
terms of the Bruhat-Tits building $\Cal B$ of $G_{ad}(\Q_p)\simeq
PGL_2(\Q_{p^2})$.
\proclaim{Theorem 5}
Suppose that $T=\text{\rm diag}(\e_1 p^{a_1},\e_2 p^{a_2},\e_3 p^{a_3})$,
where $\e_i\in \Z_p^\times$ and $0\le a_1\le a_2\le a_3$.
Associated to $T$ is a triple $\bob=(\b_1,\b_2,\b_3)$
of anticommuting automorphisms of $\Cal B$.
Then, the dual graph to any connected component of $\CZ(T,\o)$ is the
set of points in $\Cal B$ which satisfy
$$d(x, \Cal B^{\b_i}) \le \frac12 (a_i-1)\qquad\qquad i=1,\ 2,\ 3,$$
where $d(x,\Cal B^{\b_i})$ is the distance from $x$ to the fixed
point set $\Cal B^{\b_i}$ of $\b_i$. In particular, the number of
irreducible components of each connected component is the number of
vertices $|\Cal T(\bob)_0|$ of $\Cal B$ lying in the intersection
of these tubes, with specified radii, around the three fixed point
sets. Furthermore, each connected component consists of a single
projective line if and only if $a_1=a_2=1$ and $a_3$ is odd, and,
in addition, $a_3=1$ if $-\e_2\notin \Z_p^{\times,2}$.
\endproclaim

An explicit formula for the number $|\Cal T(\bob)_0|$
can be obtained in general by
a combinatorial effort, cf. for example, (8.19) and
(8.25).
On the other hand, $|\Cal T(\bob)_0|$
can be computed as a twisted orbital integral, (8.30),
and this integral can be expressed, in turn, as a representation density
of quadratic forms! Combining these facts, we obtain the formula
(Theorem~8.15)
$$|\Cal T(\bob)_0|= (1-p^{-4})^{-1} \cdot \a_p(S,p^{-1} T),\tag I.7$$
which can be viewed as giving an explicit formula for the representation
density on the right hand side. Such explicit formuli are of independent
interest.
Since the number of connected components can also be given
as an orbital integral, we obtain the following result (Proposition~8.17):
$$
\#\left\{ \matrix\text{irreducible components}\\
\text{of}\ \CZ(T,\o)\endmatrix\right\}
= 2\cdot \vol(K)^{-1}\cdot TO_T(\ph_p'')\cdot O_T(\ph^p_f).\tag I.8$$
This formula is reminiscent of the expressions arising for the number of
points of Shimura varieties over finite fields. It
also suggests that there may be a modular generating function for the numbers of
such components.

In the case where $p$ is inert and $p|T$, it is an important open
problem to find an analogue of
Theorem 3 giving the contribution of $\CZ(T,\o)$ to
an intersection multiplicity.
A result of this type was obtained in \cite{\krdrin} in the
somewhat analogous case of bad reduction of Shimura curves. Note, however, that
the p-adic uniformization which is available there is not available in our
present case. The
first degenerate cases that could be investigated are those where
each connected component of $\CZ(T,\o)$ is a
single projective line, as described in Theorem 5.

A similar analysis of the irreducible components of the degenerate
cycles can be done in the Siegel case, to which we return in
section 9. Here again, the projective lines in the supersingular
locus lying in the image of a special cycle can be described in
terms of the building of $G_{ad}(\Q_p)$, where $G$ is now a rank
$1$ twisted form of the symplectic group of rank $2$ over $\Q_p$.
The description of the irreducible components again involves an
analysis of the fixed point sets of certain anticommuting
involutions, but this analysis is now considerably more
complicated, and we draw heavily on the work of Kaiser on fixed
point sets of tori in the building \cite{\kaiser}. Nevertheless, we
obtain a rather complete description. For example, Proposition~9.8
gives the precise conditions on $T$ under which each connected
component of a degenerate cycle (now associated to a $4$-tuple of
special endomorphisms) is irreducible -- the analogue of the last
statement of Theorem 5. The same method also works in the case of
$p$-adic uniformization \cite{\krdrin}. In all of these cases, the
fact that $rk_{\Q_p}G_{ad}=1$ seems absolutely crucial.

It seems to us that our results on degenerate special cycles
constitute only the beginning of a circle of very interesting
problems. This also explains the different and much more open-ended
nature of this paper from its companion \cite{\krsiegel}. We hope that
it can serve as a basis of further investigations.

We now give an overview of the structure of this paper. Section 0
contains preliminaries on linear algebra, in particular the
exceptional isomorphism mentioned above. In section 1 we define the
Shimura variety, the associated moduli problem and introduce the
special cycles. In section 2 we extend these concepts to construct
a model at a prime of good reduction. From this point through
section 10, we assume that $p$ is inert in $\kay$. Section 3
contains an enumeration of the isogeny classes in the special fibre
of ${\Cal M}$ and their corresponding $\Q$-vector spaces of special
endomorphisms. Section 4 is a presentation of the results of Stamm
\cite{\stamm} on the structure of the supersingular locus of ${\Cal
M}\times_{{\roman{Spec}}\, \Z_{(p)}} {\roman{Spec}}\, {\Bbb F}_p$.
In section 5 we determine the $\Z_p$-module of special
endomorphisms of a supersingular Dieudonn\'e module. This is then
used in section 6 to prove our first main result (Theorem 6.1)
which characterizes the special cycles with isolated supersingular
points, cf. Theorem 2 above. In this section we also determine (by
reduction to the theorem of Gross/Keating) the length of the local
ring at an isolated point of a special cycle. Section 7 gives the
relation to Eisenstein series and proves Theorem 3 above (Theorem
7.3 and Corollary 7.4). In section 8 we give the analysis of the
set of irreducible components of $\CZ(T,\o)\cap {\Cal M}^{ss}$ and
in section 9 we treat the analogous question in the Siegel case. In
section 10 we consider the ordinary locus of special cycles. We
also investigate the compactness property of special cycles. In
section 11, we consider the case of a prime $p$ which splits in
$\kay$. In section 12, we consider special cycles on a moduli
scheme over $\Spec\ \Z[\frac{1}{N}]$ and the generating function
for their degrees. The final section 13 contains some remarks on
the interrelation of the special cycles in the various moduli
problems considered in our series of papers \cite{\annals},
\cite{\krsiegel}, \cite{\krdrin}, \cite{\tiny}.

In conclusion we wish to thank Ch.\ Kaiser for
considerable help with section 9. In particular, he corrected
some errors in an earlier version and indicated to us how to
obtain more complete results.
We also thank A.\ Langer for interesting
conversations on Hirzebruch-Zagier cycles. Stephen Kudla would like
to thank the University of Cologne for its hospitality during June
1997. Michael Rapoport would like to thank the University of
Maryland for its hospitality during March 1998. The support of the
NSF and the DFG is gratefully acknowledged.
\bigskip\noindent
{\bf Table of contents}
\roster
\item"{0.}" Preliminaries on linear algebra
\item"{1.}" Quaternionic Hilbert-Blumenthal surfaces and special cycles
\item"{2.}" Models over $\Z_{(p)}$ and special cycles
\item"{3.}" Isogeny classes and special endomorphisms
\item"{4.}" The structure of the supersingular locus
\item"{5.}" Special endomorphisms of supersingular Dieudonn\'e
models
\item"{6.}" Isolated supersingular points of special cycles
\item"{7.}" Representation densities and Eisenstein series
\item"{8.}" Components of $\CZ(T,\o)\cap {\Cal M}^{ss}$
\item"{9.}" Components in the Siegel case
\item"{10.}" Special cycles on the special fibre
\item"{11.}" The case of a split prime
\item"{12.}" A global model and generating series
\item"{13.}" On the hereditary nature of special cycles
\endroster
\bigskip\noindent
{\bf Notation:} We use $\Z_{p^2}=W({\Bbb F}_{p^2})$,
the ring of Witt vectors of the finite field $\Bbb F_{p^2}$ and
$\Q_{p^2}=\Z_{p^2}\otimes_{\Z_p}\Q_p$. We fix a unit
$\Delta\in\Z_p^{\times}-\Z_p^{\times, 2}$ and write
$$\Z_{p^2}=\Z_p[\delta] /(\delta^2-\Delta)\ \ .$$
By $\chi$ we denote the quadratic residue character on
$\Z_p^{\times}$ resp.\ ${\Bbb F}_p^{\times}$.
We also let $\Bbb F$ be a fixed algebraic closure of $\Bbb F_p$,
$W=W(\Bbb F)$, its ring of Witt vectors, $\Cal K= W\tt_{\Z}\Q$
the quotient field of $W$ and $\s$ their Frobenius endomorphism.

\subheading{\Sec0. Preliminaries on linear algebra}

In this section we begin by collecting some facts about the Clifford algebras of
$4$-dimensional quadratic spaces over a field $F$ which
will be used in the rest of the paper. In particular we will give a
realization of the accidental isomorphism mentioned in the
introduction, (I.1), comp.\ \cite{\eichler}, pp.\ 31-33. We then specialize to the
case $F=\Q$. When the signature of $V$ is $(2,2)$, we construct the data needed
to attach a Shimura variety to $V$.

Let $(V,Q)$ be a quadratic space of dimension 4 over a field $F$ of
characteristic not 2 and let $C(V)=C^+(V)\oplus C^-(V)$ be its
Clifford algebra. For a choice of an orthogonal basis $v_1,\ldots,
v_4$ of $V$ over $F$, with $Q(v_i)=a_i$, the element $\delta =
v_1v_2v_3v_4\in C^+(V)$ satisfies $\delta^2=a_1a_2a_3a_4$. Up to
multiplication by an element of $F^{\times}$, $\delta$ is
independent of the choice of basis and
${\roman{det}}(V):=\delta^2\in F^{\times}/F^{\times, 2}$ depends
only on $V$. The algebra $\kay:= F(\delta)$ is the center of
$B:=C^+(V)$ and $B$ is a quaternion algebra over $\kay$. Let $\iota$
be the main involution of $C(V)$ which reduces to the identity on
the elements of $V\subset C^-(V)$. Then $\iota$ fixes $\delta$ and
hence induces the identity on $\kay$. Furthermore it is obvious
that the $+1$- resp.\ $-1$-eigenspaces of the action of $\iota$ on
$C^-(V)$ are given by
$${\roman{span}}\, \{ v_1, v_2, v_3, v_4\}\ ,\ \text{resp.}\
{\roman{span}}\, \{ v_2v_3v_4, v_1v_3v_4, v_1v_2v_4, v_1v_2v_3\}\
\ .$$
Hence we can recover $V$ from $C(V)$ as follows.
\proclaim{Lemma 0.1}
$$V=\{ x\in C^-(V);\ x^{\iota}=x\}\ \ .$$
Furthermore, for $x\in V$,
$$Q(x)=x^2= x\cdot x^{\iota}=\nu (x)\ \ ,$$
where $\nu$ denotes the spinor norm.
\qed
\endproclaim
\proclaim{Lemma 0.2}
Let $ v_0\in V$ with $Q( v_0)=\alpha\ne 0$. Let $\sigma=
{\roman{Ad}}( v_0)$ denote the adjoint automorphism of $C(V)$
induced by $ v_0$.
\roster
\item"{(i)}" $\delta^{\sigma}=-\delta.$
\item"{(ii)}" The fixed algebra $B_0$ of $\sigma$ in $B=C^+(V)$ is
a quaternion algebra over $F$ such that $B\simeq B_0\otimes_F\kay$.
\endroster
\endproclaim

\demo{Proof}  We may suppose that $ v_0=v_4$ (by rechoosing our
basis if necessary).
 The relation $\delta v_0=- v_0\delta$ is then obvious. To check
(ii), consider the basis
 $$1,\ i= v_1v_2,\ j=v_2v_3,\ k=ij= a_2v_1v_3,\ \delta, \delta i,\
 \delta j,\
\delta k$$
 for $B$. Since $ v_0 v_i=-v_i  v_0$, for $i=1,\ 2,\ 3$, it is clear
 that $B_0$ is the span of the first $4$ basis vectors, and that
 $B=B_0\tt_FE$. Also, note that $$B_0 = ( -a_1a_2,
 -a_2a_3),\tag0.1$$
the cyclic algebra over $F$. \qed\enddemo

We introduce the algebraic group $G$ over $F$ with
$$G(R)=\{ g\in (B\otimes_FR)^{\times};\ \nu(g)= g\cdot
g^{\iota}\in R^{\times}\}\tag0.2$$ for any $F$-algebra $R$. Fix
$ v_0\in V$ with $\alpha=Q( v_0)\ne 0$ and let again $\sigma=
{\roman{Ad}}( v_0)$. Let
$$\tilde V=\{ x\in B;\ x^{\iota} =x^{\sigma}\}\ \ .\tag0.3$$
Then $G$ acts on $\tilde V$ via
$$g:x\mapsto g\cdot x\cdot g^{-\sigma}\ \ ,\ \ x\in \tilde V\ \ .$$
Indeed,
$$\align
(gxg^{-\sigma})^{\iota} &
=(g^{-\sigma})^{\iota}\cdot
x^{\iota}\cdot g^{\iota}\\ &
=\nu(g)^{-1}\cdot g^{\sigma}\cdot x^{\sigma}\cdot g^{\iota}\\
&
=g^{\sigma}\cdot x^{\sigma}\cdot g^{-1}\\
&
=(gxg^{-\sigma})^{\sigma}\ \ .
\endalign$$
On $\tilde V$ we have the quadratic form defined by
$$\tilde Q(x)= x\cdot x^{\sigma} = x\cdot x^{\iota} =\nu(x)\ \
.\tag0.4$$
This quadratic form is preserved by $G$,
$$\tilde Q(gxg^{-\sigma})= gxg^{-\sigma}\cdot (g^{\sigma}
x^{\sigma} g^{-1})= \tilde Q(x)\ \ .$$

\proclaim{Lemma 0.3}
The map $x\mapsto x\cdot v_0$ induces an isometry
$$(\tilde V, \alpha\, \tilde Q)\longrightarrow (V, Q)\ \ .$$
Under this isometry the group $G$ is identified with
$\GSpin(V,Q)$.
\endproclaim

\demo{Proof} In terms of the basis of $B$ used in the proof of
Lemma 0.2, $\tilde V$ is spanned by $1,\delta i, \delta j, \delta
k$. Hence the map is a linear isomorphism from $\tilde V$ to $V$.
For $x\in \tilde V$,
$$(x\cdot v_0)^2=x v_0\cdot x v_0 =xx^{\sigma}\cdot v_0^2
=\alpha\,\tilde Q(x)\ \ ,$$
hence the map is an isometry. Under this map the action of $G$ on
$\tilde V$ gets carried into the usual conjugation action on $V$.
In particular, $V$ is automatically preserved by this action. This
proves the last assertion.
\qed
\enddemo
\medskip\noindent
{\bf Remark 0.4.} Thus, we see that, starting with the quadratic
space $V$ and a non-isotropic vector $ v_0\in V$ with $\alpha=
Q( v_0)$, we obtain a quaternion algebra $B_0$ with $C^+(V)=
B_0\otimes_F\kay$. Conversely, let $\kay$ be a semisimple
$F$-algebra of dimension 2 with Galois automorphism $\sigma$. Let
$B_0$ be a quaternion algebra over $F$ and let
$B=B_0\otimes_F\kay$, with automorphism $\sigma=
{\roman{id}}\otimes
\sigma$. Now choose $\alpha\in F^{\times}$ and define an algebra
$C$ of dimension 16 over $F$ by
$$C=B\langle v_0\rangle\ \ ,$$
with the relations $ v_0^2= \alpha$ and $ v_0 b= b^{\sigma} v_0$, for
all $b\in B$. Extend the involution $\iota$ by $ v_0^{\iota}
= v_0$ and the automorphism $\sigma$ by $ v_0^{\sigma} = v_0$.
Finally, let
$$V=\{ x\in C^-=B\cdot v_0;\ x^{\iota} =x\}\ \ ,$$
with quadratic form $Q(x)=x^2$. Then starting with $(V,  v_0)$ we
recover $B_0$ and $\sigma$ and these two constructions are inverse
of one another in an obvious way.

From now on we take $F=\Q$ and assume that the signature of $V$ is
$(2,2)$. There is a Shimura variety associated to such a $V$, as we
now explain.

The choice of signature implies that $\delta^2 >0$, so that $\kay$ is a real
quadratic field or $\Q\oplus \Q$, and that $B$ is a totally
indefinite quaternion algebra over $\kay$.

\proclaim{Lemma 0.5} Let $ v_0\in V$ with $Q( v_0)=\alpha >0$ with associated fixed
algebra $B_0$, as above.
Let $\tau\in B_0^{\times}$ with $\tau^{\iota} =-\tau$ and $\tau^2
<0$. Then the involution
$$x\mapsto x^{\ast} =\tau x^{\iota} \tau^{-1}$$
is a positive involution of $C(V)$.
\endproclaim

\demo{Proof} Since $B$ is totally indefinite, if $\tau\in
B^{\times}$ with $\tau^{\iota} =-\tau$ and $\tau^2 <\!\!<0$ (totally
negative), then the involution $b\mapsto \tau b^{\iota}\tau^{-1}$
is a positive involution of $B$. If now $\alpha >0$ and $\tau\in
B_0^{\times}$ with $\tau^{\iota} =-\tau$ and $\tau^2 <0$, then
$${\roman{tr}}^0((b_1+b_2 v_0)\,\tau (b_1+b_2 v_0)^{\iota}
\tau^{-1})
={\roman{tr}}^0(b_1\,\tau b_1^{\iota} \tau^{-1})
+\alpha\,{\roman{tr}}^0 (b_2\,\tau b_2^{\iota}\tau^{-1}) >\!>0\
.\qed$$
\enddemo
Let $U_\Q=C(V)$, viewed as a $\Q$-vector space of dimension $16$,
and define
$$i:C(V)\tt_\Q\kay \rightarrow \End_\Q(U_\Q),\qquad i(c\tt
a)x=cxa.$$
Define a nondegenerate, $\Q$-valued, alternating form on
$U_\Q$ by
$$<x,y>\ =\ \tr^0(y^\iota \tau x),\tag0.5$$
where $\tr^0$ denotes the reduced trace on $C(V)$. Then
$$<i(c\tt a)x,y>\  =\  <x,i(c^*\tt a)y>.\tag0.6$$
Note that $c\tt a\mapsto c^*\tt a$ is a positive involution of $C\tt\kay$.
Moreover, the action of $G=\GSpin(V)$ on $U_\Q$ by right
multiplication commutes with the action of $C(V)\tt\kay$, and
preserves the form $<\ ,\ >$ up to a scalar:
$$\align
<xg,yg> &= \tr^0(g^\iota y^\iota \tau x g)\\ {}&=\tr^0(gg^\iota
y^\iota \tau x)\\ {}&=\nu(g)\, <x,y>,\endalign
$$
where $\nu(g)=gg^\iota$.

Let $\End(U_\Q,i)$ denote the commutant of $i(C(V)\tt\kay)$ in
$\End(U_\Q)$. Let ${\Cal D}$ be the space of oriented negative
$2$-planes in $V(\R)$. For $z\in {\Cal D}$ with properly oriented
orthogonal basis $z_1$,\! $z_2$, with $Q(z_i)=-1$, let $j_z =
z_1z_2\in C^+(V)_\R$. Then $j_z^2=-1$, and there is an algebra
homomorphism
$$h_z:\C\lra \End(U_\R,i)^{\text{\rm op}}\tag0.7$$
such that $h_z(\sqrt{-1}) = r(j_z)$ is right multiplication by $j_z$.
Note that $G(\R)$ is a subgroup of the invertible elements in
$\End(U_\R,i)^{\text{\rm op}}$. Via (0.7), the space $\Cal D$ can
be identified with a $G(\R)$-conjugacy class of homomorphisms
$h_z:\C^\times \rightarrow G(\R)$ satisfying the axioms required to
define a Shimura variety $Sh(G,\Cal D)$. This variety has a canonical model
over the reflex field $E(G,\Cal D)$.

Fix a point $z_0\in \Cal D$ and let
$U_\C=U_1\oplus U_2$ be the $\pm \sqrt{-1}$ eigenspaces of $h_{z_0}(\sqrt{-1})$. The
reflex field $E(G,\Cal D)$ is then the field of definition of the isomorphism
class of $U_1$ as a complex representation of $C(V)\tt\kay$.
\proclaim{Lemma 0.6} (i) $E(G,\Cal D)=\Q.$\hfill\break
(ii) For $c\in C(V)$ and $a\in \kay$,
$$\det(i(c\tt a);U_1) = N^0(c)^2 N_{\smallkay/\Q}(a)^4.$$
\endproclaim
\demo{Proof}
The space $V$ has a $\Q$-basis $v_1,\dots,v_4$ such that
the matrix for the quadratic form is $\diag(a_1,a_2,a_3,a_4)$, with
$a_1$, $a_2 >0$ and $a_3$, $a_4<0$. Let $z_j=
\frac{1}{\sqrt{|a_j|}}v_j$ and let $z$ be the oriented negative
$2$-plane spanned by $z_3$ and $z_4$. The element $j_z$ lies in
$C^+(V)\tt\Q(\sqrt{a_3a_4})$, and the $\pm \sqrt{-1}$ eigenspace of $j_z$
on $U_\C=C(V)\tt_\Q\C$ is preserved by any automorphism of $\C$
which is trivial on $\Q(\sqrt{a_3a_4},\sqrt{-1})$. In fact, for any
automorphism $\s$ of $\C$, we have either $U_1^\s= U_1$ or
$U_1^\s=U_2$. But there is an element $h\in SO(V)(\R)$ which acts
by $h:z_1 \leftrightarrow z_2$, and $h:z_3 \leftrightarrow z_4$.
There is an element $g\in G(\R)=\GSpin(V)(\R)$ whose action on $V$
by conjugation is $h$. Since $g j_z = - j_z g$, the action of $g$
on $U_\C$ by right multiplication switches the eigenspaces $U_1$
and $U_2$, and commutes with the action of $C(V)\tt\kay$. Thus
$U_1$ and $U_2$ are isomorphic as complex representations of
$C(V)\tt\kay$, and the reflex field is indeed $\Q$. Also, since
$C(V)_\C\simeq M_4(\C)$,
$$\det(i(c\tt a);U_\C) = N^0(c)^4 N_{\smallkay/\Q}(a)^8=
\det(i(c\tt a);U_1)^2.\qed$$
\enddemo

\subheading{\Sec1. Quaternionic Hilbert-Blumenthal surfaces
and special cycles}

In this section we construct the quaternionic Hilbert-Blumenthal surfaces $Sh(G,\Cal D)$ as
moduli spaces over $\Q$. We then define the special cycles, of codimension $1$ or $2$, on them by imposing
additional special endomorphisms analogous to those in the companion paper \cite{\krsiegel}.
For the standard Hilbert modular surfaces, we recover the Hirzebruch-Zagier curves.

We continue with the notation of the end of the last section. In
particular $(V,Q)$ is a quadratic space over $\Q$ of signature $(2,2)$ and
$\tau\in B= C^+(V)$ is an element with $\tau^{\iota}=-\tau$ and
$\tau^2\in\Q$ with $\tau^2 <0$. The map $x\mapsto x^{\ast}=\tau x^\iota \tau^{-1}$ is a
positive involution of $C(V)$.

We fix a compact open subgroup $K$ of $G({\Bbb A}_f)$ and consider
the following moduli problem over $\Q$. It associates
to a scheme $S\in \text{Sch}/\Q$ the set of isomorphism classes of
4-tuples $(A,\lambda,\iota,\bar{\eta})$, where
\roster
\item"{(i)}" $A$ is an abelian scheme over $S$, up to isogeny,
\item"{(ii)}" $\lambda:A\isoarrow \hat{A}$, is a $\Q^\times$-class of
polarizations on $A$,
\item"{(iii)}" $\iota:C(V)\tt\kay \lra \End^0(A)$ is a homomorphism such that
$$\iota(c\tt a)^* = \iota(c^*\tt a),$$
for the Rosati involution of $\End^0(A)$ and the involution $*$ of $C(V)$ introduced
above, and
\item"{(iv)}" $\bar{\eta}$ is a $K$--equivalence class of $C(V)\tt\kay$--equivariant isomorphisms
$$\eta:\hat{V}(A) \isoarrow U(\A_f)$$
which preserve the symplectic forms up to a scalar in
$\A_f^\times$.
\endroster
Here $\hat{V}(A) =\prod_\ell T_\ell(A)\tt\Q$.
We refer to \cite{\kottwitz} for the
precise explanation of these data.
In addition, we impose the determinant condition:
$$\det(\iota(c\tt a);\text{Lie}(A)) = N^0(c)^2 N_{\smallkay/\Q}(a)^4.
\tag1.1$$
In particular, $A$ has relative dimension $8$ over $S$.

For $K\subset G(\A_f)$ sufficiently small, this moduli problem is represented by a
smooth quasi-projective scheme $M=M_K$ over $\Q$. If $B=C^+(V)= \End(U_\Q,i)^{\text{\rm op}}$ is a division
algebra, then $M$ is in fact projective. Also, in general, $M$ is a canonical model of
the Shimura variety $Sh(G,h)_K$, and
$$M(\C) = Sh(G,h)_K(\C) \simeq G(\Q)\back {\Cal D}\times G(\A_f)/K,
\tag1.2$$ where ${\Cal D}$ is the space of oriented negative
$2$-planes in $V(\R)$, as in section 0 above.

In the case where $\kay$ is a (real) quadratic field and
$B=M_2(\kay)$, $M(\C)$ is a (union of) Hilbert-Blumenthal surfaces.
If $B$ is a division algebra, then $M(\C)$ is a (union of)
quaternionic Hilbert-Blumenthal surfaces. In the most degenerate
case, when $V$ is split over $\Q$, we have $\kay=\Q\oplus\Q$ and
$B=M_2(\Q)\times M_2(\Q)$, and $M(\C)$ is a product of modular
curves.

Next, we construct certain algebraic cycles on these varieties.
\proclaim{Definition 1.1} A special endomorphism of $U_\Q$ is an element
$j\in \End(U_\Q)$ such that, for $a\in \kay$ and $c\in C(V)$,
$$i(c\tt a) \circ j = j\circ i(c\tt a^\s),\tag sp.1$$
and
$$j^*=j,\tag sp.2$$
where $*$ denotes the adjoint with respect to $<\ ,\ >$, cf.\
(0.6).
\endproclaim

Since $\delta v = - v \delta$ for all $v\in V$, any endomorphism of
$U_\Q$ satisfying the first condition in Definition 1.1 is given by
right multiplication by an element of $C^-(V)$. Note that, for any
$c\in C(V)$,
$$<xc,y>\ =\ <x,yc^\iota>,\tag1.3$$
so that, by Lemma~0.1, a special endomorphism is precisely one
given by right multiplication by an element of $V$.
\proclaim{Definition 1.2} A special endomorphism of $(A, \lambda,
\iota, \overline\eta)\in M(S)$ is an element $j\in\End^0(A)$ such
that
$$\iota(c\otimes a)\circ j=j\circ \iota(c\otimes a^{\sigma})\ \ ,\
\ c\otimes a\in C(V)\otimes k$$
and $$j^*=j\ \ ,$$ for the Rosati involution of $\End^0(A)$.
\endproclaim
For a special endomorphism $j$ we have
$$j^2=Q(j)\cdot {\roman{id}},\ \ \hbox{with}\ Q(j)\in\Q\ \ ,$$
provided that the base scheme $S$ is connected, comp.\
\cite{\krsiegel}, Lemma 2.2. Therefore the space of special
endomorphisms of $(A, \lambda, \iota, \overline\eta)$ is a
quadratic space in this case. The positivity of the Rosati
involution implies that this quadratic space is positive-definite,
comp.\ \cite{\krsiegel}, Lemma 2.4.

Our cycles will be defined by imposing certain collections of special endomorphisms.

Fix $n$, with $1\le n\le 4$. Let $\o\subset V(\A_f)^n$ be a compact
open subset, stable under the action of $K\subset G(\A_f)$, and let
$T\in Sym_n(\Q)$ be a symmetric matrix with $\det(T)\ne 0$. {\it
We will soon assume that $T$ is positive definite and $n\le 2$.}
The cases $n=3$ or $4$ will be relevant in positive characteristic.

For the additional data $T$, $\o$, we consider the functor which
associates to each $S\in\text{Sch}/\Q$ the set of isomorphism
classes of $5$-tuples $(A,\lambda,\iota,\bar{\eta};\j)$, where
$(A,\lambda,\iota,\bar{\eta})$ is as before, and where $\j$ is an
$n$-tuple of special endomorphisms $j_1,\dots,j_n\in \End^0(A)$.
The set of endomorphisms $\eta\circ \j\circ \eta^{-1}$
corresponding to $\j$ under an isomorphism
$\eta:\hat{V}(A)\rightarrow U(\A_f)$ for $\eta\in \bar{\eta}$, are
special endomorphisms of $U_{{\Bbb A}_f}$ and hence are given as
right multiplication by elements of $V(\A_f)$. We require that
$$\eta\circ \j\circ \eta^{-1} \in \o\subset V(\A_f)^n, \tag 1.4$$
and that, in addition,
$$Q(\j)= T.\tag 1.5$$

Here, as in the rest of the paper, for $\x\in V^n$, we put
$$Q(\x)= Q((x_1, \ldots, x_n))= {\scriptstyle \frac{1}{2}} ((x_i,
x_j))\in \Sym_n(\Q)\tag1.6$$
where $(\ ,\ )$ is the bilinear form on $V$
with
$$(x,y)=Q(x+y)-Q(x)-Q(y),\tag1.7$$
so that $Q(x)=\frac12(x,x)$.
Since $\o$ is a union of $K$-orbits and since the action of $G$ on
$V$ preserves the quadratic form, conditions (1.4) and (1.5) are
independent of the choice of $\eta$ in the $K$-equivalence class
$\bar{\eta}$.

\proclaim{Proposition 1.3} The functor thus determined by
$T$ and $\o$ has a coarse moduli scheme $Z(T,\o)$. If $K$ is
sufficiently small, then $Z(T,\o)$ is a fine moduli scheme and the
natural forgetful morphism $Z(T,\omega)\to M_K$ is finite and
unramified.
\qed
\endproclaim

Of course, if we impose too many special endomorphisms, the resulting space
could be empty. The next result describes what happens in characteristic $0$.

\proclaim{Proposition 1.4} (i) $Z(T,\o)\ne \emptyset$ implies that $n=1$ or $2$ and that
$T\in Sym_n(\Q)$ is positive definite and is represented by $V(\Q)$. \hfill\break
(ii) If $T$ satisfies these conditions, fix $x\in V(\Q)^n$
such that $Q(x) =T$. Let
$G_x$ be the stabilizer of $x$ in $G$, and let
$${\Cal D}_x=\{ z\in {\Cal D};\  z\perp x\ \}.$$
For $h\in G(\A_f)$, let $Z(x,h;K)$ be the image of the map
$$ G_x(\Q)\back {\Cal D}_x\times G_x(\A_f)/(G_x(\A_f)\cap h K h^{-1})
\lra G(\Q)\back {\Cal D}\times G(\A_f)/K,$$
sending $(z,g)$ to $(z,gh)$. Also, let $\ph_\o\in S(V(\A_f)^n)$ be
the characteristic function of the compact open set $\o$. Then
$${\roman{image}}(Z(T,\o)(\C)) = \bigcup\limits_{\{ r; \varphi_{\omega}
(h_r^{-1}x)= 1\}}Z(x,h_r;K),$$ where $h_r$ runs over a set of
representatives for the double cosets $G_x(\A_f)\back G(\A_f)/K$.
\endproclaim

Under the map $h\mapsto h^{-1}x$, $G_x(\A_f)\back G(\A_f)$ is
homeomorphic to a closed subset of $V(\A_f)^n$. The intersection of
this set with the compact open subset $\o$ is thus compact, and
hence is a finite union of $K$-orbits. Thus the union on $r$ is
indeed finite.

\demo{Proof}
Suppose that $\xi=(A,\lambda,\iota,\bar{\eta};\j)\in Z(T,\o)(\C)$.
Then, $H_1(A,\Q)$ is a $16$-dimensional vector space over $\Q$,
which, by the determinant condition, is isomorphic to $U_\Q$ as a
module over $C(V)\tt \kay$. Fix such an isomorphism
$\mu:H_1(A,\Q)\isoarrow U_\Q$, and note that the complex structure
on $H_1(A,\R)\simeq U_\R$ is given by right multiplication by an
element $j_z$ for some $z\in {\Cal D}$. Then
$\mu\circ\j\circ\mu^{-1}$ is an $n$--tuple of special endomorphisms of
$U_\Q$, i.e.,
$$y=\mu\circ\j\circ\mu^{-1}\in V(\Q)^n,$$
and, by condition (1.5),
$$Q(y)=Q(\mu\circ\j\circ\mu^{-1})=T.$$
The components of $y$ must commute with $j_z$, and it is not
difficult to check that an element $y\in V\subset C(V)$ commutes
with $j_z$ if and only if $y\in z^\perp$, the orthogonal complement
of $z$ in $V_\R$, comp.\ \cite{\krsiegel}, Lemma 2.3. Note that
$z^\perp$ is a positive $2$-plane in $V_\R$. Therefore
$$y=\mu\circ\j\circ\mu^{-1}\in (\,z^\perp\,)^n,$$
and, since $\det(T)\ne 0$, we must have $n=1$ or $2$, and $T$
must be positive definite. This proves (i).

Replacing $\mu$ by $r(\gamma)\circ \mu$ for some $\gamma\in G(\Q)$,
changes $y$ to $\gamma y \gamma^{-1}$. By Witt's Theorem, we may
assume that $y=\mu\circ\j\circ\mu^{-1}=x$, our fixed  $n$-tuple of
endomorphisms, and that the complex structure is given by $r(j_z)$,
where $z\in {\Cal D}_x$. Note that $\mu$ can still be changed by an
element of $G_x(\Q)$. The composition $\eta\circ\mu^{-1}\in
\End(U(\A_f),i)$ is given by right multiplication by an element
$h\in G(\A_f)$, and condition (1.4) becomes
$$\align
\eta\circ \j\circ\eta^{-1} &= \eta\circ \mu^{-1}\circ x\circ\mu \circ\eta^{-1}\\
{}&= r(h)\circ r(x)\circ r(h^{-1})\\ {}&= r(h^{-1}xh) \in
r(\o),\endalign$$ i.e., $h^{-1}\cdot x\in\o$. Write $h = g h_r k$,
with $g\in G_x(\A_f)$ and $k\in K$ for some double coset
representative. This then gives the desired description of
${\roman{image}}(Z(T,\o)(\C))$.
\qed\enddemo

{\bf Remarks:} (i) For $n=1$ and $2$, the group theoretic cycles
$Z(x,h;K)$ and their sums over double coset representatives were
introduced in \cite{\duke}. Proposition~1.4 gives them a modular
interpretation in the case of signature $(2,2)$.
\hfill\break
(ii) For $n=1$, the cycle $Z(T,\o)$ is a divisor, rational over $\Q$, on the
surface $M$. When $B=M_2(\kay)$ for a real quadratic field $\kay$, for a suitable
choice of $T\in \Q_{>0}$ and $\o\subset V(\A_f)$, we recover the Hirzebruch-Zagier
curves $T_N$ on the Hilbert modular surface. \hfill\break
(iii) For $n=2$, the $Z(T,\o)$'s are $\Q$--rational $0$--cycles on $M$.

\subheading{\Sec2. Models over $\Z_{(p)}$ and arithmetic special cycles}

In this section, we define models of the quaternionic Hilbert-Blumenthal
surfaces $M_K$ and of the special cycles over $\Z_{(p)}$, in the case in
which $p$ is a prime of good reduction.

Fix a prime $p\ne 2$, and assume that there is a $\Z_{(p)}$-lattice
$\Lambda\subset V(\Q)$ which is self dual with respect to (1.7),
i.e.,
$$\Lambda =\{ x\in V(\Q);  (x,\Lambda)\subset \Z_{(p)}\ \}.\tag2.1$$
It follows that $p$ is unramified in $\kay$ and that $B$ is split
at each prime $\wp$ of $\kay$ over $p$.

Let $\OC=C(\Lambda)$ be the Clifford algebra of $\Lambda$. Then
$\OC$ is a $\Z_{(p)}$-order in $C(V)$ which is maximal at $p$ and
which is invariant under the main involution $\iota$. Let $\OK$ be
the ring of $\Z_{(p)}$-integers in $\kay$, and note that
$\OK\subset \OC$. We also assume that $\Lambda=\tau \Lambda
\tau^{-1},$ where $\tau\in B^\times$ is the element used to define
the positive involution $*$. It follows that $\OC$ is invariant
under $*$. Let $U =\OC\subset U_\Q$. This $\OC\tt\OK$-submodule of
$U_\Q$ is self dual with respect to $<\ ,\ >$.

Having chosen this additional data with respect to $p$, and a
compact open subgroup $K^p\subset G(\A_f^p)$, we can define a
moduli problem over $\Z_{(p)}$ which associates to
$S\in\text{Sch}/\Z_{(p)}$ the set of isomorphism classes of
4-tuples $(A,\lambda,\iota,  \overline{\eta^p})$ where
\roster
\item"{(i)}" $A$ is an abelian scheme over $S$, up to prime to $p$ isogeny,
\item"{(ii)}" $\lambda:A\isoarrow \hat{A}$, is a $\Z_{(p)}^\times$--class of
principal polarizations on $A$,
\item"{(iii)}" $\iota:\OC\tt\OK \lra \End(A)\tt\Z_{(p)}$ is a homomorphism such that
$$\iota(c\tt a)^* = \iota(c^*\tt a),$$
for the Rosati involution of $\End^0(A)$ determined by $\lambda$
and the involution $*$ of $C(V)$, and
\item"{(iv)}" $  \overline{\eta^p}$ is a $K^p$-equivalence class of
$\OC\tt\OK$-equivariant isomorphisms
$$\eta^p:\hat{V}^p(A) \isoarrow U(\A^p_f)$$
which preserve the symplectic forms up to a scalar in
$(\A^p_f)^\times$.
\endroster
Again, the action of $\iota(c\tt a)$ on $\text{Lie}(A)$ is
required to satisfy the determinant condition
$$\det(\iota(c\tt a);\text{Lie}(A)) = N^0(c)^2 N_{\smallkay/\Q}(a)^4,$$
interpreted as an identity of polynomial functions with
coefficients in $\Cal O_S$, as in \cite{\kottwitz}.

As explained in section~5 of \cite{\kottwitz}, for $K^p\subset
G(\A_f)$ sufficiently small, this moduli problem is represented by
a smooth quasi-projective scheme $\Cal M_{K^p}$ over $\Z_{(p)}$,
which is, in fact, projective if $B=C^+(V)= \End(U_\Q,i)$ is a
division algebra. Moreover, if $K_p$ denotes the intersection of $(\OC\tt\Z_p)^\times$ with
$G(\Q_p)$, and $K=K_pK^p$, then
$$M_K\simeq {\Cal M}_{K^p}\times_{\roman{Spec}
\, \Z_{(p)}}
{\roman{Spec}}\, \Q\ \ ,$$ so that $\Cal M_{K^p}$ gives a model of
$M_K$ over $\Z_{(p)}$.

We next turn to the special cycles. Again, for a fixed $n$ with
$1\le n\le 4$, we let $\o^p\subset V(\A_f^p)^n$ be a compact open
subset, stable under the action of $K^p$, and let $T\in Sym_n(\Q)$,
with $\det(T)\ne 0$. We point out that for
$(A,\lambda,\iota,\overline{\eta^p})\in {\Cal M}_{K^p}(S)$ we may
introduce the concept of a {\it special endomorphism}\/ in analogy
with Definition 1.2. If $S$ is connected the special endomorphisms
form a positive-definite quadratic space over $\Z_{(p)}$.

We consider the functor which associates to each scheme
$S\in\text{Sch}/\Z_{(p)}$, the set of isomorphism classes of
5-tuples $(A,\lambda,\iota,  \overline{\eta^p};\j)$, with
$(A,\lambda,\iota,  \overline{\eta^p})$ as before, and with $\j\in
\big(\End(A)\tt\Z_{(p)}\big)^n$ an $n$-tuple of special endomorphisms
satisfying (1.5), and the analogue of (1.4) with $\eta^p$ in place
of $\eta$ and $\o^p$ in place of $\o$.

\proclaim{Proposition 2.1} (i) If $K^p$ is sufficiently small, then the moduli
problem just defined is representable by a scheme $\Cal Z(T,\o^p)$
which maps by a finite unramified morphism to $\Cal
M_{K^p}$.\hfill\break (ii) Let $\o_p = (\Lambda\tt\Z_p)^n
\subset V(\Q_p)^n$. Then the generic fiber of $\Cal Z(T,\o^p)$ is
$$Z(T,\o_p\times\o^p) = \Cal Z(T,\o^p)\times_{\Spec\, \Z_{(p)}}
\Spec\, \Q .$$
If $n\ge 3$, then the left side is to be interpreted as the empty scheme.
\qed
\endproclaim

The reader should be warned, however, that, in contrast to the
situation for $\Cal M_{K^p}$, $\Cal Z(T,\o^p)$ is not flat over
$\Z_{(p)}$ in general. Indeed, as already in \cite{\krsiegel}, one of
the most interesting cases arises when $n\ge3$, so that the generic
fiber of $\Cal Z(T,\o^p)$ is empty.

We add here a remark on Hecke correspondences prime to $p$. Let
$g\in G({\Bbb A}_f^p)$. The Hecke correspondence on ${\Cal
M}_{K^p}$ is defined by the diagram
$$\matrix
{\Cal M}_{K^p\cap g^{-1}K^pg} &
\buildrel g\over\longrightarrow
& {\Cal M}_{gK^pg^{-1}\cap K^p}\\
\nass
\big\downarrow
&&
\big\downarrow\\
\nass
{\Cal M}_{K^p} &
--\rightarrow
& {\Cal M}_{K^p}
\endmatrix\tag2.2$$
Here the solid horizontal morphism is defined by mapping
$(A,\lambda,\iota,\overline{\eta^p})$ to
$(A,\lambda,\iota,\overline{g\circ\eta^p})$ and the vertical
arrows are defined by the natural inclusions of open compact
subgroups. For $T$ and $\o^p$ as before, we obtain a
morphism
$$\CZ(T,\o^p)\buildrel g\over\longrightarrow
\CZ(T,g\cdot\o^p)\tag2.3$$
which sends $(A,\lambda,\iota,\overline{\eta^p},\j)$ to
$(A,\lambda,\iota,\overline{g\circ\eta^p};\j)$. It is obvious that this
morphism may be combined with the natural forgetful morphisms from
the special cycles to the moduli spaces to form a commutative
diagram with three cartesian squares. We may consider this as a
prolongation of the Hecke correspondence from ${\Cal M}_{K^p}$ to
itself to a correspondence from $\CZ(T,\o^p)$ to
$\CZ(T,g\cdot\o^p)$.

From now until section 11 we will assume that $p$ is inert in
$\kay$. This implies that $\kay$ is a field (real quadratic since
the signature of $V$ is $(2,2)$).

From now on we will denote by $\o\subset V({\Bbb A}_f^p)^n$ the set
which is denoted by $\o^p$ above. We also fix an algebraic closure
${\Bbb F}$ of the residue field of $\kay$ at $p$.

\subheading{\Sec3. Isogeny classes and special endomorphisms}

In this section, we determine the $\Q$-vector
space of special endomorphisms attached to each of the isogeny classes in the
special fiber of $\Cal M$. This information determines which isogeny classes can meet the
various special cycles.
We continue to assume that $\kay$, the center of $C^+(V)$, is a quadratic field and
that $p$ is inert in $\kay$.

For a point
$\xi=(A,\lambda,\iota,\bar{\eta}^p)\in \Cal M_{K^p}(\F)= \Cal M_p$, we
consider the isogeny class of the triple $(A,\lambda,\iota)$. Thus,
$A$ is an abelian variety of dimension $8$ over $\F$ with an action
$\iota:C(V)\tt\kay\hookrightarrow \End^0(A)$.

Recall that the Clifford algebra $C(V)=C^+(V)\oplus C^-(V)$ is a central simple algebra over $\Q$,
with $C^+(V)=B$, a quaternion algebra over $\kay$.
Moreover, there is an isomorphism
$$i:C(V)\tt_\Q\kay\simeq M_2(B),\tag 3.1$$
discussed in more detail below. Hence,  up to isogeny, $A\simeq A_0^2$ where $A_0$ is a $4$--dimensional
abelian variety with $\iota_0:B\hookrightarrow \End^0(A_0)$.
Thus, the possible isogeny classes of $(A,\iota)$'s are the same as the
possible  isogeny classes of $(A_0,\iota_0)$'s, and these are
described by Theorem~5.4 of Milne, \cite{\milne}:
\roster
\item"{$I_0.$}" The supersingular class; $A_0\simeq A_{00}^4$ for a supersingular elliptic curve $A_{00}$.
\item"{$I_1.$}" The classes associated to imaginary quadratic extensions $E/\Q$; $A_0\simeq A_{00}^4$
for an ordinary elliptic curve $A_{00}$ with complex multiplication by $E$, i.e., $E\simeq \End^0(A_{00})$.
\item"{$I_2.$}" The classes associated to totally imaginary quadratic extensions $\kay'$ of $\kay$
which are not of the form $\kay\cdot E$ for any imaginary quadratic extension of $\Q$;
$A_0\simeq A_1^2$, where $A_1$ is a simple abelian surface with complex multiplication by $\kay'$, i.e.,
$\kay'\simeq \End^0(A_1)$.
\endroster
In case $I_1$, let $\kay'=\kay\cdot E$. Then, in cases $I_1$ and $I_2$, the field $\kay'$ splits $B$,
and $p$ splits in $\kay'/\kay$.

Note that type $I_0$ consists of supersingular abelian varieties while types $I_1$ and $I_2$ consist of
ordinary abelian varieties.

For a point
$\xi\in\Cal M_p$, as above, the space of special
endomorphisms of $(A,\lambda,\iota)$ is the $\Q$--vector space
$$V'= V_\xi=\{ j \in \End^0(A)\mud \iota(c\tt a)\cdot j = j\cdot \iota(c\tt a^\s),
\text{ and } j^*=j \}.\tag 3.2$$
This space depends only on the isogeny class of $(A,\lambda,\iota)$, and
carries a quadratic form defined by $j^2= Q'(j)\cdot id_A$.
In order to determine the quadratic space $(V',Q')$ in each of the cases $I_0$, $I_1$ and $I_2$,
we must first make the isomorphism (3.1) more explicit.
In particular, we must keep track of the involution $*$ and the Galois action $\s$.

As in section 1, choose $v_0\in V$, with $Q(v_0) = v_0^2=\a\in \Q^\times$, and assume that
$\a>0$.
Write
$$C(V)= C = B\oplus  B v_0.\tag 3.3$$
As in section 1, let $\s=Ad(v_0)$, and recall that $\s$ preserves $B$ and extends
the Galois automorphism $\s$ on its center $\kay$.
The fixed algebra $B_0$ of $\s$ in $B$ is a quaternion
algebra over $\Q$ with $B=B_0\tt_\Q\kay$.
As in Lemma~0.5, let $\tau\in B_0^\times$  be the element
defining the positive involution $b^* = \tau b^\iota \tau^{-1}$ of $C(V)$.
Since $v_0^\iota=v_0$, we have $v_0^*=v_0$.

\proclaim{Lemma 3.1} For a choice of the decomposition (3.3),
the isomorphism (3.1) is given explicitly by
$$i:(b_1+b_2 v_0)\tt a \mapsto \pmatrix b_1&\a b_2\\b_2^\s&b_1^\s\endpmatrix\cdot a,$$
for $b_1$ and $b_2\in B$ and $a\in \kay$.
Under this isomorphism, the involution $c\tt a\mapsto c^\iota \tt a$
becomes
$$\beta\mapsto Ad(\pmatrix \a&{}\\{}&1\endpmatrix)\big({}^t\beta^\iota\big),$$
for $\beta\in M_2(B)$.
The involution
$c\tt a\mapsto c^*\tt a$ becomes
$$\beta\mapsto Ad(\pmatrix \tau\a&{}\\{}&\tau\endpmatrix)\big({}^t\beta^\iota\big),$$
and the automorphism $c\tt a\mapsto c\tt a^\s$ becomes
$$\beta\mapsto Ad(\pmatrix {}&\a\\1&{}\endpmatrix)\big(\beta^\s\big).$$
\endproclaim
\demo{Proof} The space $U_\Q=C(V)$ is a two dimensional
right vector space over $B$, and left-right multiplication gives an algebra homomorphism
$$C\tt_\Q\kay \lra \End_B(U_\Q),\qquad c\tt a\mapsto\big(x\mapsto cxa\big).\tag 3.4$$
This is an isomorphism of simple algebras over $\Q$.
Taking $1$ and $v_0$ as a $B$-basis for $U_\Q$, we obtain the
isomorphism $C\tt\kay\rightarrow M_2(B)$ of the Lemma. For example,
$$b(x+v_0y)=b(1,v_0)\pmatrix x\\y\endpmatrix =
(1,v_0)\pmatrix b&{}\\{}&b^\s\endpmatrix\pmatrix x\\y\endpmatrix.\tag 3.5$$

Now, for example,
$$\align
i((b_1+b_2v_0)^\iota) &= i(b_1^\iota+v_0 b_2^\iota) = i(b_1^\iota+b_2^{\iota\s}v_0)\\
\nass
{}&=\pmatrix b_1^\iota&\a b_2^{\iota\s}\\b_2^\iota&b_1^{\iota\s}\endpmatrix\tag 3.6\\
\nass
{}&=Ad\big(\pmatrix \a&{}\\{}&1\endpmatrix\big) \pmatrix b_1^\iota&b_2^{\iota\s}\\\a b_2^\iota&b_1^{\iota\s}\endpmatrix\\
\nass
{}&=Ad\big(\pmatrix \a&{}\\{}&1\endpmatrix\big) {}^t i(b_1+b_2v_0)^\iota,
\endalign$$
as claimed. The other relations are similarly easily checked.
\qed
\enddemo

Note that the idempotents
$$e_1 =\frac12 \big( 1\tt1+\delta\tt\delta^{-1}\big)\qquad\text{and}\qquad e_2 =\frac12 \big( 1\tt1-\delta\tt\delta^{-1}\big)\tag 3.7$$
in the subalgebra $\kay\tt\kay\subset C\tt_\Q\kay$ have images
$$i(e_1)=\pmatrix 1&{}\\{}&0\endpmatrix\qquad\text{and}\qquad i(e_2)= \pmatrix 0&{}\\{}&1\endpmatrix\tag 3.8$$
in $M_2(B)$, and that the decomposition $A\simeq A_0^2$
is then given by
$$A = i(e_1)A\times i(e_2)A.\tag 3.9$$
Let $\lambda_0$ be the polarization of $A_0$ determined by
$$\lambda_0 = {}^t i(e_1)\circ\lambda\circ i(e_1),\tag 3.10$$
and let $\beta\mapsto \beta'$ be the corresponding Roasti involution of
$\End^0(A_0)$.

\proclaim{Proposition 3.2} The isomorphism $A\simeq A_0^2$, determined by a
decomposition (3.3), induces an isomorphism:
$$
V'=\left\{ \matrix j \in \End^0(A)\\
\nass
\iota(c\tt a)j  = j \iota(c\tt a^\s)\\
\nass
j ^*=j \endmatrix\right\}
\simeq
\left\{ \matrix j_0\in \End^0(A_0)\\
\nass
\iota_0(b)j_0 = j_0\iota(b^\s)\\
\nass
j_0'=j_0\endmatrix\right\}.
$$
on the space of special endomorphisms. Here $j^2 = \a j_0^2\cdot 1_2$,
so that $Q(j) = \a\, Q_0(j_0)$ with $j_0^2=Q_0(j_0)\,id_{A_0}$.
\endproclaim

\demo{Proof}
Having fixed the isomorphism (3.1), with $A = i(e_1)A\times i(e_2)A$, we have
$$\iota:C\tt\kay = M_2(B)\hookrightarrow M_2(\End^0(A_0)) \simeq \End^0(A),\tag 3.11$$
where the middle inclusion is induced by  $\iota_0:B\hookrightarrow \End^0(A_0)$.
Recall that this map is required to satisfy $\iota(c\tt a)^*=\iota(c^*\tt a)$.
If $j$ is
a special endomorphism, then, for $a\in\kay$,
$j\circ\iota(1\tt a) = \iota(1\tt a^\s)\circ j$,
and $j\circ\iota(a\tt1) = \iota(a\tt1)\circ j$. In matrix form, these conditions
imply that
$$j = \pmatrix {}&j_2\\j_3&{}\endpmatrix,\tag 3.12$$
with $j_2$ and $j_3\in\End^0(A_0)$.
The conditions $j\circ\iota(b\tt 1)=\iota(b\tt1)\circ j$
and $j\circ\iota(v_0\tt1)=\iota(v_0\tt1)\circ j$ imply
$j_2=\a j_3$ and $j_3 \iota_0(b) =\iota_0(b^\s)j_3$.

Finally, we must impose the condition $j^*=j$.
The polarization has the form:
$$\pmatrix\lambda_1&\lambda_2\\\lambda_3&\lambda_4\endpmatrix=\lambda :
A_0^2\simeq A\lra \hat{A} \simeq \hat{A_0}^2,\tag 3.13$$
where, for example, $\lambda_2= {}^t\iota(e_1)\lambda\iota(e_2)$, etc.
But now, since the idempotents $e_1$ and $e_2$ which
give the decomposition $A\simeq A_0^2$ lie in $\kay\tt\kay\subset C\tt\kay$,
and since the Rosati involution restricts to the identity on $\iota(\kay\tt\kay)$,
we must have $\lambda_2=\lambda_3=0$.  Since $v_0^*=v_0$, we have
${}^t\iota(v_0)\lambda =\lambda \iota(v_0)$, and hence
$$\pmatrix {}&1\\\a&{}\endpmatrix
\pmatrix \lambda_1&{}\\{}&\lambda_4\endpmatrix =
\pmatrix \lambda_1&{}\\{}&\lambda_4\endpmatrix \pmatrix {}&\a\\1&{}\endpmatrix,\tag 3.14$$
so that $\lambda_4=\a\lambda_1$.
Therefore, the Rosati involution on $M_2(\End^0(A_0))$ is given by
$$\beta \mapsto Ad(\pmatrix 1&{}\\{}&\a\endpmatrix^{-1})\big({}^t\beta'\big),\tag 3.15$$
where $x\mapsto x'$ is the positive involution of $\End^0(A_0)$ determined by
$\lambda_1$, applied
componentwise to the matrix $\beta$.
But then the condition
$$\pmatrix {}&\a j_3\\j_3&{}\endpmatrix = j
=j^* = \pmatrix {}&j_3'\\\a j_3'&{}\endpmatrix,\tag 3.16$$
is equivalent to $j_3'=j_3$.
\qed\enddemo

We are now in a position to determine the space of special endomorphisms for
each isogeny class in the special fiber. Proposition~3.2 shows that
this space can be determined in terms of the triple $(A_0,\iota_0,{}')$.
Recall that, for $b\in B$,
$$\iota_0(b)' = \iota_0(\tau b^\iota \tau^{-1}),\tag 3.17$$
and that
$\kay=\Q(\delta)$ with $\delta^2=\Delta\in\Q^\times$.

We begin with the supersingular case $I_0$, so that
$A_0\simeq A_{00}^4$ for a supersingular elliptic curve $A_{00}$.
Let $\Bbb B$ be the quaternion
algebra over $\Q$ ramified precisely at infinity and $p$, and for $B=B_0\tt_\Q\kay$, as above,
let $B_0^\vee$ be the unique (up to isomorphism) quaternion algebra over $\Q$ such
that $B_0\tt_\Q B_0^\vee \simeq M_2(\Bbb B)$.
Then we have
$$B_0\tt_\Q B_0^\vee\tt_\Q M_2(\Q) \isoarrow M_4(\Bbb B) = \End^0(A_0).\tag 3.18$$
We may assume that the embedding $\iota_0$ of $B=B_0\tt_\Q\kay$ is given by
$$\iota_0(b_0\tt a) = b_0\tt 1\tt i_0(a),\tag 3.19$$
where $i_0:\kay\hookrightarrow M_2(\Q)$ is determined by
$$i_0(\delta) = \pmatrix {}&1\\\Delta&{}\endpmatrix.\tag 3.20$$

By condition (3.17), we know the restriction of the Rosati
involution $'$ of $\End^0(A_0)$ to $\iota_0(B)$. One positive involution of
$B_0\tt B_0^\vee\tt M_2(\Q)$ which satisfies (3.17) is given by
$$\beta = b_0\tt b_1\tt c \mapsto \beta^\sh :=
\tau b_0^\iota \tau^{-1}\tt b_1^\iota \tt Ad\big(\pmatrix 1&{}\\{}&\Delta\endpmatrix\big) {}^tc.\tag 3.21$$
\proclaim{Lemma 3.3} Any positive involution $'$ of $\End^0(A_0)$ satisfying (3.17) has the form
$$\beta\mapsto \beta'= \eta \,\beta^\sh \eta^{-1},$$
where $\eta = 1\tt1\tt i_0(\eta_0)$ for $\eta_0\in \kay^\times$ with
$N_{\smallkay/\Q}(\eta_0)>0$.
\endproclaim
\demo{Proof} Any involution of $B_0\tt B_0^\vee\tt M_2(\Q)$ has the form
$\beta\mapsto \beta'=\eta \beta^\sh \eta^{-1}$ for some invertible element $\eta$. The fact that
$'$ is an involution is equivalent to the condition that $\eta^\sh \eta^{-1}=\e\in\Q^\times$.
Applying $\sh$ to this condition, we find that $\e^2=1$. So there are two basic cases
depending on $\e=\pm1$. On the other hand, since $'$ also satisfies (3.17), $\eta$
must lie in the centralizer of $B_0\tt 1\tt i_0(\kay)$, i.e., in $1\tt B_0^\vee\tt i_0(\kay)$.
We may identify $B_0^\vee\tt M_2(\Q)$ with $M_2(B_0^\vee)$ and write
$$\eta = \pmatrix \eta_1&\eta_2\\\Delta \eta_2&\eta_1\endpmatrix,\tag 3.22$$
with $\eta_1$ and $\eta_2\in B_0^\vee$. Applying $\sh$ yields
$$\eta^\sh = \pmatrix \eta_1^\iota&\eta_2^\iota\\\Delta \eta_2^\iota&\eta_1^\iota\endpmatrix = \e
\pmatrix \eta_1&\eta_2\\\Delta \eta_2&\eta_1\endpmatrix.\tag 3.23$$

If $\e=-1$, we have
$\eta_1^\iota =-\eta_1$ and $\eta_2^\iota=-\eta_2$.
\proclaim{Lemma 3.4} If $\eta^\sh =-\eta$, then the involution
$\beta\mapsto \eta \beta^\sh \eta^{-1}$ is not positive.
\endproclaim
\demo{Proof} We may as well replace $B_0^\vee$ with $B_0^\vee\tt\R=\Bbb H$.
If $\kappa\in \Bbb H$ with $\kappa^\iota=-\kappa$, then the involution
$z\mapsto \kappa z^\iota \kappa^{-1}$ is not positive. This follows from the
fact that there is an
element $z\in \Bbb H^\times$ such that $\kappa z = -z \kappa$, so that
$$\tr(z \kappa z^\iota \kappa^{-1}) = -\tr(\kappa z z^\iota \kappa^{-1}) = -2 z z^\iota <0.\tag 3.24$$

For convenience, let
$$\teta = \eta \pmatrix 1&{}\\{}&\Delta\endpmatrix = \pmatrix
\eta_1&\Delta\eta_2\\\Delta\eta_2&\Delta\eta_1\endpmatrix,\tag 3.25$$
so that, for $x\in M_2(B_0^\vee)$ or $M_2(\Bbb H)$,
$$x' = \teta \left({}^tx^\iota \right)\teta^{-1}.\tag 3.26$$
We want to consider the quantity $\tr(x x')$.
If $g\in GL_2(\Bbb H)$, then for $x\in M_2(\Bbb H)$,
$$\tr\big(g^{-1}xg \teta \left({}^t(g^{-1}x g)^\iota\right) \teta^{-1}\big)
=\tr(x (g\teta\, {}^t\!g^\iota) \left({}^t x^\iota\right) (g\teta\, {}^t\!g^\iota)^{-1}).\tag 3.27$$
Thus the positivity of the involution associated to $\teta$ is equivalent to that
of the involution associated to $g\teta\, {}^t\!g^\iota$.
But now, if $\eta_1\ne 0$,
$$\pmatrix 1&{}\\-\eta_2\Delta\eta_1^{-1}&1\endpmatrix \pmatrix
\eta_1&\Delta\eta_2\\\Delta\eta_2&\Delta\eta_1\endpmatrix
\pmatrix 1&-\eta_1^{-1}\Delta \eta_2\\{}&1\endpmatrix
=\pmatrix \eta_1&{}\\{}&\Delta(\eta_1-\Delta \eta_2\eta_1^{-1}\eta_2)\endpmatrix.\tag 3.28$$
Thus, for $x=\pmatrix x_1&0\\0&0\endpmatrix$,
we have
$$\tr(x (g\teta\, {}^t\!g^\iota) {}^t x^\iota (g\teta\,
{}^t\!g^\iota)^{-1}) = \tr(x_1\eta_1 x_1^\iota \eta_1^{-1}),\tag 3.29$$
and this quantity can be negative, as we observed above. If $\eta_1=0$, an easier
calculation yields the same conclusion,
and Lemma~3.4 is proved.
\qed\enddemo

Continuing the proof of Lemma~3.3, we must have $\e=+1$, so that $\eta_1$ and $\eta_2\in \Q$, and hence,
$\eta=i_0(\eta_0)$ for $\eta_0=\eta_1+\delta\eta_2\in \kay^\times$.
Applying (3.28) again, and taking $x=\pmatrix x_1&x_2\\x_3&x_4\endpmatrix$, we find that
$$\tr(x (g\teta\, {}^t\!g^\iota) {}^t x^\iota (g\teta\, {}^t\!g^\iota)^{-1}) =
2 x_1x_1^\iota + 2\Delta \eta_1^{-2} N(\eta_0) x_2x_2^\iota +2\Delta^{-1}\eta_1^2 N(\eta_0)^{-1}x_3x_3^\iota
+2 x_4x_4^\iota.$$
This is positive provided $x\ne 0$ and $N(\eta_0)>0$.
\qed\qed\enddemo

\proclaim{Proposition 3.5}
Suppose that $\xi=(A,\lambda,\iota,\overline{\eta^p})\in \Cal M_p^{\text{ss}}$ is
a point in the supersingular isogeny class $I_0$, and that the Rosati involution
$\lambda_0$ of $A_0$ is given as in Lemma~3.3 for $\eta_0\in \kay^\times$. Then
the space of special endomorphisms (3.2) is a
$4$--dimensional $\Q$--vector space:
$$V'\simeq
\{\pmatrix a&b\\-\Delta b&-a\endpmatrix\mud a\in \Q,\ b\in B_0^\vee,\ b^\iota=-b\ \},$$
where $B_0^\vee$ is the centralizer of $B_0$ in $M_2(\Bbb B)$, as above. The
quadratic form on this space is
$$Q_\xi(j) = \a N_{\smallkay/\Q}(\eta_0)^{-1}( a^2 + \Delta N(b) ).$$
Since $\a$ and $N(\eta_0)$ are positive, this form is positive definite.
\endproclaim
\demo{Proof} For a special endomorphism $j_0$ of $A_0$, we have
$j_0\iota_0(b)=\iota_0(b^\s)j_0$, and hence
$$j_0\in 1\tt B_0^\vee\tt \pmatrix 1&{}\\{}&-1\endpmatrix \kay.\tag3.30$$
As an element of $M_2(B_0^\vee)$, we can write
$$j_0= \pmatrix 1&{}\\{}&-1\endpmatrix \pmatrix b_1&b_2\\\Delta b_2&b_1\endpmatrix.\tag3.31$$
Then, recalling that $\eta^\sh=\eta$, the condition
$$j_0'=\eta j_0^\sh \eta^{-1} =j_0\tag3.32$$
is equivalent to
$$(j_0\eta)^\sh =j_0\eta.\tag3.33$$
Writing
$$j_0\eta = \pmatrix c_0&c_1\\-\Delta c_1&-c_0\endpmatrix,\tag3.34$$
we find that $c_0=c_0^\iota\in \Q$ and $c_1^\iota = -c_1$. This gives the
claimed description of $V'$.
Also
$$(j_0\eta)^2 = j_0^2 \eta^\s\eta,\tag3.35$$
so that
$$j_0^2 = N(\eta_0)^{-1}(c_0^2+\Delta N(c_1)),\tag3.36$$
and, recalling Proposition~3.2, we obtain the claim expression
for the quadratic form on $V'$.
\qed\enddemo

In the case $\text{I}_1$, we have $A_0\simeq A_{00}^4$ and
$\iota_0:B\hookrightarrow \End^0(A_0)\simeq M_4(E)$, where $\End^0(A_{00})=E$
and $\kay'=\kay\cdot E$ splits $B$.

\proclaim{Proposition 3.6} (i) Let $E/\Q$ be an imaginary
quadratic field, and let $E^\vee$ be the other imaginary quadratic
subfield of the biquadratic field $\kay\cdot E$. For the isogeny class
$\text{I}_1$ associated to
$E/\Q$, the space of special endomorphisms is isomorphic, as quadratic space over
$\Q$,
to $E^\vee$ with quadratic form $Q_\xi(j) = \a\a_1 N_{E^\vee/\Q}(j)$,
where the additional scalar $\a_1\in \Q^\times_{>0}$ is determined in the proof below.
\hfill\break
(ii) Let $\kay'$ be a CM extension of $\kay$ which is not a composite of $\kay$
with an imaginary quadratic field. Then the space of special endomorphisms
for the associated isogeny class $\text{I}_2$ is zero.
\endproclaim

\demo{Proof} First assume that $E$ splits the quaternion algebra $B_0$ over $\Q$.
Then
$$B_0\tt_\Q E\tt_\Q M_2(\Q)\simeq M_2(E)\tt M_2(\Q)\simeq M_4(E),\tag3.37$$
so that we can take
$$\iota_0:B=B_0\tt\kay \lra B_0\tt E\tt M_2(\Q),
\qquad b_0\tt a\mapsto b_0\tt 1\tt i_0(a),\tag3.38$$
with $i_0:\kay \hookrightarrow M_2(\Q)$ as before. The positive involution
$$\beta = b_0\tt e\tt c \mapsto \beta^\sh := \tau b_0^\iota \tau^{-1}\tt
\bar{e}\tt Ad(\pmatrix 1&{}\\{}&\Delta\endpmatrix)\big({}^tc\big)\tag3.39$$
satisfies
$\iota_0(b)^\sh =\iota_0(\tau b^\iota \tau^{-1})$.
The most general such involution has the form
$$\beta\mapsto \beta'=\eta \beta^\sh \eta^{-1},\tag3.40$$
where $\eta\in 1\tt E\tt i_0(\kay)$, the centralizer of $\iota_0(B)$.
If we view $E\tt M_2(\Q) \isoarrow M_2(E)$ and $E\tt\kay \isoarrow \kay'$
with $i_0:\kay'\hookrightarrow M_2(E)$ corresponding to $1\tt i_0$, then $\eta =i_0(\eta_0)$ with
$\eta_0\in \kay^{\prime,\times}$. Since $'$ is an involution,
$\eta^\sh\eta^{-1}=\lambda\in E^\times$, and $\lambda^\sh\lambda =\bar\lambda \lambda =1$.
Writing $\lambda = \nu \bar\nu^{-1}$, we have $(\eta\nu)^\sh=\eta \nu$.
Since $\eta$ and $\eta\nu$ give the same involution, we may assume that $\eta^\sh=\eta$.
Then, since $\eta^\sh =i_0(\bar\eta_0)$, we have $\eta_0\in \kay^\times$.
It is easy to check that the involution $'$ is positive if and only if
$N_{\smallkay/\Q}(\eta_0)>0$.

If $j_0$ is a special endomorphism, then the condition
$j_0\iota_0(b)=\iota_0(b^\s)j_0$ implies that there is a $\mu_0\in \kay^{\prime,\times}$
so that,
for $\mu=i_0(\mu_0)$,
$$j_0 =1 \tt \pmatrix 1&{}\\{}&-1\endpmatrix \mu\in 1\tt
\pmatrix 1&{}\\{}&-1\endpmatrix i_0(\kay')\subset B_0\tt M_2(E).\tag3.41$$
Since
$$\pmatrix 1&{}\\{}&-1\endpmatrix' = \eta \pmatrix 1&{}\\{}&-1\endpmatrix^\sh \eta^{-1}
= \pmatrix 1&{}\\{}&-1\endpmatrix \eta^\s \eta^{-1},\tag3.42$$
we have
$$\pmatrix 1&{}\\{}&-1\endpmatrix i_0(\mu_0) = j_0 = j_0' =
\pmatrix 1&{}\\{}&-1\endpmatrix i_0(\bar\mu_0^\s) \eta^\s\eta^{-1},\tag3.43$$
so that
$$\rho:= \mu_0 \eta_0=(\overline{\mu_0\eta_0})^\s  \in E^\vee.\tag3.44$$
Hence
$$j_0= \pmatrix 1&{}\\{}&-1\endpmatrix \eta^{-1}i_0(\rho),\tag3.45$$
and
$$j_0^2 = N_{\smallkay/\Q}(\eta)^{-1} N_{E^\vee/\smallkay}(\rho) \,1_2\tag3.46$$
as claimed, in this case.

In general, there is a quaternion algebra $B_1$ over $\Q$ such that
$B= B_1\tt_\Q\kay$ and such that $E$ splits $B_1$. Let $\s_1$ be the
automorphism of $B$ which is given by $1\tt\s$ on $B=B_1\tt\kay$.
There is then an element $h\in B^\times$ such that
$b^\s = Ad(h)(b^{\s_1})$. The fact that $\s$ and $\s_1$ have order $2$ implies that
$hh^{\s_1}\in \kay^\times$. Also, $h^\s=h h^{\s_1}h^{-1}=h^{-1}h h^{\s_1}=h^{\s_1}$.
Applying Theorem 90, we may adjust $h$
so that, in fact $h^\s=h^{\s_1}=h^\iota$. In particular, $hh^\s=hh^{\s_1}=\nu(h)\in \Q^\times$.

We then have
$$B_1\tt E\tt \kay \lra B_1\tt E\tt M_2(\Q)\simeq M_2(E)\tt M_2(\Q) \simeq M_4(E),\tag3.47$$
and an embedding
$$i_1:B = B_1\tt\kay \hookrightarrow  B_1\tt E\tt M_2(\Q),\qquad b_1\tt a\mapsto b_1\tt 1\tt i_0(a).\tag3.48$$
Note that we have
$$Ad(\pmatrix 1&{}\\{}&-1\endpmatrix) i_0(a) = i_0(a^\s),\tag3.49$$
and
$$Ad(\pmatrix 1&{}\\{}&-1\endpmatrix) i_1(b) = i_1(b^{\s_1}).\tag3.50$$

Let
$$\tau_+ =\cases \tau&\text{if $\nu(h)>0$,}\\\delta\tau&\text{ if $\nu(h)<0$,}
\endcases\tag3.51
$$
and define an involution on $B_1\tt E\tt M_2(\Q)$ by
$$\beta=b_1\tt e\tt c\mapsto \beta^\sh:= Ad(i_1(\tau_+))
\big( b_1^\iota\tt \bar e\tt Ad(\pmatrix 1&{}\\{}&\Delta\endpmatrix)({}^tc)\big).\tag3.52$$
\proclaim{Lemma 3.7} (i) The involution $\sh$ is positive and satisfies
$$i_1(b)^\sh = i_1(\tau b^\iota \tau^{-1}).$$
(ii) Any positive involution with this restriction to $B$ has the form
$$\beta\mapsto \beta' =\eta \beta^\sh \eta^{-1},$$
where $\eta=i_0(\eta_0)$ for $\eta_0\in \kay^\times$ with $N_{\smallkay/\Q}(\eta_0)>0$.
\endproclaim
\demo{Proof} The involution $\sh$
preserves the subalgebra $B_1\tt M_2(\Q)$, and it suffices to prove that its restriction to
this subalgebra is positive. Since
$$M_2(\Q) =i_0(\kay)\oplus \pmatrix 1&{}\\{}&-1\endpmatrix i_0(\kay),\tag3.53$$
we have an isomorphism of right $B$-modules
$$B_1\tt M_2(\Q) \simeq B\oplus \pmatrix 1&{}\\{}&-1\endpmatrix B.\tag3.54$$
The left regular representation yields an algebra homomorphism:
$$\align
B_1\tt M_2(\Q) \simeq B \oplus B \pmatrix 1&{}\\{}&-1\endpmatrix&\hookrightarrow M_2(B)\tag3.55\\
\nass
b_1+b_2 \pmatrix 1&{}\\{}&-1\endpmatrix &\mapsto \pmatrix b_1&b_2\\b_2^{s_1}&b_1^{s_1}\endpmatrix.\endalign$$
We write $i_2:B\rightarrow M_2(B)$ for the composition of this map with $i_1$.
The restriction of the involution
$$\beta \mapsto Ad(i_2(\tau_+))\big( {}^t\beta^\iota\big)\tag3.56$$
to $B_1\tt M_2(\Q)$ coincides with $\sh$, so it suffices to prove that (3.56) is positive on $M_2(B)$.
But we have, setting $\e=\sgn(\nu(h))$,
$$\align
&\tr\left(\pmatrix a&b\\c&d\endpmatrix ^\iota\pmatrix \tau&{}\\{}&\e\tau^{s_1}\endpmatrix
\pmatrix a^\iota&c^\iota\\b^\iota&d^\iota\endpmatrix \pmatrix \tau&{}\\{}&\e\tau^{s_1}\endpmatrix^{-1}\right)\tag3.57\\
\nass
{}&\qquad=\tr(a\tau a^\iota \tau^{-1} + b\e \tau^{\s_1}b^\iota \tau^{-1} + c\tau c^\iota \e \tau^{-\s_1} + d \tau^{\s_1}d^\iota \tau^{-\s_1}).
\endalign
$$
The quantities $\tr(a\tau a^\iota \tau^{-1})$ and
$\tr(d \tau^{\s_1}d^\iota \tau^{-\s_1}) = \tr(d^{\s_1} \tau(d^\iota)^{\s_1} \tau^{-1})^\s$
in the last expression are positive for $a$ and $d$ nonzero. On the other hand,
$$\align
\tr(b\e \tau^{\s_1}b^\iota \tau^{-1}) &= \e\tr(b h^{-1} \tau h b^\iota \tau^{-1}) \tag3.58\\
\nass
{}&= \e \nu(h)^{-1} \tr( (bh^\iota) \tau (b h^\iota)^\iota \tau^{-1}) >0,
\endalign
$$
and similarly for $\tr(c\tau c^\iota \e \tau^{-\s_1})$. This proves (i).

The centralizer of $i_1(B)$ in $B_1\tt E\tt M_2(\Q)$ is $1\tt E\tt i_0(\kay)$ and this
can be identified with $1\tt i_0(\kay')$ in $B_1\tt M_2(E)$. Note that
this algebra centralizes $i_1(\tau_+)$. We can argue as above to prove (ii).
\qed\enddemo

As in (3.41), we can write any special endomorphism in the form
$$j_0 =i_1(h)\bigg(1 \tt \pmatrix 1&{}\\{}&-1\endpmatrix \mu\bigg)\in
i_1(h)\bigg(1\tt\pmatrix 1&{}\\{}&-1\endpmatrix
i_0(\kay')\bigg)\subset B_1\tt M_2(E),\tag3.59$$ for
$\mu=i_0(\mu_0)$ and $\mu_0\in \kay'$. To take into account the
condition $j_0'=j_0$, we first note that, via (3.49) and (3.50),
$$\align
\pmatrix 1&{}\\{}&-1\endpmatrix'& = Ad(\eta)Ad(i_1(\tau_+))
\left(\pmatrix 1&{}\\{}&-1\endpmatrix\right)\tag3.60\\
\nass
{}&= \pmatrix 1&{}\\{}&-1\endpmatrix i_0(\eta_0^\s) i_1(\tau_+^{\s_1})i_1(\tau_+^{-1})i_0(\eta_0^{-1})\\
\nass
{}&= \pmatrix 1&{}\\{}&-1\endpmatrix i_0(\eta_0^\s\eta_0^{-1}) i_1(\e\tau^{\s_1}\tau^{-1}),
\endalign
$$
where $\e=\sgn(\nu(h))$, as above.
Now
$$\align
j_0' &= \eta j_0^\sh \eta^{-1}\\
\nass
{}&= Ad(\eta)Ad(i_1(\tau_+))\bigg( i_0(\bar\mu_0)
\pmatrix 1&{}\\{}&-1\endpmatrix  i_1( h^\iota )\bigg)\tag3.61\\
\nass
{}&= i_0(\bar\mu_0)\pmatrix 1&{}\\{}&-1\endpmatrix i_0(\eta_0^\s\eta_0^{-1}) i_1(\e\tau^{\s_1}\tau^{-1})
i_1(\tau h^\iota \tau^{-1})\\
\nass
{}&= \pmatrix 1&{}\\{}&-1\endpmatrix i_0(\bar\mu_0^\s\eta_0^\s\eta_0^{-1})
i_1(\e\tau^{\s_1} h^\iota \tau^{-1}).
\endalign$$
Setting this equal to
$$j_0= \pmatrix 1&{}\\{}&-1\endpmatrix i_1(h^{\s_1}) i_0(\mu_0),\tag3.62$$
recalling that $h^{\s_1}=h^\iota$, and dividing both sides by $\nu(h)$, we obtain
$$\align
i_0(\mu_0) &= i_1(h) i_0(\bar\mu_0^\s\eta_0^\s\eta_0^{-1}) i_1(\e\tau^{\s_1} h^{-1} \tau^{-1})\tag3.63\\
\nass
{}&=\e\, i_0(\bar\mu_0^\s\eta_0^\s\eta_0^{-1}) i_1(h\tau^{\s_1}h^{-1}\tau^{-1})\\
\nass
{}&=\e\, i_0(\bar\mu_0^\s\eta_0^\s\eta_0^{-1}),
\endalign
$$
since $h \tau^{\s_1} h^{-1} =\tau^\s=\tau$.
Thus, if we set
$$\rho =\cases \mu_0\eta_0 &\text{ if $\nu(h)>0$, }\\
\nass
\delta\mu_0\eta_0 &\text{ if $\nu(h)<0$,}
\endcases\tag3.64
$$
then $\rho\in E^\vee$ and
$$j_0 = i_1(h) \pmatrix 1&{}\\{}&-1\endpmatrix i_0(\rho \eta_{0,+}^{-1}),\tag3.65$$
where $\eta_{0,+}=\eta_0$ if $\nu(h)>0$ and $\eta_{0,+}=\delta\eta_0$ if $\nu(h)<0$.

Finally, $j_0^2=Q_0(j_0)\cdot id$ gives
$$Q_0(j_0) = \nu(h) N_{\smallkay/\Q}(\eta_0)^{-1}N_{E^\vee/\Q}(\rho)\cdot\cases 1&\text{ if $\nu(h)>0$,}\\
\nass
(-\Delta)^{-1}&\text{ if $\nu(h)<0$,}\endcases\tag3.66$$
as claimed. This proves the first part of Proposition~3.6.

Finally, consider case $\text{I}_2$,
where $A_0\simeq A_1^2$ for a simple abelian surface $A_1$
with $\End^0(A_1)=\kay'$ and with $\iota_0:B\hookrightarrow \End^0(A_0)\simeq M_2(\kay')$.
But then the isomorphism
$$B\tt_\kay\kay'\simeq M_2(\kay').\tag3.67$$
implies that $\iota_0(\kay)$ lies in the center of $\End^0(A_0)$, and so the
space of special endomorphisms is zero.
\qed\qed\enddemo

The results obtained so far impose restrictions on the isogeny
classes which can be met by a special cycle $\CZ(\o,T)$. Recall that we
assume that $\det(T)\ne0$. If
$\xi=(A,\lambda,\iota,\bar{\eta}^p;\j)\in \CZ(\o,T)(\F)$, then the components
of the collection $\j$
of special endomorphisms of $(A,\lambda,\iota)$ lie in $V_\xi$, and
$Q_\xi(\j)=T$.

\proclaim{Corollary 3.8} If $\xi\in \CZ(\o,T)(\F)$, then $T$ is represented by $V_\xi(\Q)$.
In particular:\hfill\break
(i) The isogeny classes of type $\text{I}_2$ do not meet any special cycle.\hfill\break
(ii) If $n\ge 3$, then $\CZ(\o,T)$ lies in the supersingular locus and is empty
unless $T>0$. \hfill\break
(iii) When $n=4$, $\CZ(\o,T)$ is empty unless $T$ is equivalent over $\Q$ to the
quadratic form of Proposition~3.5.\hfill\break
(iv)  When $n=3$, $\CZ(\o,T)$ is empty unless $T$ is represented over $\Q$ by the
quadratic form of Proposition~3.5.\hfill\break
(v) If $n=2$, then the cycle $\CZ(\o,T)$ can meet at most one nonsupersingular
isogeny class, namely, that for which $T$ is equivalent over $\Q$ to the quadratic form
$Q_\xi$ of Proposition~3.6 (i).
\endproclaim

\subheading{\Sec4. The structure of the supersingular locus}

In this section, we give a description of the supersingular locus $\mss$
in the special fiber of $\Cal M=\Cal M_{K^p}$.
The results are due Stamm \cite{\stamm}, although his presentation
is somewhat different. We continue to assume that $p$ is inert in
$\kay$. Then $\OC\tt\Z_p\simeq M_4(\Z_p)$ and $\OK\tt\Z_p\simeq
\Z_{p^2}$. Also let $W=W(\F)$ be the Witt ring of $\F$, $\K=W\tt\Q$
its quotient field, and $\s$ the Frobenius automorphism.

Fix a point
$\xi=(A,\lambda,\iota,\bar{\eta}^p)\in \Cal M^{\text{ss}}(\F)$ in the
supersingular locus. Let $A(p)$ be the p-divisible group of $A$, and
note that the inclusion
$$\iota:\OC\tt\Z_p\simeq M_4(\Z_p)\hookrightarrow \End(A(p))\tag4.1$$
yields a decomposition $A(p)\simeq {\Cal A}^4$ where ${\Cal A}$ is
a p-divisible group of dimension $2$ and height $4$ with an action
$\iota_0:\Z_{p^2}\hookrightarrow \End({\Cal A})$. Since $A$ is
supersingular, $A(p)$ and ${\Cal A}$ are formal groups. We let
$\text{\rm D}{\Cal A}$ be the contravariant Dieudonn\'e module of
${\Cal A}$. This is a free $W$-module of rank $4$ with operators
$F$ and $V$ which are $\s$- and $\s^{-1}$-linear respectively, and
with $FV=VF=p$. Let $\L= \text{\rm D}{\Cal A}\tt\Q$ be the
associated isocrystal, with $\dim_{\K}(\L)= 4$. The polarization of
$A$ gives rise to a nondegenerate alternating form $<\ ,\
>:\L\times
\L\rightarrow \K$, with $<Fx,y>\ =\ <x,Vy>^\s$. Moreover, for $a\in
\Z_{p^2}$,
$$<\iota_0(a)x,y>\ =\ <x,\iota_0(a)y>.\tag4.2$$
The fixed embedding $\Z_{p^2}\rightarrow W$ defines a
$\Z/2$-grading
$$\L_0=\{x\in \L;  \iota_0(a)x=a x, \text{for $a\in \Z_{p^2}$} \}\tag4.3$$
and
$$\L_1=\{x\in \L;  \iota_0(a)x=a^\s x, \text{for $a\in \Z_{p^2}$} \},\tag4.4$$
so that $\L=\L_0\oplus\L_1$, and $F$ and $V$ are endomorphisms of
degree $1$. Note that the two dimensional subspaces $\L_0$ and
$\L_1$ are orthogonal to each other with respect to $<\ ,\ >$. Let
$$G'_p=\{ g\in GL({\Cal L});\ < gx, gy> =\nu (g)\cdot
<x,y>,\ Fg=gF\ \text{with}\ \nu(g)\in \Q_p^{\times}\}\
\ .\tag4.5$$

We are interested in $W$-lattices $L$ in $\L$ which are stable
under $F$ and $V$ and under the action of $\Z_{p^2}$, and such that
$L^\perp = c L$ with respect to $<\ ,\ >$. The lattice $L$ is
stable under $\Z_{p^2}$ if and only if $L=L_0\oplus L_1$ with
respect to the grading. For such a lattice $L$, we have
$$L_0\supset FL_1\supset pL_0,\tag4.6$$
and
$$L_1\supset FL_0\supset pL_1,\tag4.7$$
where all inclusions have index $1$. This follows from the determinant condition,
which implies that $L_0/FL_1$ and $L_1/FL_0$ have dimension $1$ over $\F$.

Let $X$ be the set of such lattices in $\L$.

\proclaim{Definition 4.1} (i)  For a lattice $L\in X$,
the index $i$ is  critical for $L$ if $F^2L_i = p L_i$.\hfill\break
(ii) The lattice $L$ is called superspecial if both indices $0$ and $1$ are critical
for $L$.
\endproclaim

Let $X_0\subset X$ be the set of superspecial lattices.

\proclaim{Lemma 4.2} For any $L\in X$, at least one index is critical.
\endproclaim
\demo{Proof} Note that $F^2L_i=pL_i=FVL_i$ if and only if $FL_i=VL_i$. Since we have inclusions
$$\matrix FL_{i}&\subset&L_{i+1}\\
\nass
\cup&{}&\cup\\
\nass
pL_{i+1}&\subset& VL_i,\endmatrix
$$
we have that either $i$ is critical or $L_{i+1}=FL_i+VL_i$. Suppose
that there is no critical index, i.e., that $L_0=FL_1+VL_1$ and
$L_1=FL_0+VL_0$. But then, using the fact that $FL_1\supset pL_0$,
$$\align
L_0&= FL_1+VL_1\\
{}&= FL_1+V(FL_0+VL_0) = FL_1+V^2L_0\\
{}&= FL_1+V^2(FL_1+VL_1)= FL_1+V^3L_1\\
{}&=\dots\\
{}&= FL_1 +V^rL_i\quad\text{ for all $r\ge1$ }\\
{}&= FL_1,\endalign
$$
since the Dieudonn\'e module is $V$-reduced. This contradicts the
fact that $L_0/FL_1$ has dimension $1$ over $\F$.
\qed\enddemo

Note that the $\K$-vector space $\L_i$ has a $\s^2$-linear
automorphism $p^{-1}F^2$. If $i$ is a critical index for the
lattice $L\in X$, then $p^{-1}F^2$ preserves $L_i$ and defines an
$\fps$-rational structure on the $\F$-vector space $L_i/pL_i$ and
on the projective line $\Bbb P(L_i/pL_i)$.

\proclaim{Lemma 4.3}
(i) Suppose that $L_i$ is a $W$-lattice in $\L_i$ such that
$F^2L_i=pL_i$. For any line $\ell$ in the two dimensional
$\F$-vector space  $L_i/pL_i$, let $L_{i+1}= F^{-1}(\ell+pL_i)$,
and let $L=L_{i+1}\oplus L_i$. Then $L\in X$. \hfill\break
(ii) If
$L=L_0\oplus L_1\in X$ with $i$ a critical index, let $\ell\in \Bbb
P(L_i/pL_i)$ be the line $FL_{i+1}/pL_i$. Then $L$ is recovered
 from $\ell$ by the construction of (i).\hfill\break
(iii) The
lattice $L=L_{i+1}\oplus L_i\in X$ associated to $\ell\in \Bbb
P(L_i/pL_i)$ is superspecial if and only if the line $\ell$ is
rational over $\Bbb F_{p^2}$.
\endproclaim

\demo{Proof} Note that $FL_i=VL_i$. Since $L_{i+1} = F^{-1}(\ell+pL_i)\subset F^{-1} L_i$,
we clearly have $FL_{i+1}\subset L_i$ and $VL_{i+1}\subset L_i$.
Also, $L_{i+1}\supset F^{-1}pL_i = FL_i=VL_i$. Next, since the
restriction of $<\ ,\ >$ to the two dimensional space $\L_i$ is a
nondegenerate alternating form, we must have $<L_i,L_i>=c W$ for
$c=p^r$, for some $r\in \Z$, and $L_i^\perp = c^{-1}L_i$. But then,
writing $\ell+pL_i = W y+pL_i$, we have
$$\align
<\ell+pL_i,\ell+pL_i>\ &=\ <Wy+pL_i,Wy+pL_i>\\
{}&=\ <y,pL_i>\\
{}&=\ p<L_i,pL_i>\\
{}&=  pcW.\endalign
$$
But then,
$$\align
<L_{i+1},L_{i+1}>\ &=\ <F^{-1}(\ell+pL_i),p^{-1}V(\ell+pL_i)>\\
{}&= p^{-1}<\ell+pL_i,\ell+pL_i>^{\s^{-1}}\\ {}&= cW,\endalign
$$
and so $(L_{i+1})^\perp = c^{-1} L_{i+1}$. This proves (i), and
(ii) is obvious.

To prove (iii), observe that $FL_{i+1}=VL_{i+1}$ if and only if
$\ell+pL_i = VF^{-1}(\ell)+pL_i$, and that $VF^{-1} = pF^{-2}$.
\qed
\enddemo

The set $X$ can be described as the set of ${\Bbb F}$-points of a
scheme over ${\Bbb F}_{p^2}$, cf.\ \cite{\stamm}. Let $Y_i$ be the
set of $W$-lattices in ${\Cal L}_i$ which satisfy $F^2L_i= pL_i$.
Let
$$U_i= {\Cal L}_i^{pF^{-2}}\ \ .\tag4.8$$
Since $pF^{-2}$ is a $\sigma^{-2}$-linear automorphism with all
slopes equal to 0, $U_i$ is a 2-dimensional $\Q_{p^2}$-vector space
with $U_i\otimes_{\Q_{p^2}}{\Cal K}={\Cal L}_i$ and each $L_i\in
Y_i$ can be written as
$$L_i=\Lambda_i\otimes_{\Z_{p^2}}W\tag4.9$$
for the $\Z_{p^2}$-lattice $\Lambda_i=L_i^{pF^{-2}}$ of $U_i$. To
each $L_i\in Y_i$ we associate the projective line
$${\Bbb P}_{L_i}={\Bbb P}(\Lambda_i/p\Lambda_i)$$
over ${\Bbb F}_{p^2}$. At each of its $p^2+1$ ${\Bbb
F}_{p^2}$-rational points, ${\Bbb P}_{L_i}$ meets a single ${\Bbb
P}_{L_{i+1}}$.

For the following result, cf. \cite{\stamm}.
\proclaim{Proposition 4.4} There is a natural $G_p'$-equivariant surjective
map
$$\coprod_{i=0,1}\coprod\limits_{L_i\in Y_i} \Bbb P_{L_i}(\F) \lra X,$$
which induces a bijection
$$\coprod_{i=0,1}\coprod\limits_{L_i\in Y_i}
\bigg(\Bbb P_{L_i}(\F) -\Bbb P_{L_i}(\fps)\bigg)\lra X-X_0.$$
\endproclaim
For $L\in X_0$ the preimage consists of two points, one for each
$i=0$ and $i=1$.

Let $G'$ be the inner form of $G$ (cf.\ (0.2)) given by the
quaternion algebra $B'$ over $\kay$ which is ramified at the two
archimedean primes and is isomorphic to $B$ at all finite primes.
Then there exists a morphism of schemes over ${\roman{Spec}}\,
{\Bbb F}_{p^2}$ where ${\Cal M}_{K^p}^{ss}$ denotes the
supersingular locus of ${\Cal M}_{K^p}\times_{{\roman{Spec}}\,
\Z_{(p)}} \Spec\, {\Bbb F}_p$,
$$G'(\Q)\setminus \left[ \left( \coprod\limits_{i=0,1}
\coprod_{L_i\in Y_i} {\Bbb P}_{L_i}\right)\times G'({\Bbb A}_f^p)/K^p
\right]\longrightarrow {\Cal M}^{ss}\otimes_{{\Bbb F}_p} {\Bbb
F}_{p^2}\ \ .\tag4.10$$ This morphism can be identified with the
normalization of ${\Cal M}^{ss}\otimes_{{\Bbb F}_p} {\Bbb
F}_{p^2}$, cf.\ \cite{\stamm}.
It induces a bijection
$$\coprod\limits_{i=0,1}G'(\Q)\setminus G'({\Bbb A}_f)/
K_p.K^p\buildrel\sim\over\longrightarrow {\roman{Irred}}({\Cal
M}^{ss}\otimes_{{\Bbb F}_p}{\Bbb F}_{p^2})\ \ .\tag4.11$$
Here on
the right hand side is the set of geometric irreducible components
of ${\Cal M}^{ss}$ and $K_p$ is a maximal open compact subgroup of
$G'_p\cong G'(\Q_p)$ (the stabilizer of a $\Z_{p^2}$-lattice
$\Lambda\subset U_i$).

Similarly we obtain for the set of superspecial points of ${\Cal
M}^{ss}$ a bijection
$$G'(\Q)\setminus G'({\Bbb
A}_f)/K'_p.K^p\buildrel\sim\over\longrightarrow {\Cal
M}^{supersp}(\overline{\Bbb F}_p)\ \ .\tag4.12$$
Here $K'_p$ is an
Iwahori subgroup of $G'_p$. Indeed, it follows from Lemma 5.1 below
that all superspecial lattices are conjugate under $G'_p$
(existence of a standard basis) and the stabilizer of a
superspecial lattice is an Iwahori subgroup.

\subheading{\Sec5. Special endomorphisms of supersingular
Dieudonn\'e modules}

In this section we study the space $V'_p$ (cf.(5.3)) of special endomorphisms of the isocrystal $\L$
associated to the supersingular isogeny class as in section 4. After giving an explicit description
of this quadratic space over $\Q_p$, we provide a criterion for a given special endomorphism $j$ of $\L$
to induce a special endomorphism of a lattice $L\in X$ (resp. every lattice in $\Bbb P_{L_i}$ for
$L_i\in Y_i$). This information will be used to determine isolated supersingular point on special cycles
(section 6) and supersingular components of special cycles (section 8).

We keep the notations of the previous section. We begin by
considering the space of special endomorphisms of the $\K$-vector
space $\L$. Recall that $\K= W\otimes \Q$. The grading
$\L=\L_0\oplus
\L_1$ gives a grading on the endomorphism ring:
$$\End(\L)^{(0)}=\{ j\in \End(\L)\mud  j\iota_0(a) = \iota_0(a) j\},\tag5.1$$
and
$$\End(\L)^{(1)}=\{ j\in \End(\L)\mud  j\iota_0(a) = \iota_0(a^\s) j\}.\tag5.2$$
Let $\End(\L,F)$ be the algebra of endomorphisms which commute with $F$, and
let
$$V_p':=\{j\in \End(\L,F)^{(1)}\mud  j^*=j\}\tag 5.3$$
be the space of special endomorphisms of $\L$,
where $*$ denotes the adjoint with respect to the alternating form $<\ ,\ >$.
Writing $j\in \End(\L,F)^{(1)}$ as
$$j=\pmatrix {}&j_1\\j_0&{}\endpmatrix,\qquad\text{ $j_0\in \Hom(\L_0,\L_1)$,
and $j_1\in \Hom(\L_1,\L_0)$,}\tag5.4$$
we have $j^*=j$ if and only if $j_1=j_0^*$.
Since $j$ commutes with $F$, we also have $j_1=Fj_0F^{-1}$ and $j_0=Fj_1F^{-1}$,
so that
$$V_p'\simeq \{ j_0\in \Hom(\L_0,\L_1)\mud j_0=F^2j_0F^{-2},\  j_0^*=Fj_0F^{-1}\}.\tag5.5$$
If $x=j_0^*j_0^{\phantom{*}}\in \End(\L_0,F^2)$, then $x$ is self
adjoint with respect to $<\ ,\ >$ on the $2$-dimensional space
$\L_0$, and hence $x= a\cdot 1_{\L_0}$ for a scalar $a$. Thus, we
have a quadratic form on $V'_p$, defined by
$$j^2 = Q'(j)\cdot 1_\L\ \ ,\ \ Q'(j)= j^*_0\cdot j_0\ \ .\tag 5.6$$
Since the elements of $G'_p$ commute with $F$, cf\ (4.5), we have
an inclusion $G'_p\subset \End({\Cal L})^{(0),\times}$ and $G'_p$
acts on $V'_p$ by conjugation and preserves the quadratic form.

Choose a $\Q_{p^2}$-rational basis $e_1$, $e_2$ for $U_0={\Cal
L}_0^{pF^{-2}}$, with $<e_1,e_2>=p^\ell$, and let $e_3=Fe_1$ and
$e_4=p^{-1}Fe_2$, so that $<e_3,e_4>=p^\ell$ as well. We call such
a basis a {\it standard basis for}\/ $\L$. Then, we may write
$$F=\pmatrix {}&{}&p&{}\\{}&{}&{}&1\\1&{}&{}&{}\\{}&p&{}&{}\endpmatrix \s
=\pmatrix {}&\pi\\p\pi^{-1}&{}\endpmatrix\s,\tag5.7$$
where we set $\pi=\pmatrix p&{}\\{}&1\endpmatrix$.
If $j_0=\pmatrix a&b\\c&d\endpmatrix$, then $j_0^* =\pmatrix d&-b\\-c&a\endpmatrix$,
and we have $jF=Fj$ if and only if
$$j_1=p^{-1}\pi j_0^\s \pi,\qquad\text{ and } \qquad j_1^\s=p^{-1} \pi j_0 \pi.\tag5.8$$
Applying $\s$ to the first condition and comparing to the second, we have $j_0^{\s^2}=1$.
Recalling that $j_1=j_0^*$, we obtain
$$\pmatrix d&-b\\-c&a\endpmatrix = \pmatrix pa^\s& b^\s\\c^\s&p^{-1}d^\s\endpmatrix,\tag5.9$$
i.e., $d=pa^\s$, $c=-c^\s$, and $b=-b^\s$.
Thus,
$$V_p'\simeq \{x= \pmatrix a&b\\c&pa^\s\endpmatrix ;
\ a,\, b,\, c\in \Q_{p^2},\ b^\s=-b,\ c^\s=-c \}\tag 5.10$$
is a
$4$-dimensional vector space over $\Q_p$ with quadratic form
$$Q'(x) = paa^\s-bc.\tag5.11$$

Also,
$$G'_p\simeq \{g\in GL_2(\qps);\  \det(g)\in \Q_p^\times\},\tag5.12$$
and $g\in G'_p$ acts on $x\in V'_p$ by
$$g:\ x\mapsto \pi^{-1}g^\s \pi x g^{-1}.\tag5.13$$

For any lattice $L\in X$, the ring $\End(L,F)$
is an order in the algebra $\End(\L,F)$, and
$$N_L := \End(L,F)\cap V_p'\tag 5.14$$
is a $\Z_p$-lattice in $V_p'$.

\proclaim{Lemma 5.1} Suppose that $L\in X_0$ is a superspecial lattice with
$<L,L>=p^\ell W$. Then, there exists a {\it standard basis} $e_1,\,e_2,\, e_3,\, e_4$ for
$\L$ such that $L_0=We_1\oplus We_2$ and $L_1=We_3\oplus We_4$.
\qed
\endproclaim

We refer to the basis of Lemma~5.1
as a {\it standard basis} for $L\in X_0$.

\proclaim{Corollary 5.2} If $L$ is a superspecial lattice in $\L$, then with
respect to a standard basis for $L$,
$$\End(L,F)^{(0)} =
\{\pmatrix A_0&{}\\{}&A_1\endpmatrix ;\  A_0=\pmatrix a&pb\\c&d\endpmatrix,
\text{ and } A_1= \pi^{-1}A_0^\s \pi,\ a,\,b,\, c,\,d\in \Z_{p^2}\},$$
and
$$\End(L,F)^{(1)} =
\{\pmatrix {}& A_1\\A_0&{}\endpmatrix ;\  A_0=\pmatrix a&b\\c&pd\endpmatrix,
\text{ and } A_1=p^{-1}\pi A_0 \pi,\  a,\,b,\, c,\,d\in \Z_{p^2}\}.$$
In particular,
$$N_L =\End(L,F)\cap V'_p \simeq
\{x= \pmatrix a&b\\c&pa^\s\endpmatrix ;\
 a,\, b,\, c\in \Z_{p^2},\ b^\s=-b,\ c^\s=-c\ \}.$$
\endproclaim

If $L\in X_0$ is a superspecial lattice, then $p^{-1}F^2L=L$, and
we let $\Lambda$ be the set of fixed points of $p^{-1}F^2$. This is
a $\Z_{p^2}$-module of rank $4$, and $\Lambda/F\Lambda$ is an
$\fps$-vector space of dimension $2$. Note that the vectors $e_i$
of a standard basis for $L$ lie in $\Lambda$, and the images
$\bar{e}_1$ and $\bar{e}_4$ give an $\fps$-basis for
$\Lambda/F\Lambda$. There is a natural reduction map:
$$\red_L:\End(L,F)\rightarrow \End(\Lambda/F\Lambda).\tag5.15$$
\proclaim{Lemma 5.3} With respect to a fixed standard basis,
$$\red_L:\End(L,F)\rightarrow \End(\Lambda/F\Lambda)\simeq M_2(\fps)$$
is given by
$$\red_L:\pmatrix A_0&{}\\{}&A_1\endpmatrix \mapsto \pmatrix \bar{a}&{}\\{}&\bar{d}^\s\endpmatrix,$$
and
$$\red_L:\pmatrix {}& A_1\\A_0&{}\endpmatrix \mapsto
\pmatrix {}&\bar{b}^\s\\\bar{c}&{}\endpmatrix.$$
Here $a\mapsto \bar{a}$ is the reduction map $\Z_{p^2}\rightarrow \fps$,
and the other notation is as in the Corollary~5.2.
\endproclaim

Note that the grading on $M_2(\fps)$ is the checkerboard grading.

Let
$$\nn_L = \red_L(N_L) = \{ x=\pmatrix {}&-b\\c&{}\endpmatrix
;\  b,\ c\in \fps,
\text{ with } b^\s=-b,\ c^\s=-c \},\tag5.16$$
with quadratic form $\bar{Q}'(x) = bc$. Then, we have the
commutative diagram
$$\matrix N_L&\overset{Q'}\to{\lra}& \Z_p\\
\nass
\hskip -24pt\red_L\downarrow&{}&\downarrow\\
\nass
\nn_L&\overset{\bar{Q}'}\to{\lra}&\fps\endmatrix\tag5.17
$$
of quadratic spaces.

We next turn to the non-superspecial  case.
\proclaim{Proposition 5.4} (i)
Suppose that $L\in X-X_0$ is a non-superspecial  lattice with critical
index $i$, and write $L=L_i\oplus L_{i+1}$.
Then $j\in \End(\L,F)^{(1)}$, lies in $\End(L,F)^{(1)}$ if and only if
$$j_i(L_i)\subset FL_i,$$
where
$$j=\pmatrix {}&j_1\\j_0&{}\endpmatrix,$$
as above. \hfill\break
(ii) If  $L\in X-X_0$ is non-superspecial ,
and if $j\in N_L$ is a special endomorphism of $L$, so that $j\in \End(L,F)^{(1)}$ with $j^*=j$,
 then
$\ord_p(Q'(j))\ge 1$.
\endproclaim
\proclaim{Proposition 5.5} Fix an index $i$.
Suppose that $L_i\subset \L_i$ is a $W$-lattice such that
$F^2L_i=pL_i$, i.e., $\L_i\in Y_i$, in the notation of section 4.
Then, a special endomorphism $j\in V'_p$
induces an endomorphism of every lattice in $\Bbb P_{L_i}$
if and only if
$$j(L_i)\subset FL_i.$$
\endproclaim
\proclaim{Corollary 5.6} If $j\in \End(L,F)^{(1)}$ for one
non-superspecial  point $L$ in $\Bbb P_{L_i}$, then
$j\in\End(L',F)^{(1)}$ for all points $L'$ of $\Bbb P_{L_i}$.
\endproclaim
\demo{Proof of Proposition~5.4}
Since $i$ is critical,
$FL_i=VL_i$. Recall that $L_{i+1}= F^{-1}(\ell+L_i)$
where $\ell \in \Bbb P(L_i/pL_i)$ is a line, and that, since $i+1$ is not
critical,
$$L_i=FL_{i+1}+VL_{i+1}\qquad\text{ and }\qquad pL_i=FL_{i+1}\cap VL_{i+1}.\tag5.18$$
Recall that $j_0$ and $j_1$ commute with
$VF^{-1}=(p^{-1}F^2)^{-1}$. Then, if $j(L)\subset L$, we have
$j_i(L_i)\subset L_{i+1}$, and thus, $Fj_i(L_i)\subset FL_{i+1}$.
Applying $VF^{-1}$ to this, we get $Fj_i(L_i) \subset VL_{i+1}$.
Hence
$$Fj_i(L_i) \subset FL_{i+1}\cap VL_{i+1}=pL_i,\tag5.19$$
i.e., cancelling an $F$,
$$j_i(L_i)\subset VL_i=FL_i.\tag5.20$$
When $i$ is critical and $i+1$ is not critical, condition (5.20) is
in fact equivalent to the requirement that $j_i(L_i)\subset
L_{i+1}$. Indeed, otherwise $L_{i+1}=j_i(L_i)+VL_i$, hence
$VL_{i+1}=FL_{i+1}$, i.e.\ $i+1$ is critical. We must also check
that $j_{i+1}(L_{i+1})\subset L_i$. But, this is automatic, since
$$Fj_{i+1}(L_{i+1}) = j_i(FL_{i+1}) = j_i(\ell+pL_i) \subset j_i(L_i)\subset FL_i,\tag5.21$$
again by (5.20). Thus, condition (5.20) is all that is needed for
$j(L)\subset L$, and (i) of the Proposition is proved.

To prove (ii), recall that for a special endomorphism $j$,  $j_i^*j_i^{\phantom{*}}=Q'(j)\cdot
1_{L_i}$, by (5.6). Then, since, by (5.20),
$$<L_i,j_i^*j_i^{\phantom{*}}(L_i)>\ =\ <j_i(L_i),j_i(L_i)>\ \subset \ <FL_i,FL_i>\ = p<L_i,L_i>,\tag5.22$$
we have $j_i^*j_i^{\phantom{*}}(L_i)\subset pL_i$, and hence $\ord_p(Q'(j))\ge1$, as claimed.
\qed\enddemo
\demo{Proof of Proposition~5.5} If $j$ extends to every lattice in $\Bbb P_{L_i}$, then
$j$ preserves a non-superspecial  lattice with critical index $i$, and hence, by Proposition~5.4 (i),
$j_i(L_i)\subset FL_i$. Conversely, suppose that $j_i(L_i)\subset FL_i$. Let
$L=L_i\oplus L_{i+1}$, with $L_{i+1} =F^{-1}(\ell+pL_i)$, be a lattice in $\Bbb P_{L_i}$.
Then $j(L_i)\subset L_{i+1}$ if and only if $Fj(L_i)\subset FL_{i+1}=\ell+pL_i$.
But $F j(L_i) \subset F^2L_i=pL_i$, so the required inclusion follows. Similarly,
$j(L_{i+1})\subset L_i$ if and only if $Fj(L_{i+1})\subset FL_i$, but
$Fj(L_{i+1})=j(\ell+pL_i) \subset j(L_i) \subset FL_i$, as required.
\qed\enddemo

In section 8 below we will fix a special endomorphism $j\in V'_p$
of the isocrystal $\L$ and give a combinatorial description of the
set of lattices $L_i\in Y_i$ for which $j$ induces a special
endomorphism of every lattice in $\Bbb P_{L_i}$.

\subheading{\Sec6. Isolated supersingular points of special cycles}

In this section we return to the special cycles introduced in
section 2 and determine the isolated supersingular points on them. This allows us to
characterize the isolated points of intersection of our special cycles and to
obtain a formula for the `intersection multiplicity' at such points (Corollary~6.3).

We fix $n$ with $1\leq n\leq 4$ and consider the special cycle
${\Cal Z}(T,\omega)$, where $\omega\subset V({\Bbb A}_f^p)^n$ and
$T\in{\roman{Sym}}_n(\Q)$. We fix a base point $\xi_0 =(A_0 ,
\lambda_0 ,\iota_0 , \overline\eta_0 ^p)\in {\Cal M}_{K^p}^{ss}({\Bbb
F})$. As in section 4 we have a decomposition of the $p$-divisible
group $A_0 (p)= {\Cal A}_0 ^4$ and we introduce the isocrystal ${\Cal
L}=D {\Cal A}_0 \otimes \Q$. To every isogeny
$\mu:\xi\buildrel\sim\over\longrightarrow \xi_0 $ with $\xi= (A,
\lambda, \iota, \overline{\eta}^p)\in {\Cal M}_{K^p}^{ss}({\Bbb
F})$ we then associate the Dieudonn\'e lattice
$L= \mu_*(D{\Cal A})$ in ${\Cal L}$ which lies in $X$. The point $\xi$ will be
called {\it superspecial}\/ if $L\in X_0$ and {\it
non-superspecial}\/ if $L\in X-X_0$. This is independent of the
choice of the isogeny $\mu$. Suppose now that $\xi$ is the image of
$(A, \lambda, \iota,
\overline\eta^p; \j)\in {\Cal Z}(T,\omega) ({\Bbb F})$. Then $\j$
induces an $n$-tuple
$$\j_{\Cal L}\in {V'_p}^n\tag6.1$$
of special endomorphisms of ${\Cal L}$ with $\j_{\Cal L}\in
(N_L)^n$. Moreover the condition that
$Q(\eta\circ\j\circ\eta^{-1})=T$ implies that also
$$Q'(\j_{\Cal L})=T\ \ ,\tag6.2$$
with $Q'$ as in (5.6).

\proclaim{Theorem 6.1} Consider the set of supersingular points in the
image of ${\Cal Z}(T,\omega)$, for $T\in Sym_n(\Q)$ with
$\det(T)\ne0$ and for  $\o\subset V(\A_f^p)^n$.\hfill\break (i) If
$T\notin Sym_n(\Z_{(p)})$, i.e., if $T$ is not $p$-integral, then
$\CZ(T,\o)\cap\mss(\F)$ is empty.
\hfill\break
Assume that $T\in Sym_n(\Z_{(p)})$, and let $\bar{T}\in
Sym_n(\F_p)$ be its reduction modulo $p$.\hfill\break (ii) If the
rank of $\bar{T}$ is greater than $2$, or if $\bar{T}$ has rank $2$
and is anisotropic modulo its radical, then $\CZ(T,\o)\cap\mss(\F)$
is empty.
\hfill\break
(iii) If $\bar{T}\ne 0$, then $\CZ(T,\o)\cap\mss(\F)$ consists
entirely of superspecial points. \hfill\break
(iv)  Suppose that
$\bar{T}=0$ and that $\CZ(T,\o)\cap\mss(\F)\ne \phi$. If $L$ is a non-superspecial
point in the image of $\CZ(T,\o)$, and if $L\in \Bbb P_{L_i}$, then
$\Bbb P_{L_i}\subset \CZ(T,\o)$. \hfill\break
(v) Suppose that $\bar{T}=0$ and that $L=L_0\oplus L_1\in X_0$ is a superspecial
lattice corresponding to a point of $\xi\in\mss(\F)$ in the image
of ${\Cal Z}(T,\omega)$. As usual let $\Lambda
=L^{pF^{-2}}=\Lambda_0\oplus \Lambda_1$. Let $M\subset N_L$ be the $\Z_p$-submodule spanned
by the components of $\j_{\Cal L}$, and let $\frak m = \red_L(M)$
be the image of $M$ in $\nn_L\subset
\End(\Lambda/F\Lambda)$.
\roster
\item"{(a)}"  If $\frak m=0$,
then both components $\Bbb P_{L_0}$ and $\Bbb P_{L_1}$ of $\mss$
through $\xi$ (cf.\ (4.10)) lie in $\CZ(T,\o)\cap\mss(\F)$.
\item"{(b)}"  If $\frak m\ne 0$, then $\dim_\fp \frak m=1$.
A basis vector $\bar{j}\in\frak m$
is a nonzero endomorphism of
$$\Lambda/F\Lambda=(\Lambda_0/F\Lambda_1)\oplus(\Lambda_1/F\Lambda_0)$$
of degree 1 with $\bar{j}^2=0$. Therefore the condition
$$ j(L_i)\subset F(L_i)$$
holds for precisely one index $i\in\Z/2$, in which case, $\Bbb
P_{L_i}$ lies in $\CZ(T,\o)\cap\mss(\F)$ and $\Bbb P_{L_{i+1}}$
does not.
\endroster
\endproclaim

\demo{Remark} In effect, if $p\nmid T$, then $\CZ(T,\o)^{ss}$ consists of
(at most) an isolated (finite) set of superspecial points, while,
if $p\mid  T$, then the image of ${\Cal Z}(T,\omega)$ in ${\Cal
M}^{ss}$ is either empty or consists of a union of components of
$\mss$. We will give a more complete description of the components
which occur in section 8 below.
\enddemo

\demo{Proof} Suppose that
$\tilde\xi=(A,\lambda,\iota,\bar{\eta}^p;\j)\in \CZ(T,\o)$ lies
above $\xi=(A, \lambda, \iota, \overline\eta^p)\in {\Cal
M}^{ss}({\Bbb F})$. Since all of our assertions are `local' at $\xi$, we fix an
isogeny $\mu:\xi\isoarrow \xi_0$ and hence associated with $\xi$ a Dieudonn\'e lattice
$L=L_0\oplus L_1\in X$ and
$n$-tuple of special endomorphisms $\j_{\Cal L}\in N_L^n\subset
(\End(L,F)^{(1)})^n$. Let $M$ be the $\Z_p$-submodule of $N_L$
spanned by the components of $\j_{\Cal L}$. Since $Q'$ is
$\Z_p$-valued on $N_L$, it follows that $T\in Sym_n(\Z_{(p)})$.
This proves (i). (In fact, from the definitions in section 2 we
know that the whole special fibre of ${\Cal Z}(T,\omega)$ is empty
if $T$ is not $p$-integral.) Since we are assuming that $\det(T)\ne
0$, $M$ must have rank $n$ over $\Z_p$.

First suppose that $\xi$ is a non-superspecial point. Then, by (ii) of Proposition~5.4,
$\bar{T}=0$, since $Q'(j)\equiv 0 \mod p$ for every $j\in M$. This proves (iii).

If $\xi$ is a superspecial point in the image of ${\Cal Z}(T,\omega)$, then $\bar{T}$
is the matrix of inner products of the images of the components of
$\j_{\Cal L}$ under the reduction map $\red_L:N_L\rightarrow
\nn_L$. Since $\nn_L$ is a hyperbolic plane over $\fp$, it follows
that $\bar{T}$ must have rank at most $2$, and that, if the rank of
$\bar{T}$ is $2$, then $\frak m=\nn_L$ and $\bar T$, modulo its
radical, must be a hyperbolic plane. This proves (ii).

Finally, suppose that $\bar{T}=0$. If $\xi$ is a non-superspecial
point in $\CZ(T,\o)\cap \mss(\F)$, with critical index $i$, then,
by Corollary~5.6, the whole component $\Bbb P_{L_i}$ lies in
$\CZ(T,\o)\cap \mss(\F)$. This proves (iv). If $\xi$ is a superspecial point in
$\CZ(T,\o)\cap \mss(\F)$, then $\frak m$ is an isotropic subspace
of $\nn_L$ and thus has dimension $0$ or $1$. If $\frak m=0$, then,
for every $j\in M$, the condition $j(L_i)\subset FL_i$ is satisfied
for both $i=0$ and $1$. Thus, by Proposition~5.5, both components
$\Bbb P_{L_0}$ and $\Bbb P_{L_1}$ of $\mss$ through $\xi$ lie in
$\CZ(T,\o)\cap \mss(\F)$. If $\frak m$ is an isotropic line in
$\nn_L$, then there is a unique index $i$ such that $j(L_i)\subset
FL_i$ holds for all $j\in M$. Thus, in this case, $\Bbb P_{L_i}$
lies in $\CZ(T,\o)\cap \mss(\F)$, but $\Bbb P_{L_{i+1}}$ does not.
This proves (v).
\qed\enddemo

We now assume that $n\ge3$ and that $T\in
{\roman{Sym}}_n(\Z_{(p)})_{>0}$ and also that $\bar{T}\ne 0$. The
first assumption implies that $\CZ(T,\o)$ lies over the
supersingular locus $\Cal M^{\text{ss}}$ of the special fiber of
$\Cal M$. The second assumption implies that $\CZ(T,\o)=\CZ(T,\o)^{ss}$, if
nonempty, consists of a finite set of isolated superspecial points.

We fix a supersingular point $\xi=(A, \lambda, \iota,
\overline\eta^p; \j)\in {\Cal Z}(T,\omega)({\Bbb F})$. Again we
write $A(p)={\Cal A}^4$, where ${\Cal A}$ is a 2-dimensional formal
group of height 4, equipped with an action
$$\iota_0: \Z_{p^2}\longrightarrow{\roman{End}}({\Cal A})$$
and with a principal quasi-polarization
$$\lambda_{\Cal A}:{\Cal A}\buildrel\sim\over\longrightarrow
\hat{\Cal A}\ \ $$
such that $\iota_0(a)^*=\iota_0(a)$.

By the Serre-Tate Theorem, the formal completions at $\xi$ of $\Cal
M$ and of $\CZ(T,\o)$ can be interpreted as versal deformation
spaces:
$$\hat{\Cal M}_\xi = \text{Def}(\Cal A,\lambda_{\Cal A}, \iota_0),\tag6.3$$
and
$$\hat{\CZ}(T,\o)_\xi = \text{Def}(\Cal A,\lambda_{\Cal A},\iota_0;\j) =
\text{Def}(\Cal A,\lambda_{\Cal A},\iota_0;M)\tag6.4 $$
where $M$ is the $\Z_p$-submodule of $\End(\Cal A)$ spanned by the
components of $\j$. By our assumption $\overline T\ne 0$, the
latter deformation space is the spectrum of a local Artin ring
$R_{\xi}$. We let $e(\xi)$ be the length of $R_{\xi}$.

From now on, we assume that $n=3$, and we reduce the computation of
$e(\xi)$ to a result of Gross and Keating, \cite{\grosskeating}.

Since, as always, $p\ne 2$, we may choose a $\Z_p$-basis $\psi_1,\
\psi_2,\ \psi_3$ for $M$ such that the matrix for the restriction of
the quadratic form $Q_\xi$ to this basis is
$$T' = \diag(\varepsilon_1 p^{a_1},\varepsilon_2 p^{a_2},
\varepsilon_3 p^{a_3}),\tag6.5$$
with $0\le a_1\le a_2\le a_3$ and with units $\varepsilon_i$
uniquely determined modulo squares. In particular, $\psi_i^2 =
\varepsilon_i
p^{a_i}$, and $\psi_i\psi_j+\psi_j\psi_i=0$ for $i\ne j$.

Recall that $p$ is inert in $k$ so that $k_p=\Q_p(\delta)$ with
$\delta^2=\Delta\in\Z_p^{\times}$ not a square. Since $\overline
T\ne 0$, we have $a_1=0$ and we may take $\varepsilon_1=1$ or
$\varepsilon_1=\Delta$, depending on whether or not $\varepsilon_1$
is a square. If $\varepsilon_1=1$ define idempotents
$$e_0=\frac{1}{2}(1+\psi_1)\ \ ,\ \ e_1=\frac{1}{2}(1-\psi_1)\ \
.\tag6.6$$
These yield a decomposition ${\Cal A}\simeq {\Cal A}_0\times {\Cal
A}_1$ where ${\Cal A}_i=e_i{\Cal A}$ has dimension 1 and height 2.
Since $e_i^{\ast}=e_i$, the polarization $\lambda_{\Cal A}$ is of
the form
$$\lambda_{\Cal A}=\lambda_{{\Cal A}_0}\times \lambda_{{\Cal
A}_1}: {\Cal A}_0\times {\Cal A}_1\buildrel\sim\over\longrightarrow
\hat{\Cal A}_0\times \hat{\Cal A}_1\ \ .\tag6.7$$
Let $\underline{\delta}=\iota_0(\delta)\in {\roman{End}}({\Cal
A})$. Then $\underline{\delta}$ is an isomorphism with
$\underline{\delta} e_0 =e_1\underline{\delta}$ and hence may be
considered as an isomorphism ${\Cal A}_0\to {\Cal A}_1$.
Furthermore
$$\psi_ie_0= e_1\psi_i\ ,\ \ i=2,3\ \ .\tag6.8$$
We put
$$(\mu_1, \mu_2, \mu_3)=(\underline\delta, \psi_2, \psi_3)\ \
,\tag6.9$$
so that $(\mu_1, \mu_2, \mu_3)$ is a triple of isogenies
from ${\Cal A}_0$ to ${\Cal A}_1$. Since a principal
quasi-polarization on a formal group of dimension 1 and height 2
deforms automatically and using the fact that $\psi_2$ and $\psi_3$
are special endomorphisms we see that
$${\roman{Def}}({\Cal A}, \lambda_{\Cal A}, \iota_0; M)=
{\roman{Def}} ({\Cal A}_0, {\Cal A}_1; \mu_1, \mu_2,
\mu_3)\ \ ,\tag6.10$$
where the right side is the locus inside the universal deformation
space of the pair $({\Cal A}_0, {\Cal A}_1)$ to which the triple of
isogenies deforms. We note that the degree quadratic form on
${\roman{Hom}}({\Cal A}_0, {\Cal A}_1)$ has matrix equal to
$$T''=\diag(\Delta, \varepsilon_2p^{a_2}, \varepsilon_3p^{a_3})\tag6.11$$
on the subspace of rank 3 spanned by the $\mu_i$'s.

Next assume that $\varepsilon_1=\Delta$ and write
$\Delta^{-1}=a^2+b^2$ with $a,b\in\Z_p^\times$. Let
$$\psi'_1=a\underline\delta +b\psi_1\in {\roman{End}}({\Cal A})\ \
.\tag6.12$$
Then ${\psi'_1}^2={\roman{id}}$. We define idempotents
$$e_0= \frac{1}{2} (1+\psi'_1)\ ,\ e_1= \frac{1}{2}(1-\psi'_1)\ \
.\tag6.13$$
Again we obtain a decomposition ${\Cal A}={\Cal A}_0\times {\Cal
A}_1$ with ${\Cal A}_i=e_i{\Cal A}$. Since again $e_i^*=e_i$ the
polarization $\lambda_{\Cal A}$ again splits as a product. Put
$$\psi''_1=b\underline\delta -a\psi_1\ \ .\tag6.14$$
Then we have
$$\psi''_1e_0= e_1\psi''_1\ \ ,\ \ \hbox{and}\ \psi_ie_0=
e_1\psi_i\ \ ,\ \ i=2,3\ \ .\tag6.15$$ We put
$$(\mu_1, \mu_2, \mu_3)= (\psi''_1, \psi_2, \psi_3)\tag6.16$$
and again obtain a triple of isogenies from ${\Cal A}_0$ to ${\Cal
A}_1$. Since from ${\Cal A}_0, {\Cal A}_1$ and $\psi''_1$ we
recover ${\Cal A}$ with the action of $\underline\delta$ and
$\psi_1$ we see that again
$${\roman{Def}}({\Cal A}, \lambda_{\Cal A}, \iota_0; M)=
{\roman{Def}}({\Cal A}_0, {\Cal A}_1; \mu_1, \mu_2, \mu_3)\ \ .\tag6.17$$
The degree quadratic form on ${\roman{Hom}}({\Cal A}_0, {\Cal
A}_1)$ has matrix
$$T''=\diag(1, \varepsilon_2p^{a_2}, \varepsilon_3p^{a_3})\tag6.18$$
on the subspace of rank 3 spanned by the $\mu_i$'s.

We are now in a position to apply the results of Gross and Keating
\cite{\grosskeating} about the length of the artinian $W$-scheme
${\roman{Def}}({\Cal A}_0, {\Cal A}_1; \mu_1, \mu_2, \mu_3)$. By
their results this length only depends on the
$GL_3(\Z_p)$-equivalence class of the matrix $T''$. More precisely,
Proposition 5.4 of \cite{\grosskeating} implies the following
result.

\proclaim{Proposition 6.2} The length of the
local ring ${\Cal O}_{{\Cal Z}(T,\omega),\xi}$ is equal to
$e(\xi)=e_p(T)$ with $e_p(T)$ given as follows. Let $T$ be
$GL_3(\Z_p)$-equivalent to (6.5).
\hfill\break
(i) If $a_2$ is even,
$$e_p(T) = \sum_{i=0}^{a_2/2-1} (a_2+a_3-4i)p^i +\frac12 (a_3-a_2+1) p^{a_2/2}.$$
(ii) If $a_2$ is odd,
$$e_p(T) = \sum_{i=0}^{(a_2-1)/2} (a_2+a_3-4i)p^i.$$
\endproclaim

Note that the answer does not actually depend on the units $\e_1$, $\e_2$ and $\e_3$.

We conclude this section by indicating how the previous results can
be applied to the intersection problem of special cycles. The
calculus is the same as in the companion paper \cite{\krsiegel}, esp.\
section 3.

We fix integers $n_1,\ldots, n_r$ with $1\leq n_i\leq 3$ and with
$n_1+\ldots +n_r=3$. For each $i$ we choose $T_i\in
{\roman{Sym}}_{n_i}(\Q)_{>0}$ and a $K^p$-invariant open compact
subset $\omega_i\subset V({\Bbb A}_f^p)^{n_i}$. Let
$${\Cal Z} ={\Cal Z}(T_1,\omega_1)\times_{\Cal
M}\ldots\times_{\Cal M}{\Cal Z}(T_r,\omega_r,)\tag6.19$$
be the
fibre product of the corresponding special cycles. To each point
$\xi\in \Cal Z$ there is associated its {\it fundamental matrix}
$$T_{\xi}=Q((\j_1,\ldots, \j_n))\in{\roman{Sym}}_3(\Q)\ \ ,\tag6.20$$
where $(\j_1,\ldots, \j_r)$ is the 3-tuple of special endomorphisms
determined by the image of $\xi$ via its images under the projections ${\Cal Z}\to
{\Cal Z}(T_i,\omega_i)$. Since the function $\xi\mapsto T_{\xi}$ is
locally constant, we obtain a disjoint sum decomposition
$${\Cal Z}= \mathop{\coprod}\limits_{T} {\Cal Z}_T\ \ ,\tag6.21$$
where ${\Cal Z}_T$ is the union of those connected components of ${\Cal Z}$ where the value of
the fundamental matrix is equal to $T$.
Note that $T$ has the form
$$T=\pmatrix T_1&{*}&{\dots}&{*}\\{*}&T_2&{\dots}&{*}\\{\vdots}&
{\vdots}&\ddots&{\vdots}\\{*}&{*}&{\dots}&T_r\endpmatrix.\tag6.22$$
In fact we may identify
$${\Cal Z}_T={\Cal Z}(T,\omega_1\times \ldots \times \omega_r)\ \
.\tag6.23$$
Hence we may apply the previous results to obtain the following:

\proclaim{Corollary 6.3} Let $\xi\in{\Cal Z}={\Cal
Z}(T_1,\omega_1)\times_{\Cal M}\ldots\times_{\Cal M} {\Cal
Z}(T_r,\omega_r)$ with ${\roman{det}}(T_{\xi})\ne 0$. Then
$T_{\xi}\in {\roman{Sym}}_3(\Z_{(p)})_{>0}$ and $\xi$ lies over a
supersingular point of\hfill\break ${\Cal
M}\times_{{\roman{Spec}}\,
\Z_{(p)}}\! {\roman{Spec}}\,{\Bbb F}_p$. Furthermore, $\xi$ is an
isolated point of ${\Cal Z}$ if and only if $T_{\xi}\not\equiv 0\,
{\roman{mod}}\, p$. In this case the length of the local ring of
${\Cal Z}$ at $\xi$ is given by
$$e(\xi)=lg({\Cal O}_{{\Cal Z}, \xi})=e_p(T_{\xi})\ \ ,$$
with $e_p(T_{\xi})$ as in Proposition 6.2.
\endproclaim

Since the length dominates the local intersection multiplicity in
the sense of Serre we deduce from the formulas for $e_p(T)$ as in
\cite{\krsiegel}, Cor.\ 6.2.\ the following result.

\proclaim{Corollary 6.4}
The cycles ${\Cal Z}(T_1,\omega_1),\ldots, {\Cal Z}(T_r,\omega_r)$
intersect transversally at the point $\xi$ if and only if
${\roman{ord}}({\roman{det}}\, T_{\xi})=1$. In this case, the
schemes ${\Cal Z}(T_i,\omega_i)$ are regular at $\xi$ and their
tangent spaces give a direct sum decomposition of the tangent space
of ${\Cal M}$ at (the image of) $\xi$.
\endproclaim

\redefine\y{\text{\bf y}}

\subheading{\Sec7. Representation densities and Eisenstein series}

In this section we will consider the total contribution of isolated
intersection points to the intersection product of special cycles.
Since the development is analogous to the corresponding part of the
companion paper \cite{\krsiegel}, sections 7-10, we will be brief.

For special cycles ${\Cal Z}(T_1,\omega_1),\ldots, {\Cal
Z}(T_r,\omega_r)$ as at the end of the last section we define the
contribution of the isolated intersection points to the total
intersection number as
$$\langle {\Cal Z}(T_1,\omega_1),\ldots, {\Cal Z}(T_r,\omega_r)
\rangle_p^{\roman{proper}} =\sum_{\xi}e(\xi)\ \ .\tag7.1$$ Here
$\xi$ runs over the isolated points of ${\Cal
Z}(T_1,\omega_1)\times_{\Cal M}\ldots\times_{\Cal M} {\Cal
Z}(T_r,\omega_r)$. In the special case $r=1$ we have the cycle
${\Cal Z}(T,\omega)$ which lies over the supersingular locus of
${\Cal M}$ and consists of isolated points if and only if
$T\in{\roman{Sym}}_3(\Z_{(p)})$ is not divisible by $p$. In this
case we use the notation
$$\langle {\Cal Z}(T,\omega)\rangle_p= \sum_{\xi\in {\Cal
Z}(T,\omega)} e(\xi)\ \ .\tag7.2$$ By (6.21) and (6.23) we have
$$\langle {\Cal Z}(T_1,\omega_1),\ldots, {\Cal Z}(T_r,\omega_r)
\rangle_p^{\roman{proper}} =\sum_T \langle {\Cal Z}(T,\omega)\rangle_p\ \ ,\tag7.3$$
where the summation is over $T\in
{\roman{Sym}}_3(\Z_{(p)})_{>0}$ which have diagonal blocks
$T_1,\ldots, T_r$ as in (6.22) and for which $p\nmid T$. Hence it
suffices to determine the quantity (7.2).

By Proposition 6.2, the length of the local ring $\Cal O_{\Cal Z(T,\o), \xi}$
depends only on $T$ so that we may write
$$\langle{\Cal Z}(T,\omega)\rangle_p=
e_p(T)\cdot\vert {\Cal Z}(T,\omega)({\Bbb F})\vert\ \ .\tag7.4$$
Here $e_p(T)$ is given by Proposition~6.2. Thus it only remains to
determine $\vert {\Cal Z}(T,\omega)({\Bbb F})\vert$.

As at the end of section 4, let $B'$ be the definite quaternion
algebra over $\kay$ which is isomorphic to $B$ at all finite places
of $\kay$, and let $G'$ be the corresponding inner form of $G$,
defined by (0.2) with $B'$ in place of $B$. Then $G'(\Q_p)\simeq
G'_p$. Let $K'_p$ be the stabilizer of a superspecial lattice
$L_o\in X$, an Iwahori subgroup of $G'_p$ -- cf. (5.7) and
Corollary~5.2. Then the set of superspecial points in ${\Cal
M}^{ss}({\Bbb F})$ is in bijective correspondence with the double
coset space, comp.\ (4.12),
$$G'(\Q)\setminus \bigg(G'(\Q_p)/K'_p\times G({\Bbb A}_f^p)/K^p\bigg)\ \
.\tag7.5$$

Let $V'$ be the quadratic space over $\Q$ which is
positive definite, and with
$$V'(\Q_p)=V'_p\ \ \text{and}\ \ V'({\Bbb A}_f^p)= V({\Bbb A}_f^p)\
\ .\tag7.6$$
We note that $V'(\Q_p)$ and $V(\Q_p)$ have identical determinants
and opposite Hasse invariants. Also,  the quaternion algebra
$C^+(V')$ associated in section 0 to $V'$ is isomorphic to $B'$.
Indeed, it suffices to check this locally, using the fact that
$C^+(V')=B'_0\otimes_{\Q}\kay$ for the quaternion algebra $B_0'$
over $\Q$ given by Lemma~0.2. For the archimedean places the assertion
now follows from the positive definiteness of $V'$. At the place
$p$, both $B'$ and $C^+(V')$ split since $p$ is inert in $\kay$.
For the finite places not dividing $p$ the claim is obvious.

The
superspecial lattice $L_o$ (base point) then defines a lattice in
$V'_p$,
$$V'(\Z_p)= {\roman{End}}(L_o, F)\cap V'_p\ \ .\tag7.7$$
Let $\Omega'_T$ be the fiber over $T$ of the map defined by the
quadratic form on $V'$,
$$Q':{V'}^3\longrightarrow {\roman{Sym}}_3(\Q)\ \ .\tag7.8$$
Then $G'(\Q)$ acts transitively on $\Omega'_T(\Q)$ and the
stabilizer of a point $\y\in\Omega'_T(\Q)$ is $ Z'(\Q)$, the kernel
of the projection of $G'(\Q)$ to $SO(V')$, comp.\ \cite{\krsiegel},
Remark 7.4. The usual procedure therefore gives the following
expression for the cardinality of ${\Cal Z}(T,\omega)({\Bbb F})$,
cf.\ \cite{\krsiegel}, Lemma 7.1 and Theorem 7.2.

\proclaim{Proposition 7.1}
Let $K'=K'_p.K^p\subset G'({\Bbb A}_f)$. Let
$$\varphi_f^p= {\roman{char}}(\omega)\ ,\ \
\varphi'_p={\roman{char}} (V'(\Z_p))^3$$
and
$$\varphi'_f= \varphi'_p\cdot\varphi_f^p\in S(V'({\Bbb A}_f))^{K'}\
\ .$$
For arbitrary choices of a
base point $\y\in\Omega'_T(\Q)$ and of Haar measures on $G'({\Bbb
A}_f)$ and on $Z'({\Bbb A}_f)$), introduce the orbital integral
$$O_T(\varphi'_f)= \int\limits_{Z'({\Bbb A}_f)\setminus
G'({\Bbb A}_f)}\varphi'_f(g^{-1}\y)dg.$$
Then
$$\vert {\Cal Z}(T,\omega)({\Bbb F})\vert
={\roman{vol}}(K')^{-1}\cdot {\roman{vol}}(Z'(\Q)\setminus
Z'({\Bbb A}_f))\cdot O_T(\varphi'_f)\ \ .\tag7.9$$
\endproclaim
We therefore have obtained an explicit expression for (7.4).

We now want to compare this expression with the derivative of the $T$-th Fourier coefficient
of a certain Eisenstein series which is associated to our data as follows,
\cite{\annals}, \cite{\krsiegel}. (See Part I of \cite{\annals} for a more
extensive description of such incoherent Eisenstein series.)
We let
$$\lambda_p:S(V(\Q_p)^3)\longrightarrow I_3(0,
\chi_{V_p})\tag7.10$$
be the usual map into the induced representation of $Sp_{6,\Q_p}$
defined by $\lambda_p(\varphi_p)(g)=(\omega(g)\varphi_p)(0)$ where
$\omega=\omega_{\psi}$ denotes the action of $Sp_{6, \Q_p}$ on
$S(V(\Q_p)^3)$ via the local Weil representation defined by the
fixed additive character $\psi$. Let
$$\Phi_p=\lambda_p(\varphi_p)\ \ \text{with}\
\varphi_p={\roman{char}}\, V(\Z_p)^3\ \ .\tag7.11$$
Here $V(\Z_p)=\Lambda\otimes \Z_p$ in the notation of (2.1). We
define $\Phi'_p=\lambda'_p(\varphi'_p)$ with
$\varphi'_p={\roman{char}}\, V'(\Z_p)^3$ in the analogous way. We
complete $\Phi_p$ into an incoherent standard section
$$\Phi(s)= \Phi_{\infty}^2(s)\cdot \Phi_f^p(s)\cdot \Phi_p(s)\ \
,\tag7.12$$ where $\Phi^2_{\infty}(s)$ is associated in the usual
way to the Gaussian in $S(V'(\R)^3)$ and where
$\Phi_f^p(0)=\lambda_f^p(\varphi_f^p)$with $\varphi_f^p$ as in
Proposition 7.1. Then for $h\in Sp_{6,\R}$ and each $T\in
{\roman{Sym}}_3(\Q)_{>0}$ which is represented by $V({\Bbb A}_f^p)$
but not by $V(\Q_p)$ we have (comp.\ \cite{\krsiegel}, section 8),
$$\align
E'_T(h,0,\Phi) &
={\roman{vol}}(SO(V')(\R)\cdot pr(K'))\cdot W_T^2(h)\tag7.13\\
&\qquad
\times\frac{W'_{T,p}(e,0,\Phi_p)}{W_{T,p}(e,0,\Phi'_p)}
\cdot {\roman{vol}}(K')^{-1}{\roman{vol}} (Z(\Q)\setminus
Z({\Bbb A}_f))\cdot O_T(\varphi'_f)\ ,
\endalign$$
provided that $\omega$ is {\it locally centrally symmetric,}\/
i.e.\ invariant under the action of $\mu_2({\Bbb A}_f^p)$ so that
$\varphi'_f$ is locally even. Here $pr(K')$ denotes the image of
$K'$ under the projection map $pr: G'({\Bbb A}_f)\to SO(V')({\Bbb
A}_f)$. Also, $W_T^2(h)$ is an archimedean factor defined in
analogy with $W_T^{\frac{5}{2}}(h)$ in \cite{\krsiegel}, (8.24),
comp.\ \cite{\annals}, section 7.

We will see in a moment that the value of the Whittaker functional
in the denominator is non-zero.

Before stating the next proposition we recall that a nonsingular
$T\in {\roman{Sym}}_3(\Q_p)$ is represented by precisely one of the
quadratic spaces $V(\Q_p)$ and $V'_p=V'(\Q_p)$, \cite{\annals}, Proposition~1.3.

\proclaim{Proposition 7.2}
Let $T\in {\roman{Sym}}_3(\Q_p)$ with ${\roman{det}}(T)\ne 0$.
\roster
\item"{(i)}" If $W'_{T,p}(e,0,\Phi_p)\ne 0$, then $T\in
{\roman{Sym}}_3(\Z_p)$.
\item"{(ii)}" If $T\in {\roman{Sym}}_3(\Z_p)$ and if $T$ is
represented by $V(\Q_p)$ and $T\not\equiv 0\ {\roman{mod}}\, p$, then
$$W_{T,p}(e,0,\Phi'_p)= 2\cdot \gamma(V_p')\cdot p^{-4}(p^2-1)\ \ .$$
\item"{(iii)}" If $T\in {\roman{Sym}}_3(\Z_p)$ is represented by
$V'_p=V'(\Q_p)$ and $T\not\equiv 0\ {\roman{mod}}\, p$, then
$$W'_{T,p}(e, 0, \Phi_p)= {\roman{log}}\, p\cdot\gamma(V_p)\cdot (1+p^{-2})\cdot
(1-p^{-2})\cdot e_p(T)\ \ .$$
\endroster
The factors $\gamma(V_p)$ and $\gamma(V'_p)$ are eighth roots of unity given explicitly
by Proposition~A.4 of \cite{\annals}. Moreover,
$$\gamma(V_p)/\gamma(V'_p) = -1.$$
\endproclaim

{\bf Remark.} By  (1.16) of \cite{\annals},
the quadratic space $V'_p$ represents $T =\text{\rm diag}(\e_1 p^{a_1},\e_2 p^{a_2},\e_3 p^{a_3})$,
with $0\le a_1\le a_2\le a_3$ and $\e_i$ units
precisely when:
$$ -1= (-1)^{a_1+a_2+a_3} \chi(-1)^{a_1+a_2+a_3 + a_1a_2+a_2a_3+a_3a_1}
\chi(\e_1)^{a_2+a_3}\chi(\e_2)^{a_1+a_3}\chi(\e_3)^{a_1+a_2}.$$
Here $\chi(x)=(x,p)_p$ so that $\chi(\Delta)=-1$.

\demo{Proof}
We use the well-known relation, reviewed in the Appendix to \cite{\annals}, between the values of
$W_{T,p}(e,s,\Phi_p)$ at integer values of $s$ and representation
densities. For a suitable basis the quadratic form on $V(\Z_p)$ has
matrix
$$S= {\roman{diag}}(1,-1,1,-\Delta)\ \ .\tag7.14$$
For $r\ge 0$, let $S_r=S\perp H_{2r}$, where $H_{2r}$ denotes the
split quadratic form of rank $2r$ over $\Z_p$. Then
$W_{T,p}(e,r,\Phi_p)=0$ for $T\in {\roman{Sym}}_3(\Q_p)\setminus
{\roman{Sym}}_3(\Z_p)$ and
$$W_{T,p}(e,r,\Phi_p)= \gamma(V_p)\cdot \alpha_p(S_r, T)\ \
,\tag7.15$$ where $\alpha_p(S_r, T)$ is the classical
representation density of $T$ by $S_r$ and $\gamma(V_p)$ is the
factor appearing in \cite{\annals}, Cor.\ A.1.5. Here we have used the fact that $S$ is
unimodular. There is a rational function $A_{S,T}(X)$ such that
$$\alpha_p(S_r,T)= A_{S,T}(p^{-r})\ \ ,\tag7.16$$
i.e.\ $\alpha_p(S_r, T)$ is a rational function of $X=p^{-r}$.
Taking (7.15) and (7.16) together we obtain
$$W'_{T,p}(e,0,\Phi_p)=-{\roman{log}}\, p\cdot\gamma(V_p)\cdot
\frac{\partial}{\partial X} \{ A_{S,T}(X)\}\Big\vert_{X=1}\ \
.\tag7.17$$
To calculate $\alpha_p(S_r,T)$ we use Kitaoka's formulas in the
form given in \cite{\krsiegel}, section 10. We use repeatedly the
standard reduction formula
$$\alpha_p(N\perp M, N\perp L)= \alpha_p(N\perp M, N)\cdot
\alpha_p(M,L)\ \ ,\tag7.18$$
valid provided the quadratic form $N$ is unimodular. In our case
let
$$\tilde T={\roman{diag}}(\varepsilon_2p^{a_2},
\varepsilon_3p^{a_3}),\qquad \ T_{\varepsilon_1}=
{\roman{diag}}(\varepsilon_1, \varepsilon_2p^{a_2},
\varepsilon_3p^{a_3})\ \ .\tag7.19$$
Here $\varepsilon_i\in\Z_p^{\times}$, $a_i\ge 0$. Using the
reduction formula (7.18) we obtain
$$\alpha_p(S_r, T_{\varepsilon_1})=\alpha_p(S_r,
\varepsilon_1)\cdot\alpha_p(\tilde S_r, \tilde T)\ \ .$$
Here $\tilde S={\roman{diag}}(1,-1,-\varepsilon_1\Delta)$ and
$\tilde S_r=\tilde S\perp H_{2r}$ as before. Using the reduction
formula again we obtain
$$\alpha_p(H_{2r+4}, T_{\varepsilon_1\Delta})= \alpha_p(H_{2r+4},
\varepsilon_1\Delta)\cdot\alpha_p(\tilde S_r, \tilde T)\ \ ,$$
and hence
$$\alpha_p(S_r, T_{\varepsilon_1})=
\frac{\alpha_p(S_r,\varepsilon_1)}{\alpha_p(H_{2r+4},
\varepsilon_1\Delta)} \cdot\alpha_p(H_{2r+4},
T_{\varepsilon_1\Delta})\ \ ,\tag7.20$$ and
$$A_{S,T_{\varepsilon_1}}(X)=
\frac{A_{S,\varepsilon_1}(X)}{A_{H_4, \varepsilon_1\Delta}(X)}
\cdot A_{H_4, T_{\varepsilon_1\Delta}}(X)\ \ .\tag7.21$$
But we have \cite{\yang}, Theorem 3.2,
$$\alpha_p(H_{2r+4}, \varepsilon_1\Delta) = 1-p^{-2}X\ \
,\tag7.22$$ and
$$\alpha_p(S_r, \varepsilon_1)= 1+p^{-2} X\ \ ,\tag7.23$$
with $X=p^{-r}$. For the other factor on the right hand side we use
\cite{\krsiegel}, Proposition 10.3, which is, in turn, a reformulation of
a result of Kitaoka \cite{\kitaoka}:
\par\noindent
{\bf Case} $a_2$ {\bf even:} Then
$$\frac{\alpha_p(H_{2r+4},
T_{\varepsilon_1\Delta})}{(1-p^{-2}X)(1-p^{-2}X^2)}
=\sum_{\ell=0}^{\frac{a_2}{2}-1} p^{\ell}(X^{2\ell}+\chi(T)\cdot
X^{a_2+a_3-2\ell}) +p^{\frac{a_2}{2}}\cdot
X^{a_2}\cdot\sum_{j=0}^{a_3-a_2} (\sigma X)^j\tag7.24$$ where
$\sigma =\chi(-\varepsilon_1\Delta\varepsilon_2)$ and where
$\chi(T)$ is equal to 1 if $a_3$ is even and is equal to $\sigma$
if $a_3$ is odd. Hence (comp.\ loc.\ cit.) the value of (7.24) at
$X=1$ is equal to zero if and only if $\chi(T)=
\sigma=-1$. In this case the value of the derivative of
(7.24) at $X=1$ is equal to
$$-\sum_{\ell=0}^{\frac{a_2}{2}-1}p^{\ell}(a_2+a_3-4\ell)-
p^{\frac{a_2}{2}}\cdot \left( \frac{a_3-a_2+1}{2}\right)
=-e_p(T_{\varepsilon_1})\ \ ,$$
as a comparison with Proposition 6.2 shows.
\par\noindent
{\bf Case} $a_2$ {\bf odd:} In this case
$$\frac{\alpha_p(H_{2r+4},T_{\varepsilon_1\Delta})}{(1-p^{-2}X)
\cdot (1-p^{-2}X^2)}
=\sum_{\ell=0}^{(a_2-1)/2} p^{\ell} (X^{2\ell}+\chi(T)\cdot
X^{a_2+a_3-2\ell})\tag7.25$$ with $\chi(T)$ equal to
$\chi(-\varepsilon_1\Delta\cdot\varepsilon_3)$ if $a_3$ is even and
equal to $\chi(-\varepsilon_2\varepsilon_3)$ if $a_3$ is odd. This
expression vanishes at $X=1$ if and only if $\chi(T)=-1$ and in
this case the value of the derivative of (7.25) at $X=1$ is equal
to
$$-\sum_{\ell=0}^{(a_2-1)/2} p^{\ell}(a_2+a_3-4\ell)
=-e_p(T_{\varepsilon_1})\ \ .\tag7.26$$
In both cases we therefore obtain
$$\frac{\partial}{\partial X} (A_{H_4,
T_{\varepsilon_1\Delta}}(X))\Big\vert_{X=1} =-(1-p^{-2})^2\cdot
e_p(T_{\varepsilon_1})\ \ ,\tag7.27$$ provided that
$T_{\varepsilon_1}$ is not represented by $H_4$. Taking into
account (7.21)-(7.23) we obtain
$$\align
\frac{\partial}{\partial X}\{
A_{S,T_{\varepsilon_1}}(X)\}\Big\vert_{X=1} &
=-\frac{1+p^{-2}}{1-p^{-2}}\cdot (1-p^{-2})^2\cdot
e_p(T_{\varepsilon_1})\tag7.28\\ &
=-(1+p^{-2})(1-p^{-2})\cdot e_p(T_{\varepsilon_1})\ \ ,
\endalign$$
provided that $T_{\varepsilon_1}$ is not represented by $S$, i.e.\
by $V(\Q_p)$. Taking into account (7.17) we obtain the formula in
(iii).

We proceed in a similar way to prove (ii). For a suitable choice of
a basis for $V'(\Z_p)$ the quadratic form has matrix
$$S'={\roman{diag}}(1,-1,p,-p\Delta)\ \ .\tag7.29$$
Using the reduction formula (7.18) twice we obtain
$$\alpha_p(S', T_{\varepsilon_1})= \frac{\alpha_p(S',
\varepsilon_1)}{\alpha_p({\tilde S}', -\varepsilon_1\Delta)}
\cdot\alpha_p ({\tilde S}', T_{-\varepsilon_1\Delta})\ \
,\tag7.30$$ where
$${\tilde S}'={\roman{diag}}(1,\Delta,p,-p\Delta)\tag7.31$$
is the quadratic form given by the norm form on the maximal order
of the quaternion division algebra over $\Q_p$. Appealing to
\cite{\yang}, Theorem 3.2 we have
$$\alpha_p(S', \varepsilon_1)= 1-p^{-1},\qquad  \alpha_p({\tilde S}',
-\varepsilon_1\Delta)= 1+\chi(-1)p^{-1}\tag7.32$$
and by \cite{\grosskeating}, Proposition~6.10,          
$$\alpha_p({\tilde S}', T_{-\varepsilon_1\Delta})=
2(1+\chi(-1)p^{-1})\cdot (p+1)\ \ .\tag7.33$$
Inserting (7.32) and
(7.33) into (7.30) we therefore obtain
$$
\alpha_p(S', T_{\varepsilon_1})
=2p^{-1}\cdot (p^2-1)\ \ .\tag7.34
$$ We have the relation
$$\align
W_T(e,0,\Phi'_p) &
=\gamma(V'_p)\cdot\vert S'\vert^{\frac{3}{2}}\cdot\alpha_p(S',
T)\tag7.35\\ &
=2\cdot\gamma(V'_p)\cdot p^{-4}\cdot (p^2-1)\ \ ,
\endalign$$
which proves (ii).
\qed\enddemo
We now plug in the expressions obtained in Proposition 7.2 into
(7.13) and obtain
$$\align
E'_T(h,0,\Phi) &
=-{\roman{vol}}(SO(V')(\R)\cdot pr(K'))\cdot
W_T^2(h)\cdot\tag7.36\\ &
\qquad\cdot{\roman{log}}\,
p\cdot\frac{(1+p^{-2})(1-p^{-2})}{2\cdot (p^2-1)} \cdot p^4\cdot
e_p(T)\cdot{\roman{vol}}(K')^{-1}\cdot\\ &\qquad
\cdot{\roman{vol}}(Z(\Q)\setminus Z({\Bbb A}_f))\cdot
O_T(\varphi'_f)\\ &
=-{\roman{vol}}(SO(V')(\R)\cdot pr(K'))\cdot
W_T^2(h)\cdot\frac{1}{2}{\roman{log}}\, p\cdot (p^2+1)\cdot\langle
{\Cal Z}(T,\omega)\rangle_p
\endalign$$
Taking into account the fact that the index of the Iwahori subgroup $K'_p$
in the maximal compact subgroup $K_p$ is equal to $p^2+1$, we
finally obtain the following theorem.

\proclaim{Theorem 7.3}
Let $T\in{\roman{Sym}}_3(\Q)_{>0}$ be represented by $V({\Bbb
A}_f^p)$ but not by $V(\Q_p)$. Also assume that $\omega$ is locally
centrally symmetric.
\roster
\item"{(i)}" If $T\not\in {\roman{Sym}}_3(\Z_{(p)})$, then ${\Cal
Z}(T,\omega)=\emptyset$ and $E'_T(h, 0, \Phi)=0$.
\item"{(ii)}" If $T\in{\roman{Sym}}_3(\Z_{(p)})$ is not divisible by
$p$, then ${\Cal Z}(T,\omega)$ is zero-dimensional and
$$E'_T(h,0,\Phi)=-\frac{1}{2}{\roman{vol}}(SO(V')(\R))\cdot
W_T^2(h)\cdot{\roman{vol}}(pr(K))\cdot{\roman{log}}\, p\cdot\langle
{\Cal Z}(T,\omega)\rangle_p\ .$$
\endroster
\endproclaim

To apply this theorem to the intersection problem of special cycles
consider the situation of the end of section 6. Let
$$W=W_1+\ldots + W_r\tag7.37$$
be the decomposition of the standard 6-dimensional symplectic space
into symplectic subspaces of dimension $2n_i$ compatible with the
fixed symplectic basis, and let
$$\iota: H_{1,{\Bbb A}}\times\ldots\times H_{r,{\Bbb
A}}\longrightarrow Sp_{6,{\Bbb A}}\tag7.38$$ be the corresponding
homomorphism of metaplectic groups, covering the embedding
$$\iota: Sp(W_{1,{\Bbb A}})\times\ldots\times Sp(W_{r, {\Bbb
A}})\hookrightarrow Sp_{6,{\Bbb A}}\ \ .\tag7.39$$ Restricting to
the archimedean place, for $(h_1,\ldots, h_r)\in H_{1,\R}\times
\ldots\times H_{r,\R}$ we have
$$W_T^2(\iota(h_1,\ldots, h_r))=W_{T_1}^2 (h_1)\cdot\ldots\cdot
W^2_{T_r}(h_r)\ \ ,\tag7.40$$ where $T$ has diagonal blocks
$T_1,\ldots, T_r$.

\proclaim{Corollary 7.4}
We have the following identity, provided that
$\omega_1,\ldots,\omega_r$ are locally centrally symmetric.
$$\align
\sum_TE'_T(\iota(h_1,\ldots, h_r),0,\Phi)
&
=-\frac{1}{2}{\roman{vol}}(SO(V')(\R))\cdot
W^2_{T_1}(h_1)\ldots W^2_{T_r}(h_1)\cdot\\ &\qquad
\cdot{\roman{vol}}(K)\cdot{\roman{log}}\, p\cdot\langle {\Cal
Z}(T_1,\omega_1), \ldots, {\Cal
Z}(T_r,\omega_r)\rangle_p^{\roman{proper}}
\endalign$$
Here the summation runs over $T\in {\roman{Sym}}_3(\Z_{(p)})_{>0}$
with diagonal blocks $T_1,\ldots, T_r$, and such that (i) $T$ is represented
by $V({\Bbb A}_f^p)$ but not by $V(\Q_p)$ and (ii)
$T\not\equiv 0\ {\roman{mod}}\, p$.
Also $\Phi$ is determined as in (7.12) with $\omega
=\omega_1\times \ldots\times \omega_r$.
\endproclaim

We note that as in \cite{\annals} the left side of this expression
is part of the Fourier coefficient corresponding to $(T_1,\ldots,
T_r)$ of the pullback of $E'(g,0,\Phi)$ to $H_{1,{\Bbb A}}\times
\ldots\times H_{r, {\Bbb A}}$.

\subheading{\Sec8. Components of $\CZ(T,\o)\cap{\Cal M}^{ss}$}

When $T\in {\roman{Sym}}_n(\Z_{(p)})_{>0}$ is divisible by $p$,
i.e., when $\bar{T}=0$, Theorem~6.1 says that
$\big(\text{image}(\CZ(T,\o)\big)\cap\mss(\F)$ is a union of certain components
$\Bbb P_{L_i}$ of $\mss(\F)$. In this section, we describe the
configurations of components which occur, in terms of the matrix
$T$. We again consider the isocrystal ${\Cal L}$ associated to a
fixed base point $\xi=(A,\lambda,\iota,\overline{\eta}^p)\in\mss(\F)$.

Recall, from (5.1), that $V'_p =\{ j\in \End(\L,F)^{(1)}\mud
j^*=j\}$ is the space of special endomorphisms of the graded,
polarized isocrystal $(\L,<\ ,\ >,\iota)$. Suppose that $T\in p
{\roman{Sym}}_n(\Z_{(p)})$ and that $\j\in (V'_p)^n$ with
$Q'(\j)=T$. We want to determine the set
$$X(\j):=\{\  L\in X\mud  \j\in \big(N_L\big)^n \ \}.\tag 8.1$$
Obviously this only depends on the $\Z_p$-submodule of $V'_p$
spanned by the components of $\j$. Thus we can work with a $\Z_p$-basis for this submodule
for which the restriction of the quadratic form $Q'$ is diagonal.

We begin by
determining, for $j\in V'_p$, the set
$$X(j)=\{ L\in X\mud\ j\in N_L\}\ \ .\tag8.2$$
We will always assume that $j$ is not isotropic, i.e.\ $j^2\ne 0$.

Recall from the end of section~4 the set
$$Y_i=\{\ L_i\subset \L_i\mud\   F^2L_i=pL_i\ \},\tag 8.3$$
and let
$$Y_i(j)=\{\ L_i\in Y_i\mud\  j(L_i)\subset FL_i
 \}.\tag8.4$$
Then, combining Proposition~5.5 and Proposition~4.4,
we have:
\proclaim{Proposition 8.1} There is a natural surjection
$$\coprod_{i=0,1}\coprod_{L_i\in Y_i(j)} \Bbb P_{L_i}(\F)\lra X(j).$$
\endproclaim

Thus, we must give a description of the sets $Y_i(j)$, for $i=0,1$.
These sets are invariant under homothety and $FY_0(j)=Y_1(j)$,
since $j$ commutes with $F$. Thus, it suffices to consider
$Y_0(j)$.

Let
$$U=U_0=\big(\L_0\big)^{p^{-1}F^2},\tag8.5$$
be the $2$-dimensional $\Q_{p^2}$-vector space of fixed points of
the $\s^2$-linear endomorphism $p^{-1}F^2$. This space comes
equipped with a nondegenerate alternating form $<\ ,\ >$ valued in
$\qps$. Moreover, $G'_p$ preserves $U$ and restriction yields an
isomorphism
$$G'_p\simeq\{g\in GL(U)\mud \det(g)\in \Q_p^\times\},$$
cf.\ (5.7).

The set $Y_0$ is then precisely the set of $W$-lattices in
$\L_0=U\tt_\qps\K$ of the form $L_0=\LL\tt_{\Z_{p^2}}W$, where
$\LL= (L_0)^{p^{-1}F^2}$ is a $\zps$-lattice in $U$. Let $\CB_0$ be
the set of $\zps$-lattices $\Lambda$ in $U$ up to homothety, i.e.,
the set of vertices of the building $\CB$ for the group $PGL(U)$.
The natural surjection
$$pr:Y_0\rightarrow\CB_0\tag 8.6$$
is given by $L_0=\LL\tt W\mapsto [\LL]$.

Let $j\in V'_p$. The endomorphism $F^{-1}j$ of $\L$ has degree $0$,
is $\s^{-1}$-linear, and commutes with $F$. It therefore induces a
$\sigma$-linear endomorphism of $U=\L_0^{p^{-1}F^2}$,
$$\b = F^{-1}j|_U.\tag 8.7$$
Obviously $\b$ conversely determines $j$. The condition $j^*=j$ is
equivalent to $\b=\b^*$ where $*$ is defined by
$$<\b x, y > \ = \ <x,\b^* y>^\s.\tag 8.8$$
Indeed, for $x$ and $y\in U$,
$$\align
<\b x,y>&\ = \ <F^{-1}j x,y> \ = \  p^{-1}<j x,Fy>^{\s}\\
\nass
{}&\ = \ <x,j p^{-1}Fy>^\s \ = \  <x,F^{-1}j y>^\s \ = \  <x, \b y>^\s.\endalign$$
We will refer to the elements of the space
$$V''_p:=\{ \beta\in \End_{\Q_p}(U)\mud\  \text{$\b$ is $\s$\snug-linear,
and $\b^*=\b$}\ \}\tag8.9$$ as {\it special endomorphisms}\/ of
$(U,<\ ,\
>)$. Of course, as $\Q_p$-vector spaces, $V'_p\simeq V''_p$ via the
map (8.7). Define a quadratic form $Q''$ on $V''_p$ by
$$\b^2 =
Q''(\b)\cdot 1_U,\tag8.10$$
and note that
$$Q'(j)=p Q''(\b),\tag8.11$$
since
$$\align
Q''(\b)<x,y>^\s&\ = \ <\b x,\b y> \ = \  <F^{-1}j x,F^{-1} j y>\\
{}& \ = \  p^{-1}< j x,j y>^\s
\ = \  p^{-1}Q'(j)<x,y>^\s.\endalign$$
If $e_1, e_2$ is a $\qps$-basis for $U$ for which $<e_1,e_2>=1$,
then a $\s$-linear endomorphism can be written in the form
$$\b=\pmatrix a&b\\c&d\endpmatrix \s,$$
and
$$\b^* = \pmatrix d^\s&-b^\s\\-c^\s&a^\s\endpmatrix \s,$$
where $\sigma^2=1$ on $U$. Thus $\b$ is special if it has the form
$$\b=\pmatrix a&b\\c&a^\s\endpmatrix \s,\tag8.12$$
for $a$, $b$, and $c\in\qps$ with $b^\s=-b$ and $c^\s=-c$.
Then $Q''(\b)=aa^\s-bc$.

If $\b\in V''_p$ and $j\in V_p'$ are related by (8.7), then
$L_0=\LL\tt_{\zps}W\in Y_0(j)$ if and only if $\b(\LL)\subset \LL$.
Our description of $Y_0(j)$ depends on the following simple
observation which goes back to Kottwitz and Tate (comp.\
\cite{\krdrin}, Lemma 2.4). A special endomorphism $\beta$ of $U$
induces an involution on ${\Cal B}$ (since we are always assuming
that $\beta^2=Q''(\beta)\cdot 1_U$ with $Q''(\beta)\ne 0$). Let
${\Cal B}^{\beta}$ denote the fixed point set of $\beta$. Also let
$d(x,y)$ denote the distance in the building.

\proclaim{Lemma 8.2}
If $\b$ is a special endomorphism of $U$, then
$$\align
\b(\LL)\subset \LL
&
\quad \iff \quad d([\beta(\LL)],[\LL]) \le
\ord_p(Q''(\b))=\ord_p(\det\beta)\\
&
\quad\iff\quad d([\Lambda], {\Cal B}^{\beta})\leq
\frac{1}{2}\cdot\ord_p(\Q''(\beta))\ \ .
\endalign$$
\endproclaim

For a special endomorphism $\b$ of $U$, let
$$\CT(\b)=\{\ x\in \CB\ \mud \ d(x,\CB^{\b})\le \frac{1}{2}\ord_p(\det\beta)\ \}
\tag 8.13$$
be the closed tube of radius $\frac{1}{2}\cdot\ord_p(\det\beta)$ around
the fixed point set. Let $\CT(\beta)_0$ be the set of vertices in
$\CT(\beta)$.

\proclaim{Corollary 8.3} For a special endomorphism $\b$ of $U$
corresponding to a special endomorphism $j\in V'_p$, via (8.7),
$$Y_0(j) = \pr^{-1}(\CT(\b)_0)\simeq \Z\times\CT(\beta)_0,
$$
where $\pr$ is the projection of (8.6). \qed
\endproclaim

Returning to our $n$-tuple of special endomorphisms $\j$, let
$\bob=(\b_1,\dots,\b_n)$ be the corresponding tuple of special
endomorphisms of $U$.
As noted above, the sets $X(\j)$ and $Y_0(\j)$ depend only on the $\Z_p$-span
of the components of $\j$. Thus, by changing the $\Z_p$-basis if necessary, we may assume that
the matrix $T$ is of the form
$$T=Q'(\j)= \diag(\varepsilon_1 p^{a_1},\dots,\varepsilon_n
p^{a_n})\tag8.14$$
with $1\le a_1\le\dots \le a_n$. It follows by (8.11) that
$$Q''(\bob) = \diag(\varepsilon_1 p^{r_1},\dots,\varepsilon_n
p^{r_n})=p^{-1}T,\tag8.15$$
with $r_i=a_i-1$.
\proclaim{Corollary 8.4} Let
$$\CT(\bob):=\CT(\b_1)\cap \dots\cap \CT(\b_n).$$
Then
$$Y_0(\j) = \pr^{-1}( \CT(\bob)_0)\simeq\Z\times{\Cal T}(\bob)_0.
\qquad\qed$$
\endproclaim

Let us now return to $X(j)$. Suppose we are given a path in the
building $\CB$ with consecutive vertices $x_0,\ x_1,\ \dots,x_r,\
\dots$. Let $\LL_0\supset
\LL_1\supset\dots\supset \LL_r\supset\dots$ be a corresponding
sequence of lattices in $U$ with $(\LL_r:\LL_{r+1})=1$. Construct a
sequence of lattices
$$L_0=\LL_0\tt W,\ L_1=F^{-1}(\LL_1\tt W),\  \dots L_r=F^{-r}(\LL_r\tt W)
\dots,\tag8.16$$
where $L_r\in Y_0$ if $r$ is even and $L_r\in Y_1$ if $r$ is odd.
Associated to this sequence of lattices is a `chain' of $\Bbb
P^1$'s
$$\Bbb P_{L_0},\ \Bbb P_{L_1},\ \dots,\Bbb P_{L_r},\ \dots\tag8.17$$
in $X$. These $\Bbb P^1$'s cross at the sequence of superspecial lattices:
$$\align
L_0\oplus L_1 &= \Bbb P_{L_0}\cap\Bbb P_{L_1},\\
\nass
L_2\oplus L_1 &= \Bbb P_{L_2}\cap\Bbb P_{L_1},\\
\nass
L_2\oplus L_3 &= \Bbb P_{L_2}\cap\Bbb P_{L_3},\tag8.18\\
\nass
\dots\\
\nass
L_{2t}\oplus L_{2t\pm1}&=\Bbb P_{L_{2t}}\cap\Bbb P_{L_{2t\pm1}},\\
\nass
\dots,
\endalign
$$
which can be viewed as indexed by the edges in the path. Moreover,
if the path lies in $\CT(\b)$, then the chain of $\Bbb P^1$'s lies
in $X(j)$.\hfill\break A given path also gives rise to additional
chains of $\Bbb P^1$'s obtained from the first by applying powers
of $F$. These chains are all disjoint (consider the index in a
fixed lattice in $U$). More generally, any connected subset of
$\CB$ can be viewed as the dual graph to a connected curve of $\Bbb
P^1$'s in $X$, uniquely determined up to the action of $F$.


Applying these considerations to $\CT(\bob)$, we obtain the following.
\proclaim{Proposition 8.5} The set $X(\j)$ is a disjoint union of copies, indexed
by $\Z$, of a connected union of $\Bbb P^1$'s. The dual graph of a
connected component is naturally isomorphic to $\CT(\bob)$, where
$\bob$ and $\j$ are related by (8.7). In particular, the number of
$\Bbb P^1$'s in a connected component is $|\CT(\bob)_0|$. The
number of crossing points in a connected component is
$\vert\CT(\bob)_1|$, where $\CT(\bob)_1$ denotes the set of edges
contained in $\CT(\bob)$.
\endproclaim

It thus remains to obtain a better understanding of the sets
$\CT(\bob)$. To this end we collect some facts on the action of
special endomorphisms on the building $\CB$. We first note that the
building $\CB$ has the property that every pair of points $x$ and
$y$ is joined by a unique geodesic and that every path from $x$ to
$y$ contains that geodesic.

\proclaim{Lemma 8.6} (i) For a special endomorphism $\b$ of $U$, the
fixed point set $\CB^{\b}$ is connected.\hfill\break (ii) If $x\in
\CB$, then the midpoint of the unique geodesic in $\CB$ joining $x$
and $\b (x)$ is fixed by $\b$ and is the point of $\CB^{\b}$
nearest to $x$. If $x\in \CT(\b)$, then this geodesic lies in
$\CT(\b)$.\hfill\break (iii) The set $\CT(\b)$ contains the
geodesic joining any two of its points.\hfill\break (iv) If
$\b_1,\dots,\b_m$ are special endomorphisms of $U$, then
$\CT(\b_1)\cap\dots\cap\CT(\b_m)$ is connected.
\endproclaim
\demo{Proof} If $x$ and $y$ are fixed points of $\b$ in $\CB$, then the
unique geodesic in $\CB$ joining $x$ and $y$ is also fixed. This
proves (i) since non-emptyness follows from (ii). The assertion
(ii) follows from the fact that $\beta$ induces an involution of
$\CB$ (comp.\ \cite{\krdrin}, Lemma 2.4.) and was in fact used in
the proof of Lemma 8.2 above.  The set $\CT(\b)$ is connected since
each point of it is connected to $\CB^{\b}$ by a geodesic in
$\CT(\b)$, and $\CB^{\b}$ is connected. Since the geodesic joining
any pair of points is a subset of any path joining them, (iii)
follows. Finally, (iv) is immediate from (iii).
\qed\enddemo

\proclaim{Lemma 8.7} Let $\b$ be a special endomorphism of $U$.\hfill\break
(i) If $\ord_p(Q''(\b))$ is even, then $\CB^{\b}\simeq
\CB(PGL_2(\Q_p))$.\hfill\break (ii) If $\ord_p(Q''(\b))$ is odd,
then $\CB^{\b}$ is a single point, which is a midpoint of an edge.
\endproclaim
\demo{Proof} The centralizer $G_\b$ in $GL_2(\qps)$ of the $\s$-linear endomorphism $\b$
is a $\Q_p$-form of $GL_2$, and obviously ${\Cal
B}(PG_{\beta})\subset\CB^{\beta}$. If $\ord_p(Q''(\b))$ is even,
then $G_\b\simeq GL_2(\Q_p)$, while, if $\ord_p(Q''(\b))$ is odd,
then $G_\b\simeq \Bbb B_p^\times$. In the latter case
$\CB(PG_{\beta})=\CB^{\beta}$ is a midpoint. Let us show that also
in the first case $\CB(PG_{\beta})=\CB^{\beta}$. After correcting
$\beta$ by a power of $p$ we may assume that $\beta^2=1_U$. If
$\Lambda$ is a representative of a vertex of $\CB^{\beta}$, then
$\beta$ induces a $\sigma$-linear automorphism
$\overline{\beta}:\Lambda/p\Lambda\to \Lambda/p\Lambda$. But
$\overline{\beta}$ defines an ${\Bbb F}_p$-rational structure on
${\Bbb P}(\Lambda/p\Lambda)$ which cannot have more than $p+1$
${\Bbb F}_p$-rational points. The assertion follows since
$\CB^{\beta}$ is connected, cf.\ Lemma 8.6, (i).
\qed\enddemo

\proclaim{Lemma 8.8} Suppose that $\b_1,\dots,\b_n$ is a collection of
special endomorphisms which anticommute, i.e., such that
$\b_i\b_j+\b_j\b_i=0$ for $i\ne j$. The set of common fixed points
$\CB^{\b_1}\cap\dots\cap\CB^{\b_n}$ is nonempty and connected.
\endproclaim
\demo{Proof}
The set $\CB^{\b_i}$ is preserved by $\b_j$ for each $j$. If
$x_1\in \CB^{\b_1}$, then the midpoint $x_{12}$ of the geodesic
from $x_1$ to $\b_2 (x_1)$ lies in $\CB^{\b_1}\cap\CB^{\b_2}$, by
(i) and (ii) of Lemma~8.6. The midpoint $x_{123}$ of the geodesic
joining $x_{12}$ and $\b_3 (x_{12})$ lies in
$\CB^{\b_1}\cap\CB^{\b_2}\cap
\CB^{\b_3}$, etc.
\qed\enddemo

\proclaim{Lemma 8.9} For a collection of special endomorphism $\b_i$, as in
Lemma~8.8, and a point $x\in \CB$, let $x_i\in \CB^{\b_i}$ be the
point closest to $x$, and let $d_i=d(x,x_i)=d(x,\CB^{\b_i})$.
Assume that $d_1\le d_2\le\dots
\le d_n$. Then
$$\align
x_1&\in \CB^{\b_1}\\ x_2&\in \CB^{\b_1}\cap \CB^{\b_2}\\ x_3&\in
\CB^{\b_1}\cap \CB^{\b_2}\cap\CB^{\b_3}\\
\nass
&\dots\\
\nass
x_n&\in \CB^{\b_1}\cap\dots\cap\CB^{\b_n}.\endalign
$$
Moreover, if, for some $i$, $d_i=d_{i+1}$, then $x_i=x_{i+1}$. Finally,
for any $i$ with $1\le i\le n$,
the geodesic from $x$ to $x_i$ is
$$[x,x_i]=[x,x_1][x_1,x_2]\dots [x_{i-1},x_i].\qquad\qed$$
\endproclaim

Now return to the set $\CT(\bob)$ of Corollary 8.4, where the
$\b_i$'s anticommute and where $0\le r_1\le \dots\le r_n$, as in
(8.14).

{\it First suppose that $r_1$ is odd.}\/ Then $\CB^{\b_1}$ is the
midpoint of an edge, and $\CB^{\b_1}\subset \CB^{\b_i}$ for all
$i$. Since $r_1\le r_i$, for all $i$,
$$\CT(\b_1)\subset \CT(\b_i)$$
for all $i$, and hence $\CT(\bob)$ is simply the ball of radius
$r_1/2$ centered at the point $\CB^{\b_1}$.

{\it Next suppose that $r_1=2t$ is even.}\/ There is a unit
$u\in\zpsx$ such that $\b_0:=u^{-1} p^{-t}\b_1$ satisfies
$\b_0^2=1$. Note however that $\b_0$ is no longer a special
endomorphism. We may then write $U=U^0\tt_{\Q_p}\qps$, where $U^0$
is the fixed point set of $\b_0$. This gives a natural isomorphism
$\CB^{\b_1}=\CB^{\b_0}\simeq
\CB(PGL(U^0))$, which is given on vertices by $[\LL]\mapsto
[\LL^0]$, for $\LL=\LL^0\tt_{\Z_p}\zps$, where $\LL^0$ is a
$\Z_p$-lattice in $U^0$.

We fix a basis $e_1,e_2$ for $U^0$ and write
$\beta_1=up^t\cdot\sigma=up^t\cdot\beta_0$. For $i\geq 2$, define
the matrix $\gamma_i\in GL_2(\Q_{p^2})$ by $\beta_i=\gamma_i
u\delta\sigma$. Then the relation $\beta_1\beta_i=-\beta_i\beta_1$
is equivalent to $\gamma_i\in GL_2(\Q_p)$.
Next we claim that $\gamma_i^2$ is a scalar matrix. Indeed, we have
$$\b_i^2 = (\gamma_i u\delta \s)(\gamma_i u\delta \s).$$
Since $\b_1\b_i=-\b_i\b_1$, we have
$$u\s (\gamma_i u \delta \s) = -(\gamma_i u \delta \s)u\s,  \qquad i.e., \qquad
u \gamma_i u^\s \delta^\s = - \gamma_i u \delta u^\s.$$
But $\delta = -\delta^\s\cdot(\text{\rm scalar})$, and hence
$$u \gamma_i = \gamma_i u.$$
Therefore,
$$\b_i^2 = \gamma_i^2 u \delta u^\s \delta^\s.$$
But $u u^\s = p^{-2t}\b_i^2$ is a scalar matrix, hence, since
$$\b_i^2 = \gamma_i^2\cdot u u^\s (-\Delta)$$
is a scalar matrix, so is $\gamma_i^2$, as claimed.

Moreover,
$$\gamma_i\gamma_j+\gamma_j\gamma_i =0,\ \text{ for $i\ne j$, and }
\qquad \gamma_i^2 =
- \Delta^{-1}\varepsilon_1^{-1}\varepsilon_i p^{r_i}.$$
Under the isomorphism $\CB^{\b_1}\simeq \CB(PGL(U_0))=:\CB^0$ we
have
$$\CB^{\b_1}\cap\CB^{\b_i} \simeq \CB^{0,\gamma_i},$$
where $\CB^{0,\gamma_i}$ is the fixed point set of $\gamma_i$ in
$\CB^0$.

\proclaim{Lemma 8.10} Let $\gamma\in GL_2(\Q_p)$ with $\tr(\gamma)=0$.
Then
$$
\CB^{0,\gamma} =\cases \text{ an apartment} &\text{ if $\gamma^2\in \Q_p^{\times,2}$,}\\
\text{ a vertex }&\text{ if $\gamma^2\notin \Q_p^{\times,2}$, and
$\ord_p(\gamma^2)$ is even,}\\
\text{ a midpoint}&\text{ if $\gamma^2\notin \Q_p^{\times,2}$, and
$\ord_p(\gamma^2)$ is odd.}
\endcases$$
Furthermore $\gamma$ generates a split, unramified elliptic, or
ramified elliptic Cartan subgroup of $GL_2(\Q_p)$ in the three
cases respectively.
\endproclaim
\demo{Proof} The matrix $\gamma$ has two eigenlines in $\Bbb P^1(
\overline{\Q}_p)$. In the first
case, these lie in $\Bbb P^1(\Q_p)$, the boundary of $\CB^0$, and
the (infinite) geodesic joining them is the apartment fixed by
$\gamma$. The other two cases are clear.
\qed\enddemo

Corresponding to the cases we call $\beta_i$ $(i\geq 2)$ split,
unramified elliptic or ramified elliptic.

\proclaim{Lemma 8.11}
Assume that $r_1$ is even. If $\beta_i$ is ramified elliptic (resp.\
unramified elliptic) for some $i\geq 2$, then $\beta_j$ is not
unramified elliptic ( resp.\ ramified elliptic) for any $j\geq 2$.
\endproclaim

\demo{Proof}
Otherwise one fixed point set would be a midpoint of an edge and
the other a vertex in $\CB^0$, contradicting Lemma 8.8.
\qed
\enddemo

\proclaim{Lemma 8.12} If $n=3$ or $4$, then $\CB^{\bob}:=
\CB^{\b_1}\cap\dots\cap\CB^{\b_n}$
is a point. This point is a midpoint if at least one $r_i$ is odd
and is a vertex if all $r_i$ are even.
\endproclaim
\demo{Proof} This is immediate if any of the $r_i$'s is odd, so suppose that each $r_i$
is even. We may assume that $\CB^{\b_1}\cap\CB^{\b_2}\simeq
\CB^{0,\gamma_2}$ is an apartment associated to the two eigenlines of
$\gamma_2$ in $U^0$. Since $\gamma_2\gamma_3=-\gamma_3\gamma_2$,
these eigenlines are switched by $\gamma_3$, so that $\gamma_3$
preserves $\CB^{0,\gamma_2}$ and acts by a reflection through some
unique fixed point. This point is a vertex $\CB^{0,\gamma_2}$,
since $\ord_p(\gamma_3^2)$ is even. If $n=4$, $\gamma_4$ also fixes
this vertex.
\qed\enddemo

\proclaim{Proposition 8.13}
If $n=3$ or 4, then $\CT(\bob)_0$ consists of a single vertex if
and only if $r_1=r_2=0$ and if $r_i$ is even for $i\geq 3$ and
furthermore $r_3=0$ if $\chi(-\varepsilon_1\varepsilon_2)=-1$.
Equivalently, by Proposition 8.5, each connected component of
${\Cal Z}(T,\omega)\cap{\Cal M}^{ss}$ is irreducible of dimension 1
if and only if $T$ is $GL_3(\Z_p)$-equivalent to ${\roman{diag}}(p,
\varepsilon_2p, \varepsilon_3p^{a_3})$ with $a_3$ odd and with $a_3=1$ if
$\chi(-\varepsilon_2)=-1$ resp.\ $GL_4(\Z_p)$-equivalent to
${\roman{diag}}(p, \varepsilon_2p,
\varepsilon_3p^{a_3}, \varepsilon_4p^{a_4})$ with $a_3$ and $a_4$
odd and $a_3=1$ if $\chi(-\varepsilon_2)=-1$.
\endproclaim

\demo{Proof}
If there exists $i\geq 1$ with $r_i$ odd, then $\CB^{\bob}$
consists of the midpoint of an edge and both vertices of this edge
belong to $\CT(\bob)_0$. Hence, if $\CT(\bob)_0$ consists of a
single vertex, all $r_i$ have to be even. The unique vertex $x_0$
in $\CB^{\bob}$ belongs to $\CT(\bob)_0$. If $r_1\geq 2$, then any
vertex with distance 1 from $x_0$ in $\CB$ also belongs to
$\CT(\bob)_0$, hence $r_1=0$. A similar argument applied to
$\CB^0=\CB^{\beta_1}$ shows that $r_2=0$. However, in this case
$\beta_2$ is unramified elliptic iff $\chi(-\Delta
\varepsilon_1\varepsilon_2)=-1$. For such a $\b_2$,  $\CB^{\beta_1}\cap
\CB^{\beta_2}=\CB^{0,\gamma_2}$ consists of the single vertex
$x_0$ and $(\CT(\beta_1)\cap \CT(\beta_2))_0= \{ x_0\}$. If
$\chi(-\Delta\varepsilon_1\varepsilon_2)=1$, then
$\CB^{\beta_1}\cap \CB^{\beta_2}=\CB^{0,\gamma_2}$ is an apartment
in $\CB^0$. To exclude the vertices on this apartment from
$\CT(\bob)_0$ we must have $r_3=0$. The result follows.
\qed
\enddemo

It seems to us that the cases enumerated in Proposition 8.13 are
the simplest to consider when one wants to determine the
contribution of ${\Cal Z}(T,\omega)$ to the intersection product of
special cycles ${\Cal Z}(T_1,\omega_1),\ldots, {\Cal
Z}(T_r,\omega_r)$ in the case of excess intersection. The hope
would be to obtain a result similar to Corollary 7.4.

In general, the determination of the number of irreducible
components within one connected component of ${\Cal
Z}(T,\omega)\cap{\Cal M}^{ss}$, or equivalently of $\vert
\CT(\bob)_0\vert$ is a tedious exercise. Assume again that $r_1\leq
r_2\leq\ldots\leq r_n$. If $r_1$ is odd, then $\CT(\bob)$ is a ball
of radius $r_1/2$ around the midpoint of an edge. Hence in this
case
$$\vert\CT(\bob)_0\vert = 2(1+p^2+p^4+\ldots + p^{2(r_1-1)/2})\ \
.\tag8.19$$

Next assume that $r_1$ is even. Assume that $n=3$.

By Lemma 8.11, there are then seven cases:
\roster
\item"{\bf(1)}" $\gamma_2$ is split, $\gamma_3$
is unramified elliptic, i.e., $r_1$, $r_2$ even,
$\chi(-\varepsilon_1\varepsilon_2)=-1$,
$\chi(-\varepsilon_1\varepsilon_3)=1$,
\item"{\bf(2)}" $\gamma_2$ is split,
$\gamma_3$ is ramified elliptic, i.e., $r_1$ even, $r_2$ odd,
$\chi(-\varepsilon_1\varepsilon_2)=-1$,
\item"{\bf(3)}" $\gamma_2$ is unramified elliptic,
$\gamma_3$ is split, i.e., $r_1$, $r_2$ even,
$\chi(-\varepsilon_1\varepsilon_2)=1$,
$\chi(-\varepsilon_1\varepsilon_3)=-1$,
\item"{\bf(4)}" $\gamma_2$ is ramified elliptic,
$\gamma_3$ is split, i.e., $r_1$ odd, $r_2$ even,
$\chi(-\varepsilon_1\varepsilon_3)=-1$,
\item"{\bf(5)}" $\gamma_2$ and $\gamma_3$ are unramified elliptic,
i.e., $r_1$, $r_2$ even, $\chi(-\varepsilon_1\varepsilon_2)=1$,
$\chi(-\varepsilon_1\varepsilon_3)=1$,
\item"{\bf(6)}" $\gamma_2$ and $\gamma_3$ are ramified elliptic,
i.e., $r_1$, $r_2$ odd,
\item"{\bf(7)}" $\gamma_2$ and $\gamma_3$ are split,
i.e., $r_1$, $r_2$ even, $\chi(-\varepsilon_1\varepsilon_2)=-1$,
$\chi(-\varepsilon_1\varepsilon_3)=-1$.
\endroster
\medskip\noindent
{\it Let us determine the cardinality $\vert\CT(\bob)_0\vert$ in
the case (1).} At the fixed vertex $x_0=\CB^{\bob}$, there are
$p^2+1$ edges, $p^2-p$ of which lie outside of $\CB^0=\CB^{\b_1}$
and $p+1$ of which lie in $\CB^0$. If we move along a path
beginning with an edge running out of $\CB^0$, then our distances
$d_1$, $d_2$, and $d_3$ from $\CB^{\b_1}=\CB^0$, $\CB^{\b_2}$ and
$\CB^{\b_3}$ increase at the same rate. Hence, we remain inside of
$\CT(\bob)$ if and only if we move a distance at most $r_1/2$. The
number of vertices in $\CT(\bob)$ reached in this fashion is
therefore:
$$ 1 + (p^2-p)\big( 1+ p^2 +p^4+\dots+ p^{2(r_1/2-1)}\big).\tag8.20$$
The leading $1$ in this expression is the contribution of the
vertex $x_0$. Among the initial edges in $\CB^0$, there are two
which lie in the apartment $\CB^{0,\gamma_2}$, and $p-1$ which lie
outside of it. Suppose that we move a distance $j\ge1$, beginning
along one of these latter $p-1$ edges, and arrive at a vertex $y$.
The point $y$ has distances $d_1=0$, and $d_2=d_3=j$ from
$\CB^{\b_1}$, $\CB^{\b_2}$ and $\CB^{\b_3}$ respectively. If we
then move out of $\CB^0$ along one of the $p^2-p$ available initial
edges at $y$, we may move at most an additional
$\min(r_1/2,r_2/2-j)$ steps. There are
$$1 + (p^2-p)\big( 1+ p^2 +p^4+\dots+ p^{2(\min(r_1/2,r_2/2-j)-1)}\big)
\tag8.21$$
points of $\CT(\bob)$ reached in this fashion, and so the number of
points of $\CT(\bob)$ reached along the $p-1$ initial edges lying
in $\CB^0$ but outside of $\CB^{0,\gamma_2}$ is
$$(p-1)\sum_{j=1}^{r_2/2}
1 + (p^2-p)\big( 1+ p^2 +p^4+\dots+ p^{2(\min(r_1/2,r_2/2-j)-1)}\big).
\tag8.22$$
Finally, suppose that we begin along one of the two initial edges
in the apartment $\CB^{0,\gamma_2}$, and move a distance $k\ge 1$
in the apartment, arriving at a point $y$ with distances
$d_1=d_2=0$ and $d_3=k$. There are $p-1$ edges at $y$ which lie in
$\CB^0$ but outside $\CB^{0,\gamma_2}$, and $p^2-p$ edges lying
outside of $\CB^0$. The number of vertices of $\CT(\bob)$ reached
along paths beginning on the $p-1$ initial edges is $p-1$ times
$$\sum_{j=1}^{\min(r_2/2,r_3/2-k)}
1+(p^2-p)\big( 1+ p^2 +p^4+\dots+ p^{2(\min(r_1/2,r_2/2-j,r_3/2-k-j)-1)}
\big).\tag8.23$$
The number of vertices of $\CT(\bob)$ reached
along paths beginning on the $p^2-p$ `vertical' initial edges is
$$1+(p^2-p)\big( 1+ p^2 +p^4+\dots+ p^{2(\min(r_1/2,r_3/2-k)-1)}\big).
\tag8.24$$
Here the leading $1$ counts the vertex $y$ itself.
Combining these contributions, we obtain
$$
\align
&|\CT(\bob)_0|\tag8.25\\
\nass\nass
{}& =  1 + (p^2-p)\big( 1+ p^2 +p^4+\dots+ p^{2(r_1/2-1)}\big)\\
{}&\qquad+(p-1)\sum_{j=1}^{r_2/2}
1 + (p^2-p)\big( 1+ p^2 +p^4+\dots+ p^{2(\min(r_1/2,r_2/2-j)-1)}\big)\\
{}&\qquad+ 2\sum_{k=1}^{r_3/2}
1+(p^2-p)\big( 1+ p^2 +p^4+\dots+ p^{2(\min(r_1/2,r_3/2-k)-1)}\big)\\
\nass
{}&\qquad+ 2(p-1)\sum_{k=1}^{r_3/2}\\
{}&\qquad\sum_{j=1}^{\min(r_2/2,r_3/2-k)}
1+(p^2-p)\big( 1+ p^2 +p^4+\dots+ p^{2(\min(r_1/2,r_2/2-j,r_3/2-k-j)-1)}\big).
\endalign
$$

Similar explicit expressions may be obtained in the remaining 6 cases; we will not give them here.

We conclude this section by establishing a link between the
quantity $\vert\CT(\bob)_0\vert$ and the representation density of
certain quadratic forms.

The group $G'_p$ acts transitively on the set of $\zps$-lattices in
$U$. Fix a $\zps$-lattice $\LL_0$ in $U$, and let $K'_p$ be the
stabilizer of $\LL_0$ in $G'_p$. The stabilizer of the point
$[\LL_0]\in \CB$ is $K'_pZ'_p$, where $Z'_p$ is the group of scalar
matrices in $G'_p$. Let
$$N_0 =\{ \b\in V''_p;\  \b(\LL_0)\subset \LL_0\}.\tag8.26$$
Note that,
for $g\in G'_p$ and $\b\in V''_p$,
$$\b(g(\LL_0))\subset g(\LL_0)\iff g^{-1}\beta g\in N_0.\tag8.27$$
This condition depends only on the coset $gK'_pZ'_p$. Let
$\ph^{\prime,0}_p\in S(V_p'')$, the Schwartz space of $V''_p$, be
the characteristic function of the lattice $N_0$, and let
$$\phn=\ph^{\prime,0}_p\tt\dots\tt\ph^{\prime,0}_p\in S((V''_p)^n)
\tag8.28$$
be the characteristic function of $N_0^n$.
Let
$\bob\in (V''_p)^n$ be an $n$-tuple of special endomorphisms of $U$
such that $Q''(\bob)=p^{-1}T$. Then,
$$[\LL]=[g(\LL_0)]\in \CT(\bob)_0 \iff \phn(g^{-1}\bob g)=1.\tag8.29$$
It follows that
$$\align
|\CT(\bob)_0| &= |\{ gK'_pZ'_p;\  \phn(g^{-1}\bob g)=1\}|\tag
8.30\\
\nass
{}&= \int_{Z'_p\back G'_p} \phn(g^{-1}\bob g)\, dg,
\endalign$$
where the measure on $Z'_p\back G'_p$ is taken so that
$\vol(Z'_p\back K'_pZ'_p)=1$. Note that, if we write $\b\in V''_p$
in the form
$$\b = \b_0 \s = \pmatrix a&b\\c&a^\s\endpmatrix \s,\tag8.31$$
as in (8.12) above, then $\b_0\in G'_p$, and the action of $g\in
G'_p\subset GL_2(\qps)$ is given by twisted conjugacy of $\b_0$:
$$g^{-1}\b g = g^{-1}\b_0 g^\s \s.\tag8.32$$
Thus $|\CT(\bob)_0|$ is given by a kind of twisted orbital
integral. Of course, the integral in (8.30) will not be finite, in
general, if $G'_p(\bob)$, the pointwise stabilizer of $\bob$ is not
compact modulo the center.
\proclaim{Lemma 8.14} If $n=3$ or $4$, then
$G'_p(\bob)=Z'_p$, and the quantities in (8.30) are finite.
\endproclaim
\demo{Proof} Note that $G'_p/Z'_p\simeq SO(V''_p)$, and that, since
$\det(T)\ne 0$, the components of $\bob$ span an $n$-dimensional
non-degenerate subspace of $V''_p$. The group $G'_p(\bob)/Z'_p
\simeq SO(V''_p)(\bob)$ acts trivially on this span, and hence is
trivial (via, $\det=1$ when $n=3$). Thus, by Witt's Theorem, the
map $g\mapsto g^{-1}\bob g$ gives a bijection
-- in fact a homeomorphism -- of
$G'_p/Z_p$ to the hyperboloid
$$\Omega_T=\{ x\in (V''_p)^n;\ Q''(x)=p^{-1}T\}.\tag8.33$$
Since $\Omega_T$is closed in $(V''_p)^n$ and since $N_0^n$ is
compact, the support of the function $g\mapsto \phn(g^{-1}\bob g)$
is compact in $G'_p/Z'_p$, and the integral in (8.30) is finite.
\qed\enddemo

Finally, let $S$
be the matrix for the quadratic form $Q''$ on the rank $4$ $\Z_p$-lattice $N_0$, with
respect to some $\Z_p$-basis.
Then the classical representation density of $p^{-1}T$ by $S$ is
$$\align
\a(S,p^{-1}T) &= \lim_{t\rightarrow\infty} p^{t(\frac{n(n+1)}{2}-4n)}
\int\limits_{\scr x\in U^n\atop\scr S[x]-T\ \equiv\  0 \,(\text{\rm mod } p^t)}
\phn(x)\, dx\tag8.34\\
\nass
{}&= \int_{\Omega_T} \phn(x)\, d\mu_T(x)\\
\nass
{}&= \int_{G'_p/Z'_p} \phn(g^{-1}\bob g)\, d_Tg.
\endalign
$$
Here $dx$ is the Haar measure on $V''_p$ which is self dual with
respect to the bilinear form associated to $Q''$, and  $d\mu_T$ is
the measure on $\Omega_T$ which is defined by the limit in the
first line. Note that $d\mu_T$ is invariant under
$G'_p/Z'_p=SO(V''_p)$. Finally, $d_Tg$ is the Haar measure on
$G'_p/Z'_p$ coming from $d\mu_T$ under the homeomorphism
$G'_p/Z_p\simeq \Omega_T$. It only remains, then, to determine
$\vol(K'_pZ'_p/Z'_p,d_Tg)$.

Consider the case $n=3$. Let $\o$ be a top degree, translation
invariant, differential form on $(V''_p)^3$ whose associated
measure $dx$ is self dual with respect to the bilinear pairing
coming from $Q''$. Let $\eta$ be a top degree differential form on
${\roman{Sym}}_3$ such that the volume of ${\roman{Sym}}_3(\Z_p)$
for the associated measure is $1$. Let $\mu:(V''_p)^3\rightarrow
{\roman{Sym}}_3$ be the map taking $x$ to $Q''(x)$ or, rather, its
algebro-geometric counterpart. Note that $\mu$ is equivariant for
the action of $GL(3)$, and that, for $h\in GL(3)$,
$h^*(\o)=\det(h)^{4}\o$, and $h^*(\eta)=
\det(h)^{4}\eta$. There is a differential form $\tau$ on $(V''_p)^3$ of
degree $6$ such that
$$\o=\tau\wedge \mu^*(\eta),\tag8.35$$
and $h^*(\tau)= \tau$. Moreover, $\tau$ can be taken to be
invariant under the action of $SO(V''_p)$. Then, for any $T\in
\mu(V''_p)^3\subset {\roman{Sym}}_3(\Q_p)$ with $\det(T)\ne0$, the measure
$d\mu_T$ on $\Omega_T$, defined above, is precisely the measure
defined by the restriction of $\tau$ to $\Omega_T$. The invariance
properties of $\tau$ imply that the pullback $i^*_x(\tau)$ to
$SO(V)$ via a map $i_x:SO(V)\rightarrow (V''_p)^3$, $g\mapsto
g\cdot x$, where $\det(\mu(x))\ne0$, is independent of $x$. Thus,
the measure $d_Tg$ and the volume
$\kappa_p^{-1}:=\vol(K'_pZ'_p/Z_p,d_Tg)$ are independent of $T$.

\proclaim{Theorem 8.15} For any $T\in p {\roman{Sym}}_3(\Z_p)$ with $\det(T)\ne
0$, and for any $\bob\in (V''_p)^3$ with $T=p Q''(\bob)$,
$$|\CT(\bob)_0| = (1-p^{-4})^{-1}\cdot \a_p(S,p^{-1}T),$$
where $\a_p(S,p^{-1}T)$ is the representation density of the quadratic
form $p^{-1}T$ by the quadratic form $Q''$ on $N_0$.
\endproclaim

\demo{Proof}
By what was said above, both sides of this identity differ by a
multiplicative constant independent of $T$. To evaluate this
constant we take $T={\roman{diag}}(p,p,p)$. By Proposition 8.13 we
have $\vert\CT (\bob)_0\vert =1$. On the other hand, putting
$T'=p^{-1}T={\roman{diag}}(1,1,1)$ we have
$\alpha_p(S,T')=1-p^{-4}$, as is checked easily using the reduction
formulae.
\enddemo

{\bf Remark 8.16.} For a suitable basis for the lattice $N_0$ in
the $4$-dimensional $\Q_p$-vector space $V''_p$, we have
$$S=  \diag(1,1,1,\Delta),\tag8.36$$
where $\Delta\in \Z_p^\times$ is a {\it nonsquare}. In the case in
which $S$ is a hyperbolic space of dimension $4$, i.e., if $\Delta$
were a square, explicit formulas for the representation densities
$\a_p(S,T')$ where $T'=\diag(\varepsilon_1 p^{r_1},\varepsilon_2
p^{r_2},\varepsilon_3 p^{r_3})$ with $0\le r_1\le r_2\le r_3$, were
found by Kitaoka, comp.\ section 7. In our case, Theorem~8.15
reduces the computation of the densities, to the problem of
counting the number of vertices $|\CT(\bob)_0|$. As explained above
(at least in one of seven cases) this can be done explicitly. The
use of these explicit computations in the building is, however,
limited since one wants to determine more generally the
representation densities $\alpha_p(S_r, T')$ where $S_r=S\perp
H_{2r}$ for any $r\geq 0$, comp.\ section 7 above. In these more
general cases this combinatorial method seems difficult to handle.
However, just as in the case of the usual (twisted) orbital
integrals \cite{\langlands} it should not be forgotten.

We end this section with a global result. Let us fix $T\in
{\roman{Sym}}_3(\Z_{(p)})_{>0}$ and $\o\in V({\Bbb A}_f^p)^3$. For
$i=0,1$ let us consider the following subset $I(T,\o)$ of
$$\Omega'_T(\Q)\times Y_i\times G'({\Bbb A}_f^p)/K^p\ \
.\tag8.37$$
Here $G'$ is the usual inner form of $G$ (comp.\ beginning of
section 7) and $\Omega'_T$ is the hyperbola defined in (7.8). The
subset consists of the triples $(\underline x, L, gK^p)$ such that
\roster
\item"{(i)}" $g^{-1}\underline x g\in\o$
\item"{(ii)}" $L\in Y_i(\underline j)$, where $\underline
j=\underline x(p)$ is the $p$-component of $\underline x$.
\endroster
The set of irreducible components of the supersingular locus
$\CZ(T,\o)\cap{\Cal M}^{ss}$ breaks up into two disjoint subsets,
one for $i=0$ and one for $i=1$, which are interchanged by the
Frobenius in ${\roman{Gal}}({\Bbb F}_{p^2}/{\Bbb F}_p)$ (comp.\ end
of section 4). For each subset we have
$${\roman{Irred}}(\CZ(T,\o)^{ss})_i= G'(\Q)\setminus I(T,\o)\ \
.\tag8.38$$
The group $G'(\Q)/Z'(\Q)$ acts simply transitively on
$\Omega'_T(\Q)$. Let us fix a base point $\x$ and let us put
$Y_i(\j)=Y_i(\x(p))$. Let us also set
$$I(\x,\o)=\{ g\in G'({\Bbb A}_f^p)/K^p;\ g\x g^{-1}\in\o\}\ \
.\tag8.39$$
Then we may write $${\roman{Irred}}(\CZ(T,\o)^{ss})_i=
Z'(\Q)\setminus [Y_i(\j)\times I(\x,\o)]\ \ .\tag8.40$$ Let
$$Z'(\Q)^0=\{ z\in Z'(\Q);\ {\roman{ord}}_p(z)=0\}\ \
.\tag8.41$$
Then $Z'(\Q)^0$ acts trivially on $Y_i(\j)$. We have
$$Z'(\Q)= <p>\times Z'(\Q)^0\ \ ,\tag8.42$$
and $p$ acts faithfully on $Y_i(\j)$. We may identify the quotient
space for this action with $\CT(\j)_0$, comp.\ Corollary 8.4. Hence
we obtain
$${\roman{Irred}}(\CZ(T,\o)^{ss})_i= \CT(\j)_0\times
(Z'(\Q)^0\setminus I(\x,\o))\ \ .\tag8.43$$ The cardinality of the
second factor in (8.43) is given by an orbital integral, comp.\
Prop.\ 7.1,
$$\vert Z'(\Q)^0\setminus I(\x,\o)\vert =
{\roman{vol}}(K^p)^{-1}\cdot O_T (\varphi_f^p)\ \ .\tag8.44$$
Similarly, the cardinality of the first factor is given by a kind
of twisted orbital integral. For this fix a base lattice
$\Lambda_0\subset U$ and introduce the lattice of $V''_p$,
$$V''(\Z_p)= \{ \beta\in V''_p;\
\beta(\Lambda_0)\subset\Lambda_0\}\ \ .\tag8.45$$
Let
$$\varphi''_p={\roman{char}}\, V''(\Z_p)^3\ \ .\tag8.46$$
Then, setting $\bob =F^{-1}\j$, we have
$$\vert\CT(\j)_0\vert =
{\roman{vol}}(K_p)^{-1}\cdot\int\limits_{Z'(\Q_p)\setminus
G'(\Q_p)} \varphi''_p(g^{-1}\bob g)dg\ \ .\tag8.47$$ Taking (8.44)
and (8.47) together we obtain the following formula.
\proclaim{Proposition 8.17}
The number of irreducible components of the supersingular locus
$\CZ(T,\o)\cap{\Cal M}^{ss}$ is equal to
$$2\cdot{\roman{vol}}(K)^{-1}\cdot TO_T(\varphi''_p)\cdot
O_T(\varphi_f^p)\ \ .$$ Here $TO_T(\varphi''_p)$ denotes the
integral appearing in (8.47).
\endproclaim
Something analogous can, of course, also be done to write the
cardinality of the set of superspecial points of $\CZ(T,\o)$ as a product of
an orbital integral and a kind of twisted orbital integral. We do
not know whether these expressions have a global significance,
i.e.\ are related to automorphic forms (Eisenstein series?) or some
(relative?) trace formula.

\subheading{\Sec9. Components in the Siegel case}

The purpose of this section is to show that the method of the
previous section may also be applied to the case considered in the
companion paper \cite{\krsiegel}, with similar results. The notations
in this section will be different from those used in the rest of
the paper. As mentioned in the introduction, this section
benefitted from the corrections and suggestions of Ch.~Kaiser.
Also, it is based on his results in \cite{\kaiser}.

We begin by briefly recalling the situation considered in sections
4 and 5 of \cite{\krsiegel}. Recall from loc.\ cit.\ the polarized
isocrystal $({\Cal L},<\ ,\ >)$. Then $\dim_{\Cal K}{\Cal L}= 4$
and all slopes of $F$ are equal to $1/2$. Let
$$U=({\Cal L})^{p^{-1}F^2}\ \ .\tag9.1$$
Then $U$ is a $\Q_{p^2}$-form of ${\Cal L}$ and is equipped with
the symplectic form $<\ ,\ >$. Furthermore, $U$ is a $B'_p$-vector
space of dimension 2. Here
$$B'_p=\Q_{p^2}+\Q_{p^2}\Pi,\ \ \Pi^2=p,\ \ \Pi
a=a^{\sigma}\Pi\tag9.2$$ is the quaternion division algebra over
$\Q_p$ which acts on $U$ via $\Pi\mapsto F$. On $U$ we have the
$B'_p$-hermitian form (with respect to the main involution of
$B'_p$),
$$(\ ,\ ):U\times U\longrightarrow B'_p\tag9.3$$
defined by
$$(x,y)=-< x,\Pi y>\delta +<
x,y>\cdot\delta\Pi\ \ ,$$ for our fixed element
$\delta\in\Z_{p^2}^{\times}$ for which $\delta^{\sigma}=-\delta$.
The twisted form of the symplectic group may be identified with
$$G'(\Q_p)= \{ g\in GL_{B'_p}(U);\ (gx,gy)=\nu(g)\cdot (x,y),\
\nu(g)\in \Q_p^{\times}\}\ \ .$$
Recall the set of {\it distinguished lattices} in ${\Cal L}$,
$$\tilde X=\{ \tilde L\subset{\Cal L};\ \tilde L\supset F\tilde
L\supset p\tilde L,\ F\tilde L=c\cdot \tilde L^{\perp}\ \hbox{for
some}\ c\in {\Cal K}^{\times}\}\ \ .\tag9.4$$ A distinguished
lattice $\tilde L$ is stable under $p^{-1}F^2$ and defines
therefore a lattice $\tilde\Lambda =(\tilde L)^{p^{-1}F^2}$ in $U$.
This lattice is stable under the maximal order ${\Cal
O}'_p=\Z_{p^2} +\Z_{p^2}\Pi$ in $B'_p$. In this way we obtain a
one-to-one correspondence between $\tilde X$ and the set of ${\Cal
O}'_p$-lattices $\tilde\Lambda$ in $U$ such that
$$\tilde\Lambda^{\perp}=c\cdot\Pi\cdot\tilde\Lambda\ \ ,\ \
\hbox{some}\ c\in \Q_{p^2}^{\times}\ \ .\tag9.5$$
Here
$$\tilde\Lambda^{\perp}=\{ x\in U;\ (x,\tilde\Lambda)\subset {\Cal
O}'_p\}=\{ x\in U;\ <x,\tilde\Lambda >\subset\Z_p\}\ \ .\tag9.6$$
Recall the set of special endomorphisms of ${\Cal L}$,
$$V'_p=\{ j\in\End({\Cal L}, F);\ j^*=j,\ {\roman{tr}}(j)=0\}\ \
.\tag9.7$$
It is equipped with the quadratic form $q(j)=j^2$. For $\tilde L\in
\tilde X$ we consider
$$N_{\tilde L}
=\{ j\in V'_p;\ j(\tilde L)\subset\tilde L\}
=\{ j\in V'_p;\ j(\tilde\Lambda)\subset \tilde\Lambda\}\ \ .\tag9.8
$$
Here again $\tilde\Lambda =(\tilde L)^{p^{-1}F^2}$. As in the
previous section we fix $T\in{\roman{Sym}}_n(\Z_{(p)})$ with
$\det(T)\ne 0$ and $\j\in (V_p)^n$ with $q(\j)=T$. Consider the set
$$\tilde X(\j)=\{ \tilde L\in \tilde X;\ \j\in N^n_{\tilde L}\}\ \
.\tag9.9$$
Then, by theorems 5.12 and 5.13 of \cite{\krsiegel}, there is a
one-to-one correspondence between $\tilde X(\j)$ and those
projective lines ${\Bbb P}_{\tilde L}$ of special lattices whose
projection to the supersingular locus ${\Cal M}^{ss}$ lies entirely
in the image of the special cycle ${\Cal Z}(T,\omega)$, cf.\
\cite{\krsiegel}, section 5. Therefore it suffices to determine
$\tilde X(\j)$ resp.\ for a single $j\in V'_p$ the corresponding
set $X(j)$. Consider the building ${\Cal B}$ of $G'_{ad}(\Q_p)$. It
is a tree. Its vertices correspond to $B'_p$-homothety classes
$[\Lambda]$ of ${\Cal O}'_p$-lattices in $U$ which are fixed under
the involution induced by $\Lambda\mapsto\Lambda^{\perp}$, comp.\
\cite{\kaiser}, section 5. There are two types of vertices. One of
them comes from a distinguished lattice, i.e.\ is represented by a
lattice $\tilde\Lambda$ satisfying (9.5). We denote by ${\Cal B}_0$
the set of vertices of this type. The other is represented by a
lattice $\Lambda$ with
$$\Lambda^{\perp}=c\cdot\Lambda\ \ ,\ \ \hbox{some}\ c\in
\Q_{p^2}^{\times}\ \ ,\tag9.10$$
i.e.\ $\Lambda=(L)^{p^{-1}F^2}$ comes from a {\it superspecial
lattice} $L\subset{\Cal L}$ in the sense of section 4 of
\cite{\krsiegel}. An edge in ${\Cal B}$ is represented by two
lattices of opposite types with
$$\Pi\Lambda\subset\Lambda'\subset\Lambda\ \ .\tag9.11$$
Any $j\in V'_p$, with $q(j)\ne 0$ induces an involution of ${\Cal
B}$. Since we have a theory of elementary divisors for ${\Cal
O}'_p$-lattices and since ${\Cal B}$ is a tree the proof of
\cite{\krdrin}, Lemma 2.4 (i.e.\ of Kottwitz and Tate) applies
and gives the following result.

\proclaim{Lemma 9.1}
Let $j\in V'_p$ with $q(j)\ne 0$. For a ${\Cal O}'_p$-lattice
$\Lambda$ in $U$ representing a vertex $[\Lambda]\in {\Cal B}$ we
have
$$\align
j(\Lambda)\subset\Lambda &
\Longleftrightarrow d([j(\Lambda)], [\Lambda])\leq
2\ {\roman{ord}}_pq(j)\\ &
\Longleftrightarrow d([\Lambda], {\Cal B}^j)\leq
{\roman{ord}}_pq(j)\ \ .
\endalign$$
\endproclaim
We used here that $q(j)^2=Nm^0(j)$ when $j$ is considered as an
element of $GL_{B'_p}(U)\cong GL_2(B_p^{'{\roman{op}}})$.

For $j\in V'_p$ we introduce
$$\CT(j)= \{ x\in{\Cal B};\ d(x,{\Cal B}^j)\leq
{\roman{ord}}_p q(j)\}\tag9.12$$ and $$\CT(j)_0=\CT(j)\cap
\CB_0\tag9.13$$
Denoting by
$${\roman{pr}}:\tilde X\longrightarrow \CB_0\ \ ,\ \ \tilde
L\longmapsto [\tilde\Lambda]$$ the natural projection we therefore
obtain the following result.

\proclaim{Corollary 9.2}
For $j\in V'_p$ with $q(j)\ne 0$ we have
$$\tilde X(j)= {\roman{pr}}^{-1}(\CT(j)_0)\simeq \Z\times \CT(j)_0\
\ .$$
\endproclaim
On the right hand side the first factor $\Z$ enumerates the
connected components of $\tilde X(j)$.
\hfill\break
Similarly, if $\j\in(V'_p)^n$ with non-singular $T=q(\j)$ we may
change the basis of the $\Z_p$-span of the components such that the
components $j_i$ satisfy $q(j_i)\neq 0$. Then we put
$$\CT(\j)= \CT(j_1)\cap\ldots\cap \CT(j_n)\ \ \hbox{and}\ \
\CT(\j)_0= \CT(\j)\cap \CB_0\tag9.14$$
By Lemma 9.1 this is independent of the choice of basis.

\proclaim{Corollary 9.3}
$$\tilde X(\j)={\roman{pr}}^{-1}(\CT(\j)_0)=\Z\times \CT(\j)_0\ \
.$$
\endproclaim

It therefore remains to obtain a better understanding of the sets
$\CT(\j)$, especially of the fixed point sets in $\CB$ of special
endomorphisms.

We review the results of Kaiser \cite{\kaiser}. To state these
results we will call $j\in V'_p$ with $q(j)\ne 0$
\par\noindent
{\it split}, if $q(j)\in\Q_p^{\times,2}$
\hfill\break
{\it unramified elliptic}, if
$q(j)\in\Q_p^{\times}-\Q_p^{\times,2}$ and ${\roman{ord}}_pq(j)$
even
\hfill\break
{\it ramified elliptic}, if $q(j)\in \Q_p^{\times}-\Q_p^{\times,2}$
and ${\roman{ord}}_pq(j)$ odd.

We will also need the following notation. Let $E$ be a finite
extension of $\Q_p$. Then we will denote by $\CB(PGL_2(E))^+$ the
simplicial complex which is the first barycentric subdivision of
the building of $PGL_2(E)$. Thus there are two kinds of vertices in
$\CB(PGL_2(E))^+$: those that come from $\CB(PGL_2(E))$ which lie
on $q+1$ edges, where $q$ is the cardinality of the residue field,
and those that come from midpoints which lie on 2 edges. The first
type will be called {\it primary} and the second type {\it
secondary.}

\proclaim{Lemma 9.4} {\rm (Kaiser):}
Case by case we have for the fixed point set $\CB^j$ of a special
endomorphism $j$:
\hfill\break
{\rm if $j$ is split,} then $\CB^j$ consists of a single
superspecial vertex, cf.\ (9.10).
\hfill\break
{\rm if $j$ is unramified elliptic,} the $\CB^j$ is isomorphic to
$\CB(PGL_2(\Q_{p^2}))^+$ by an isomorphism which takes primary
vertices to distinguished vertices (cf.\ (9.5)) and secondary
vertices to superspecial vertices (cf.\ (9.10)).
\hfill\break
{\rm if $j$ is ramified elliptic,} then $\CB^j$ is isomorphic to
$\CB(PGL_2(\Q_p(j))$, where $\Q_p(j)$ is the ramified extension
generated by $j$.
\endproclaim

Here the isomorphisms are equivariant for the action of
$G'(\Q_p)_j\cong GL_2(\Q_p(j))$. The Lemma follows from
\cite{\kaiser}, Lemmas 5.2.2 and 5.2.3.\ together with Lemma 5.3.2.,
at least in the elliptic cases. The split case is analogous. Let us
give our version of Kaiser's isomorphisms.

{\it Let $j$ be split.} Correcting $j$ by an element of
$\Q_p^{\times}$ (which does not change the fixed point set) we may
assume $j^2=1$. The corresponding eigenspace decomposition of $U$
has the form
$$U=U_1\oplus U_{-1}\ \ .\tag9.15$$
Here $U_1$ and $U_{-1}$ are both $B'_p$-vector spaces of dimension
1 which are orthogonal to each other for the symplectic form $<\ ,\
>$, hence each of them is equipped with a symplectic form. A ${\Cal
O}'_p$-lattice $\Lambda$ in $U$ is fixed by $j$ iff
$$\Lambda=\Lambda_1\oplus \Lambda_{-1}\ \ ,\ \ \text{where}\
\Lambda_1=\Lambda\cap U_1\ \ ,\ \ \Lambda_{-1}=\Lambda\cap U_{-1}\
\ .\tag9.16$$
Since ${\roman{dim}}_{B'_p}U_{1}= {\roman{dim}}_{B'_p}U_{-1}=1$, we
may enumerate the lattices on the right hand side of (9.16) by
$\Z^2$ in such a way that passage to the dual lattice corresponds
to multiplication by $-1$,
$$\Lambda\longmapsto \Lambda^{\perp} \Longleftrightarrow
(n,m)\longmapsto (-n,-m)\ \ .\tag9.17$$ The identity
$\Lambda^{\perp}=\Pi^r\Lambda$ for $[\Lambda]\in\CB^j$ therefore
implies that
$$(-n,-m)= (n+r,m+r)\ \ ,\ \text{i.e.}\ \ r=-2n=-2m\ \ .$$
Hence the homothety class of $\Lambda$ is uniquely determined and,
since $r$ is even, $[\Lambda]$ is superspecial.
\par\noindent
{\it Next let $j$ be unramified elliptic.} Correcting $j$ by an
element of $\Q_p^{\times}$ we may assume $j^2=\Delta$. We obtain
the corresponding eigenspace decomposition,
$$U=U_1\oplus U_{-1}\ \ ,\tag9.18$$
where
$$U_1=\{ x\in U;\ j(x)=\delta x\}\ ,\ U_{-1}=\{ x\in U;\
j(x)=-\delta x\}\ \ .$$ Each of these subspaces is a 2-dimensional
vector space over the unramified quadratic extension $E=\Q_p(j)$,
they are orthogonal to each other, and ${\roman{deg}}\,\Pi=1$ for
this $\Z/2$-grading of $U$. The isomorphism
$$\CB(PGL_E(U_1))^+\buildrel\sim\over\longrightarrow
\CB^j\tag9.19$$
sends now
$$[pN\subseteq M\subseteq N]\ \ \text{to}\ \ [N\oplus \Pi^{-1}M]\ \
.$$
Here the brackets on the left indicate the class for the
equivalence relation for which
$$\align
pN\subseteq M\subseteq N\sim pN'\subseteq M'\subseteq
N'\Leftrightarrow \exists\, \alpha\in E^{\times}:(M',N') &
=(\alpha M,\alpha N)\\
\text{or}\ \ \ (M',N')
&
=(\alpha pN, \alpha M)\ \ .
\endalign$$

{\it Finally, let $j$ be ramified elliptic.} After correcting $j$
by an element in $\Q_p^{\times}$ we may assume that
$$j^2=\left\{\matrix
p\\
\Delta p
\endmatrix
\right.\tag9.20$$
Put
$$\tilde j=\left\{\matrix
\Pi^{-1}j\\
\varepsilon\cdot\delta^{-1}\Pi^{-1}j
\endmatrix
\right.\ \ .\tag9.21$$
Here $\varepsilon\in\Q_{p^2}$ is a fixed element with
$Nm_{\Q_{p^2}/\Q_p}(\varepsilon)=-1$. Then $\tilde j$ induces a
$\sigma$-linear automorphism of the $\Q_{p^2}$-vector space $U$
with $\tilde j^2=1$. We put
$$U_1=U^{\tilde j}\ \ ,\ \ U=U_1\otimes_{\Q_p}\Q_{p^2}\ \
.\tag9.22$$
Then $U_1$ is a vector space of dimension 2 over $E=\Q_p(j)$ and
the isomorphism
$$\CB(PGL_E(U_1))\buildrel\sim\over\longrightarrow \CB^j\tag9.23$$
sends
$$[M]\ \ \text{to}\ \ [M\otimes_{\Z_p}\Z_{p^2}]\ \ .$$
We next investigate the relative position of the fixed point sets
of two special endomorphisms $j,j'$ which are orthogonal to each
other, i.e., $jj'=-j'j$. We distinguish the following mutually
exclusive cases.

{\bf First case:} {\it One of $j,j'$ is split.} Then $\CB^{j,j'}$
is a single superspecial vertex, cf.\ Lemma 9.4 above.

{\bf Second case:} {\it One of $j,j'$ is unramified elliptic.} Let
us assume that $j$ is unramified elliptic, and let us consider the
corresponding eigenspace decomposition (9.18). With respect to this
decomposition we have ${\roman{deg}}\, \Pi^{-1}j'=0$ and the
induced endomorphism
$$\Pi^{-1}j':\ U_1\longrightarrow U_1\tag9.24$$ is $\tau$-linear,
where $\tau$ is the generator of
${\roman{Gal}}(E/\Q_p)={\roman{Gal}}(\Q_p(j)/\Q_p)$. Furthermore we
have a bijection between fixed point sets
$$(\CB(PGL_E(U_1))^+)^{\Pi^{-1}j'}\buildrel\sim\over\longrightarrow
\CB^{j,j'}\ \ .\tag9.25$$
Indeed, if $pN\subseteq M\subseteq N$ is a representative of a
vertex of $\CB(PGL_E(U_1))^+$, then the corresponding vertex
$[N\oplus
\Pi^{-1}M]$ of $\CB$ is fixed under $j'$ iff
$$\align
[\Pi^{-1}j'(M)\oplus j'(N)]= & [N\oplus \Pi^{-1}M]\Leftrightarrow\\
& [p\cdot\Pi^{-1}j'(M)\subseteq p\cdot\Pi^{-1}j'(N)\subseteq
\Pi^{-1}j'(M)] =[pN\subseteq M\subseteq N]\ .
\endalign$$
We also note that
$$(\Pi^{-1}j')^2=p^{-1}q(j')\ \ .\tag9.26$$
Hence we see that if $j'$ is unramified elliptic, then $\Pi^{-1}j'$
fixes a unique secondary vertex of $\CB(PGL_E(U_1))^+$, hence
$\CB^{j,j'}$ consists of a single superspecial vertex of $\CB$. If
$j'$ is ramified elliptic, then after correcting $j'$ by an element
in $\Q_p^{\times}$ we may assume that $j^{'2}$ is equal to $p$ or
to $\Delta p$, comp.\ (9.20). As in (9.21) we put
$$\tilde j'=\left\{ \matrix
\Pi^{-1}j'
&\text{if}\ j^{'2}=p\hfill\\
\varepsilon\cdot\delta^{-1}\Pi^{-1}j'
&
\text{if}\ j^{'2}=\Delta p\ \ .\hfill
\endmatrix\right.
\tag9.27$$
Then $\tilde j'$ is a $\tau$-linear automorphism of $U_1$ with
$\tilde j^{'2}=1$. We may thus write
$$U_1= U_{1,1}\otimes_{\Q_p}E\ \ ,\ \ U_{1,1}=U_1^{\tilde j'}\ \
.\tag9.28$$
Since the extension $E$ of $\Q_p$ is unramified we have bijections
(unramified descent)
$$\CB(PGL_{\Q_p}(U_{1,1}))^+\buildrel\sim\over\longrightarrow (\CB
(PGL_E(U_1))^+)^{\tilde j}\buildrel\sim\over\longrightarrow
\CB^{j,j'}\ \ .\tag9.29$$
Hence in this case the common fixed point set of $j,j'$ is
isomorphic to $\CB(PGL_2(\Q_p))^+$. Under this isomorphism primary
vertices correspond to distinguished vertices and secondary
vertices to superspecial vertices.

{\bf Third case:} {\it Both $j,j'$ are ramified elliptic.} After
correcting $j$ and $j'$ by elements of $\Q_p^{\times}$, we may
assume that $q(j), q(j')\in\{ p,\Delta p\}$. Let us define $\tilde
j$ as in (9.21) and put $U_1=U^{\tilde j}$. Then $U_1$ is a
2-dimensional vector space over $E=\Q_p(j)$. Define now
$$\tilde j'=\left\{\matrix
\delta^{-1}\Pi^{-1}j'
&
\text{if}\ q(j)=p\hfill\\
\varepsilon\cdot\Pi^{-1}j'
&
\text{if}\ q(j)=\Delta p\ \ .\hfill
\endmatrix\right.\tag9.30$$
Then
$$\tilde j\cdot \tilde j'=\tilde j'\cdot \tilde j\tag9.31$$
Hence $\tilde j'$ induces an automorphism of $U_1$ which is
$\tau$-linear $(<\tau>={\roman{Gal}}(E/\Q_p))$ and such that
$$\tilde j^{'2}=\left\{\matrix
-\Delta^{-1}
&
\text{if}\ q(j)=q(j')=p\hfill\\
-1
&
\text{if}\ q(j)=p,\ q(j')=\Delta p\hfill\\
-1
&
\text{if}\ q(j)=\Delta p,\ q(j')=p\hfill\\
-\Delta
&
\text{if}\ q(j)=\Delta p,\ q(j')=\Delta p\hfill
\endmatrix\right.\tag9.32$$
For fixed $j$, precisely one possible $j'$ has the property that
$\tilde j^{'2}\equiv 1\, {\roman{mod}}\, \Q_p^{\times,2}$ and for
the other possibility for $j'$ we have $\tilde j^{'2}\not\in
Nm_{E/\Q_p}(E^{\times})$. For $\tilde j^{'2}\in\Q_p^{\times,2}$, we
can correct $\tilde j'$ by an element in $\Q_p^{\times}$, so as to
obtain a $\tau$-linear automorphism $\tilde j''$ of $U_1$ with
$(\tilde j'')^2=1$. Hence
$$U_1=U_{1,1}\otimes_{\Q_p}E\ \ ,\ \ U_{1,1}=U_1^{\tilde j''}\ \
.\tag9.33$$
In this case we obtain bijections
$$\CB(PGL_{\Q_p}(U_{1,1}))^+\buildrel\sim\over\longrightarrow
\CB(PGL_E(U_1))^{\tilde j''}\buildrel\sim\over\longrightarrow
\CB^{j,j'}\tag9.34$$
(note that we are dealing here with a {\it ramified} descent from
$E$ to $\Q_p$, which explains the barycentric subdivision). Under
this identification primary vertices correspond to superspecial
vertices and secondary vertices correspond to distinguished
vertices. Indeed, ${\Cal B}^{j,j'}$ is contained in the fixed point
locus of $E_1^{\times}$ where $E_1= \Q_p(jj')$ is the unramified
quadratic field extension of $\Q_p$ generated by $jj'$, hence the
assertion follows from \cite\kaiser, Lemma 5.3.2.\ a). Since
$\chi((jj')^2)=-\chi(\tilde j^{'2})= -1$, the element $jj'$ indeed
generates a field extension.
\hfill\break
Now assume that $\tilde j^{'2}\not\in \Q_p^{\times, 2}$, i.e.\ that
$(jj')^2=c^2\in \Q_p^{\times,2}$. Put $j_0=c^{-1}(jj')$. We
consider the eigenspace decomposition corresponding to $j_0$,
$$U=U_1\oplus U_{-1}\tag9.35$$
This time, in contrast to (9.15), since $j_0^{\ast}=-j_0$, both
$U_1$ and $U_{-1}$ are isotropic and $(\ ,\ )$ induces a perfect
pairing between these one-dimensional $B'_p$-vector spaces. A
${\Cal O}_p$-lattice $\Lambda$ in $U$ is fixed by $j_0$ iff
$$\Lambda=\Lambda_1\oplus \Lambda_{-1}\ \ ,\ \ \hbox{with}\
\Lambda_1=U_1\cap \Lambda\ \hbox{and}\ \Lambda_{-1}=U_{-1}\cap
\Lambda\ \ .\tag9.36$$
Let $M\subset U_1$ be a ${\Cal O}_{B'_p}$-lattice and let
$M^{\bot}\subset U_{-1}$ be the dual lattice. Then
$$\align
{\Cal B}^{j_0} &
=\{ \Pi^{\Z}\cdot M\} \times \{ \Pi^{\Z}\cdot M^{\bot}\} /
B_p^{\times}\tag9.37\\ &
=(\Z\oplus\Z) /\Z\ \ ,
\endalign$$
with the diagonal action of $\Z$. The action of $j$ on ${\Cal
B}^{j_0}$ is induced by
$$j:(m,n)\mapsto (n+i, m-i+2)\ \ ,\tag9.38$$
since ${\roman{ord}}_p(j^2)=1$. The integer $i$ is independent of
the choice of $M$. It follows that
$${\Cal B}^{j,j'}={\Cal B}^{j_0, j}=\{ (i-1, 0)\}\in (\Z\oplus \Z)
/\Z\ \ .\tag9.39$$
Furthermore, this vertex is superspecial if $i-1$ is even and
distinguished if $i-1$ is odd. Consider now the hermitian form on
$U_1$,
$$\{ x,y\}= (x,j(y))\ \ .\tag9.40$$
Now any hermitian form on a one-dimensional $B'_p$-vector space
admits a selfdual $O_{B'_p}$-lattice. If $M'=\Pi^kM$ is selfdual,
then
$$O_{B'_p}=\{ M', M'\} = (M', j(M'))= (\Pi^kM, \Pi^{k-i+2}\cdot
M^{\bot} )= \Pi^{2k-i+2}\cdot O_{B'_p}\ \ .$$ It follows that $i$
is even and hence ${\Cal B}^{j,j'}$ consists in this case of a
single distinguished vertex.

Let us consider now a 4-tuple of special endomorphisms $j_1,\ldots,
j_4$ with $q(j_i)\in\Q_p^{\times}$ and which are pairwise
orthogonal to each other. The following Lemma is the analogue of
Lemma 8.12.

\proclaim{Lemma 9.5}
$\CB^{\j}:= \CB^{j_1}\cap\ldots\cap\CB^{j_4}$ is a vertex. More
precisely, ${\Cal B}^{\bold j}$ is a superspecial vertex if at
least one $j_i$ is split or at least two $j_i$ are unramified
elliptic and ${\Cal B}^{\bold j}$ is a distinguished vertex if at
least three $j_i$ are ramified elliptic.
\endproclaim

\demo{Proof} The cases where at least one $j_i$ is split, or at
least two $j_i$ are unramified elliptic have already been taken
care of. In these cases the common fixed point is a superspecial
vertex.
\hfill\break
Let us now assume that $j_1,j_2, j_3$ are all ramified elliptic. We
proceed as in the third case above, starting with $j=j_1$ and
$j'=j_2$. Therefore we introduce $\tilde j=\tilde j_1$ as in (9.21)
and $U_1=U^{\tilde j}$. We also introduce $\tilde j'=\tilde j_2$ as
in (9.30). If $\tilde j_2^2\not\in\Q_p^{\times,2}$, then the common
fixed point set of $j_1,j_2$ is a distinguished vertex which is
then equal to $\CB^{\j}$. If $\tilde j_2^2\in\Q_p^{\times,2}$ we
rescale $\tilde j_2$ to obtain $\tilde j''_2$ with $(\tilde
j''_2)^2=1$, and write $U_1=U_{1,1}\otimes_{\Q_p}E$, for
$U_{1,1}=U_1^{\tilde j''_2}$ and $E=\Q_p(j_1)$. Finally we
introduce $\tilde j_3$ as in (9.30). Then $\tilde j_3$ commutes
with $\tilde j_1$ and anticommutes with $\tilde j''_2$. Now $j_1$
acts on $U_1$; we set $U_{1,-1}=j_1(U_{1,1})$. Then we obtain a
$\Z/2$-grading of $U_1$,
$$U_1=U_{1,1}\oplus U_{1,-1}\ \ .\tag9.41$$
Then $U_{1,-1}$ ist the $(-1)$-eigenspace of the $\Q_p$-linear
endomorphism $\tilde j''_2$ of $U_1$. Since $\tilde j_3$
anticommutes with $\tilde j''_2$, it has degree one with respect to
this grading. The primary vertices of ${\Cal
B}(PGL_{\Q_p}(U_{1,1}))= {\Cal B}^{j_1, j_2}$ are represented by
lattices $M\subset U_{1,1}$ and such a lattice is represented by
the lattice
$$M\oplus j_1(M)= M\otimes_{\Z_p}O_E\tag9.42$$
in $U_1$. Now
$$\tilde j_3(M\oplus j_1(M))=\tilde j_3j_1(M)\oplus \tilde j_3(M)=
j_1\tilde j_3(M)\oplus j_1(j_1 \tilde j_3(M))\ \ .$$ Hence since
$j_1 \tilde j_3$ acts on $U_{1,1}$ with ${\roman{ord}}_p((j_1
\tilde j_3)^2)=1$, we see that ${\Cal B}^{j_1,j_2,j_3}= ({\Cal
B}(PGL_{\Q_p}(U_{1,1}))^+)^{j_1,\tilde j_3}$ consists of a single
secondary vertex. By the third case above, this corresponds to a
distinguished vertex of ${\Cal B}$.
\qed\enddemo

{\bf Remark 9.6.} We thus see that in most cases it is sufficient
to intersect the fixed point loci of 3 special endomorphisms to get
down to a vertex. The only case when 4 are needed is if one of them
is unramified elliptic and the others are ramified elliptic. Let us
discuss the case when $j_1$ is unramified elliptic and
$j_2,j_3,j_4$ are ramified elliptic. We introduce the
$\delta$-eigenspace $U_1$ for $j_1$ which is a 2-dimensional vector
space over $E=\Q_p(j_1)$. The special endomorphism $j_2$ induces
the $\tau$-linear endomorphism of $U_1$
$$\tilde j_2:=\left\{\matrix
\Pi^{-1}j_2
&
\text{if}\ q(j_2)=p\hfill\\
\varepsilon\delta^{-1}\Pi^{-1}j_2
&
\text{if}\ q(j_2)=\Delta p\ \ .\hfill
\endmatrix\right.
\tag9.43$$
Of course, here we again scaled $j_2$ by an element in
$\Q_p^{\times}$. As in (9.28) we have $U_1=U_{1,1}\otimes_{\Q_p}E$
with $U_{1,1}=U_1^{\tilde j_2}$. As in (9.30) above we put for
$i=3,4$
$$\tilde j_i=\left\{\matrix
\delta^{-1}\Pi^{-1}j_i
&
\text{if}\ q(j_2)=p\hfill\\
\varepsilon\cdot\Pi^{-1}j_i
&
\text{if}\ q(j_2)=\Delta p\ \ .\hfill
\endmatrix\right.\tag9.44$$
Then $\tilde j_3$ and $\tilde j_4$ are endomorphisms of the
2-dimensional $\Q_p$-vector space $U_{1,1}$ which anticommute and
whose squares are scalar matrices of even valuation. Furthermore,
we have
$$\tilde j_i^2 =\left\{\matrix
-\Delta^{-1}p^{-1}q(j_i)\hfill\\
&
\equiv -\Delta\cdot q(j_2)\cdot q(j_i)\ {\roman{mod}}\
\Q_p^{\times,2}\ \ .\\
-p^{-1}q(j_i)\hfill
\endmatrix\right.\tag9.45$$
Hence $\tilde j_i$ is unramified elliptic if
$\chi(-q(j_2)q(j_i))=1$ and is split if $\chi(-q(j_2)q(j_i))=-1$.
Two possibilities occur: either one of $\tilde j_3$, $\tilde j_4$
is unramified elliptic in which case this involution of ${\Cal
B}(PGL_2(\Q_p))= {\Cal B}(PGL_{\Q_p}(U_{1,1}))$ has a primary
vertex as its fixed point set, or both $\tilde j_3$, $\tilde j_4$
are split and have apartments as their fixed point sets which
intersect in a primary vertex. In either case the common fixed
point set is a primary vertex which by the second case preceding
Lemma 9.5.\ is mapped to a distinguished vertex of ${\Cal B}$.

{\bf Remark 9.7:} Note that very subtle relations hold between $T$
and $\CT(\j)_0$. For instance we see that if $\varepsilon_1=1$ and
$a_1=0$, then $\CT(\j)_0$ is empty as it is contained in the ball
of radius 0 around a superspecial vertex. Of course, we already
knew this from \cite{\krsiegel} via Corollary 9.3 above, since if $T$
represents 1, then $\tilde X(\j)=\emptyset$ (only isolated points
in the special cycle $\CZ(T,\o))$. A more surprising conclusion is
that the following combination is excluded:
\hfill\break
$\chi(\varepsilon_1)=1$, $a_1$ even
\hfill\break
$\chi(\varepsilon_2)=-1$, $a_2$ even
\hfill\break
$a_3$ and $a_4$ odd and $\chi(-\varepsilon_3\varepsilon_4)=1$.
\hfill\break
Indeed, the fixed point set of $j_1$ would be a superspecial vertex
and the common fixed point set of $j_2,j_3,j_4$ would be a
distinguished vertex, a contradiction to the fact that $\CB^{\j}$
is non-empty. It seems, however, that all these cases can already
be excluded on the basis of the following remark. If $T\in
{\roman{Sym}}_4(\Q_p)$ with ${\roman{det}}(T)\ne 0$ is represented
by $V'_p$ then
$$-1=\varepsilon(V'_p)=\varepsilon(T)\cdot ({\roman{det}}(T),
-{\roman{det}}(V'_p))_p=\varepsilon(T)\cdot ({\roman{det}}(T), -1)_p\ \ ,$$
(\cite{\annals}, (1.16)) where $\varepsilon(V'_p)$ resp.\
$\varepsilon(T)$ denotes the Hasse invariant of $V'_p$ resp.\ $T$.
This remark is to be compared with Remark 2.10 in \cite{\krdrin}
(in this case, too, the exclusion principle by types of vertices is
implied by the exclusion principle by the representability of $T$
by the relevant quadratic form).

\proclaim{Proposition 9.8}
Let $T\in{\roman{Sym}}_4(\Z_{(p)})_{>0}$ be $GL_4(\Z_p)$-
equivalent to
$$T\sim{\roman{diag}}(\varepsilon_1p^{a_1}, \varepsilon_2p^{a_2},
\varepsilon_3p^{a_3}, \varepsilon_4p^{a_4})\ \ ,\ \ 0\le
a_1\le\ldots \le a_4\ \ .$$  Then each connected component of the
supersingular locus in the image of $\CZ(T,\o)$ is irreducible of
dimension one, or equivalently $\vert\CT(\j)_0\vert=1$ for any
$\j\in V_p^4$ with $q(\j)=T$, if and only if the following
conditions hold:
\hfill\break
a) Assume that $a_1$ is even. Then
\roster
\item"{(i)}" $a_1=0$  and $\chi(\varepsilon_1)=-1$
\item"{(ii)}" $a_2,a_3,a_4$ are all odd
\item"{(iii)}" if $\chi(-\varepsilon_2\varepsilon_3)=1$ then
$a_3=1$ and if $\chi(-\varepsilon_2\varepsilon_3)=-1$ then $a_4=1$.
\endroster

b) Assume that $a_1$ is odd. Then
\roster
\item"{(i)}" $a_1=1$
\item"{(ii)}" if $\chi(-\varepsilon_1\varepsilon_2)=1$ then
$a_2=1$ and if $\chi(-\varepsilon_1, \varepsilon_2)=-1$ then
$a_3=1$.
\endroster

\endproclaim
\demo{Proof}
a) Let $y$ be the common fixed vertex of $j_1,\ldots, j_4$. If
$a_1\ge 2$, then the ball with radius 2 around $y$ is contained in
$\CT(\j)$. Since this ball contains more than one distinguished
vertex we see that if $\vert\CT(\j)_0\vert=1$, we must have $a_1\le
1$. If $a_1=0$ and $\chi(\varepsilon_1)=1$, then
$\CT(\j)=\emptyset$.  Hence if $\vert\CT(\j)_0\vert=1$, then
condition (i) holds. Also (ii) holds. Indeed, otherwise by our
previous results $y$ would be a superspecial vertex which
corresponds to the midpoint of an edge in $\CB^{j_1}\simeq
\CB(PGL_2(\Q_{p^2}))^+$. The two primary vertices $x_1,x_2$
of this edge have distance 1 to $y$. Let $\CT(\j)_0=\{ x\}$. Since
$a_1=0$, we have $x\in\CB^{j_1}$ and $d(x_1, \CB^{j_i})=d(x_2,
\CB^{j_i})\le d(x, \CB^{j_i})$ for any $i\ge 2$ and hence $\{
x_1,x_2\}\subset\CT(\j)_0$, a contradiction.  Hence (ii) holds and
$y$ is a distinguished vertex, and $\CB^{j_1,j_2}\simeq
\CB(PGL_2(\Q_p))^+$. Furthermore, the fixed point set of $j_3$ in
$\CB^{j_1, j_2}$ is either $\{ y\}$ or an apartment in
$\CB(PGL_2(\Q_p))^+$ depending on whether
$\chi(-\varepsilon_2\varepsilon_3)=1$ or
$\chi(-\varepsilon_2\varepsilon_3)=-1$, cf.\ third case before
Lemma 9.5 (comp.\ Remark 9.6 above). In the first case $a_3\le 1$
because otherwise all distinguished vertices of distance 2 from $y$
in $\CB(PGL_2(\Q_p))^+$ would be contained in $\CT(\j)_0$. In the
second case, the fixed point set of $j_4$ in $\CB^{j_1,j_2}$ is an
apartment in $\CB(PGL_2(\Q_p))^+$ which intersects the previous
apartment in $y$ and hence $a_4\le 1$. We have therefore shown that
$\vert\CT(\j)_0\vert
=1$ implies all conditions (i) -- (iii).
\hfill\break
Conversely, if conditions (i) -- (iii) hold, then $\CB^{\j}$ is a
distinguished vertex $y$ and $y\in\CT(\j)_0$. Let $x\in\CT(\j)_0$.
Since $a_1=0$, we have $x\in \CB^{j_1}$. For any $i\ge 2$, looking
back at the second case before Lemma 9.5, the closest point to $x$
in $\CB^{j_i}$ cannot be a midpoint of an edge in $\CB^{j_1}\cong
\CB(PGL_2(\Q_{p^2}))^+$, hence either $x\in\CB^{j_1,j_i}$ or $d(x,
\CB^{j_i})\ge 2$. The first possibility cannot hold for all $i$ if
$x\ne y$. But for $i$ where the second possibility holds, this
implies $a_i\ge 2$. The condition (iii) implies that we are in the
first case of (iii) and $i=4$. But then $x\in\CB^{j_1,j_2,j_3}= \{
y\}$.

b) Let again $y$ be the common fixed vertex of $j_1,\ldots, j_4$.
As before the case $a_1\geq 3$ is excluded since otherwise the ball
with radius 3 would be contained in ${\Cal T}({\bold j})$ and hence
${\Cal T}({\bold j})_0$ would have more than one element. For the
same reason $y$ has to be distinguished and ${\Cal T}({\bold j})_0=
\{ y\}$. Consider ${\Cal B}^{j_1}\simeq {\Cal
B}(PGL_2(\Q_p(j_1))$. Considering the number of neighbours of a
vertex contained in ${\Cal B}^{j_1}$ in ${\Cal B}$ resp.\ ${\Cal
B}^{j_1}$ we see that all neighbours of a superspecial vertex in
${\Cal B}^{j_1}$ have to lie in ${\Cal B}^{j_1}$. Hence ${\Cal
T}(j_1)_0= {\Cal B}_0^{j_1}$ (set of distinguished vertices in
${\Cal B}^{j_1}$). Since ${\Cal T}({\bold j})_0=\{ y\}$, it follows
that $a_2=1$, since otherwise the distinguished vertices in the
ball with radius 2 around $y$ would be contained in ${\Cal
T}(\j)_0$. It follows now by our previous results that if
$\chi(-\varepsilon_1 \varepsilon_2)=1$ then ${\Cal B}^{j_1, j_2}=
\{ y\}$ and if $\chi(-\varepsilon_1\varepsilon_2)=1$ then ${\Cal
B}^{j_1,j_2}\cong {\Cal B}(PGL_2(\Q_p))^+$ where primary vertices
correspond to superspecial vertices in ${\Cal B}$ and secondary
vertices correspond to distinguished vertices in ${\Cal B}$. In the
first case ${\Cal T}(j_1)_0\cap {\Cal T}(j_2)_0= {\Cal B}_0^{j_1,
j_2}=\{ y\}$ and we are in the case of the first clause of (ii). In
the second case it follows that $a_3=1$ because otherwise the
secondary vertices in the ball of radius 2 around $y$ in ${\Cal
B}(PGL_2(\Q_p))^+$ would all be contained in ${\Cal T}(\j)_0$. In
this case ${\Cal T}(j_1)_0\cap {\Cal T}(j_2)_0\cap {\Cal T}(j_3)_0=
{\Cal B}_0^{j_1, j_2, j_3}=\{ y\}$. This proves the necessity of
conditions (i) and (ii), and the sufficiency is proved in a similar
way.
\qed\enddemo

{\bf Remark:} The cases enumerated in Proposition 9.8 should be
considered first when investigating cases of excess intersection of
special cycles.

We finally remark that, as in Theorem 8.15, there is a close
relation between $\vert\CT(\j)_0\vert$ and the representation
density of $T$ by $V'_p$. Also the number of irreducible components
of the image of $\CZ(T,\o)$ in ${\Cal M}^{ss}$ may be expressed as
a product of two orbital integrals, one prime to $p$, the other
$p$-adic. As in Proposition 8.16, the global significance of such
an expression is unclear.

\subheading{\Sec10. Special cycles on the special fibre}

So far most of our detailed results concerned the intersection of
the supersingular locus in the special fibre with (the image of) a
special cycle $\CZ(T,\o)$. In this section we will take a more
global view and also consider the intersection of $\CZ(T,\o)$ with
the ordinary locus.
\par
We will use the following notations. By ${\Cal M}^{\roman{ord}}$
resp.\ $\CZ(T,\o)^{\roman{ord}}$ we denote the open complement of
the supersingular locus. Also let
$$\overline{\Cal M} ={\Cal M}\times_{\roman{Spec}\,\Z_{(p)}}
{\roman{Spec}}\, {\Bbb F}_p\ \ ,\ \ \overline{\CZ}(T,\o)=
\CZ(T,\o)\times_{\roman{Spec}\,\Z_{(p)}}{\roman{Spec}}\, {\Bbb
F}_p\tag10.1$$ and $\overline{\Cal M}^{\roman{ord}}$,
$\overline{\CZ}(T,\o)^{\roman{ord}}$ denote the special fibres.
\par
Let $(A,\lambda,\iota,\overline{\eta}^p)\in {\Cal
M}^{\roman{ord}}({\Bbb F})$. The action of ${\Cal O}_C\otimes
\Z_{p^2}=M_4(\Z_{p^2})$ on the $p$-divisible group $A(p)$ gives a
decomposition $A_0(p)={\Cal A}^4$, where ${\Cal A}$ is a
$p$-divisible group of dimension 2 and height 4, comp.\ section 4.
Let $M$ be the Dieudonn\'e module of ${\Cal A}$. It comes equipped
with the action of $\Z_{p^2}$,
$$\iota_0:\Z_{p^2}\to{\roman{End}}(M)\ \ .\tag10.2$$
We consider the decomposition of $M$ into its \'etale and its
infinitesimal component,
$$M=M^0\oplus M^1\ \ .\tag10.3$$
Each summand is $\Z_{p^2}$-stable and in fact
$$M^0\cong (\Z_{p^2}\otimes W,{\roman{id}}\otimes\sigma_W)\ \ ,\ \ M^1
\cong (\Z_{p^2}\otimes
W,p\cdot {\roman{id}}\otimes\sigma_W)\ \ .\tag10.4$$ The
alternating form $<\ ,\
>$ induced by the polarization identifies $M^1$ with the Serre-dual
of $M^0$.

Now let $j$ be a special endomorphism of
$(A,\lambda,\iota,\overline{\eta}^p)$. It induces an endomorphism
of $M$ of the form
$$j=j^0\oplus j^1\ \ ,\tag10.5$$
where $j^1$ is the dual to $j^0$ and where
$$j^0:M^0\to M^0\tag10.6$$
is linear for the $W$-action but antilinear for the
$\Z_{p^2}$-action. Denoting as usual by $\sigma$ the non-trivial
automorphism of $\Z_{p^2}$, and using the identification (10.4) we
may write
$$j^0=a\sigma\ \ ,\ \ a\in\Z_{p^2}\ \ ,\tag10.7$$
and then
$$Q(j)=Nm(a)\ \ .\tag10.8$$
We may therefore state the following fact.

\proclaim{Lemma 10.1}
Let $T\in {\roman{Sym}}_n(\Z_{(p)})_{>0}$ with
$\overline{\CZ}(T,\o)^{\roman{ord}}\ne\emptyset$. Then $n\leq 2$
and $T$ is represented by the norm form on $\Z_{p^2}$. In
particular, the $p$-adic ordinal of each diagonal entry of $T$ is
even.\qed
\endproclaim

We next use the theory of Serre-Tate canonical coordinates to
investigate the infinitesimal deformation of
$\CZ(T,\o)^{\roman{ord}}$ at $(A,\lambda,\iota,\overline{\eta}^p)$.
The $p$-divisible group ${\Cal A}$ is an extension
$$0\to L\otimes \hat{\Bbb G}_m\to {\Cal A}\to L'\otimes
\Q_p/\Z_p\to 0\tag10.9$$
where $L$ and $L'$ are free $\Z_{p^2}$-modules of rank one. The
polarization allows us to identify $L'$ with the $\Z_p$-dual of
$L$. By Serre-Tate theory we may identify the formal deformation
space of ${\Cal A}$ with
$${\roman{Hom}}(L',L)\otimes \hat{\Bbb G}_m\ \ .\tag10.10$$
Inside this formal torus of relative dimension 4 the locus where
the action $\iota_0$ and the polarization lifts is given by the
formal subtorus of relative dimension 2 (we identify $L'$ with the
dual of $L$),
$${\Cal T}={\roman{Sym}}^2_{\Z_{p^2}}(L)\otimes \hat{\Bbb G}_m\ \
.\tag10.11$$
Fix a generator $e$ of the $\Z_{p^2}$-module $L$. Then $e\otimes e$
is a generator of the $\Z_{p^2}$-module
${\roman{Sym}}^2_{\Z_{p^2}}(L)$ and to any multiple $B(e\otimes
e)\in {\roman{Sym}}^2_{\Z_{p^2}}(L)\otimes {\Bbb G}_m$  corresponds
the symmetric bilinear form on $L$ with values in $\hat{\Bbb G}_m$,
$$\varphi_B(x,y)= {\roman{tr}}_{\Z_{p^2}/\Z_p}(xBy)\ \ .\tag10.12$$
The condition that an endomorphism $j$ of the special fibre ${\Cal
A}$ extends to the deformation corresponding to $\varphi:L'\to
L\otimes {\Bbb G}_m$ (cf.\ (10.10)) is
$$j\circ \varphi = \varphi\circ j\ \ ,\tag10.13$$
where $j$ denotes the induced maps on $L$ resp.\ $L'$. For a
deformation lying in ${\Cal T}$ and corresponding to $B(e\otimes
e)$ this translates into
$$\varphi_B(j(z_1), z_2)= \varphi_B(z_1, j^*(z_2)),\ \ z_1,
z_2\in L\ \ .\tag10.14$$ Let us now consider a special endomorphism
$j$ of $(A,\lambda,
\iota, \overline{\eta}^p)$. Then $j$ induces an endomorphism of
$L$ which sends
$$xe\longmapsto A\cdot x^{\sigma}e\ \ ,\tag10.15$$
for some well-defined scalar $A\in \Z_{p^2}$. Then since $j^*=j$,
the locus in ${\Cal T}$ where $j$ deforms is given by the condition
$$\varphi_B(Ax^{\sigma},y)= \varphi_B(x,Ay^{\sigma})\ \ ,\ \
\text{i.e.}$$
$$(AB)^{\sigma}= AB\ \ ,\ \ \text{i.e.}\ \ AB\in \hat{\Bbb
G}_m\ \ .\tag10.16$$ Let us write
$$B=a+b\delta\ \ ,\ \ A=\alpha+\beta\delta\ \ \text{with}\
a,b\in\hat{\Bbb G}_m\ \ ,\ \ \alpha,\beta\in\Z_p\ \ .$$ Then the
identity (10.16) is equivalent to
$$a\beta +b\alpha =0\ \ .\tag10.17$$
Let us choose multiplicative coordinates on the formal torus ${\Cal
T}$ dual to the entries $a$ and $b$ of $B$. Then, by (10.17) the
locus inside $\CT$ to which $j$ deforms is given by the equation
$$q_b^{\alpha}=q_a^{-\beta}\ \ ,\tag10.18$$
which we wish to consider locally around the base point
$q_a=q_b=1$.

\proclaim{Proposition 10.2}
Let $t\in\Z_{(p)}$ with ${\roman{ord}}_p(t)=0$. Then
$\CZ(t,\o)^{\roman{ord}}$ is a smooth relative curve over
${\roman{Spec}}\, \Z_{(p)}$.
\endproclaim
\demo{Proof}
It suffices to check this locally around a point
$(A,\lambda,\iota,\overline{\eta}^p)\in
\CZ(t,\o)^{\roman{ord}}({\Bbb F})$. Let us introduce local
coordinates $T_a$ and $T_b$ around the origin of $\CT$,
$$q_a=1+T_a\ \ ,\ \ q_b=1+T_b\ \ .\tag10.19$$
Then the leading term in the deformation equation (10.18) of $j$
with $Q(j)=t$ is
$$\alpha T_b+\beta T_a=0\ \ .\tag10.20$$
On the other hand,
$$Q(j)={\roman{Nm}}\, A=\alpha^2-\Delta\beta^2\ \ .\tag10.21$$
Hence the hypothesis implies that (10.20) defines a non-trivial
linear equation in the tangent space of $\CT$ and the claim
follows.
\qed
\enddemo

\proclaim{Corollary 10.3}
Let $t\in\Z_{(p)}$ with ${\roman{ord}}_pt =2s\ge 0$. Then
$\CZ(t,\o)^{\roman{ord}}$ is a relative divisor over
${\roman{Spec}}\, \Z_{(p)}$ of the form
$p^{2s}\cdot\CZ(t_0,\o)^{\roman{ord}}$ where
$\CZ(t_0,\o)^{\roman{ord}}$ is a smooth relative divisor over
${\roman{Spec}}\, \Z_{(p)}$.
\endproclaim
\demo{Proof}
Let $(\alpha,\beta)\in\Z_p^2$ be the vector corresponding to the
special endomorphism $j$ with $Q(j)=t$. By assumption $t=p^{2s}t_0$
with ${\roman{ord}}_pt_0=0$ and hence
$$(\alpha,\beta)= p^s\cdot (\alpha_0, \beta_0)\ \ ,\ \ (\alpha_0,
\beta_0)\in \Z_p^2\ \ .\tag10.22$$
By the previous Proposition, the locus inside $\CT$ corresponding
to $(\alpha_0,\beta_0)$ is a formal subtorus of relative dimension
1 of $\CT$. The locus inside $\CT$ corresponding to
$(\alpha,\beta)$ is given by (10.18), i.e.\ is the inverse image of
the locus corresponding to $(\alpha_0,\beta_0)$ under the isogeny
of degree $p^{2s}$ given by multiplication by $p^s$,
$$p^s:\CT\longrightarrow \CT\ \ .\tag10.23$$
The result follows.
\qed
\enddemo

We next consider the intersection of special divisors. By Corollary
3.8 we know that for $n\ge 3$
$\overline{\CZ}(T,\o)^{\roman{ord}}=\emptyset$ and that for $n=2$
the cycle $\overline{\CZ}(T,\o)^{\roman{ord}}$ consists of finitely
many points which are all in one and the same isogeny class of type
$I_1$. Let $T\in{\roman{Sym}}_2(\Z_{(p)})_{>0}$ be
$GL_2(\Z_p)$-equivalent to
$$T'= {\roman{diag}}\,(\varepsilon_1p^{2s_1},
\varepsilon_2p^{2s_2})\ \ .\tag10.24$$
Let $j_1,j_2$ be special endomorphisms of
$(A,\lambda,\iota,\overline{\eta}^p)\in {\Cal M}^{\roman{ord}}(\Bbb
F)$ with $Q(j_1,j_2)=T'$. As before, we associate to $j_1$ and $j_2$
scalars in $\Z_{p^2}$,
$$A_1=\alpha_1+\beta_1\delta\ \ ,\ \ A_2=\alpha_2+\beta_2\delta\ \
.\tag10.25$$
The condition $j_1j_2=-j_2j_1$ translates into
$$A_1A_2^{\sigma}=-A_2A_1^{\sigma}\Leftrightarrow
\alpha_1\alpha_2-\Delta\beta_1\beta_2=0\ \ .\tag10.26$$

Let us first assume that $s_1=s_2=0$. Then the local equations for
$\CZ(T,\o)$ in $\CT$ in terms of the multiplicative coordinates
$q_a,q_b$ resp.\ the local coordinates $T_a,T_b$ at the base point
are given by
$$\alpha_1T_a+\beta_1T_b=0\ \ ,\ \ \alpha_2T_a+\beta_2T_b=0\ \
.\tag10.27$$
The determinant of this system of linear equations is equal to
$\alpha_1\beta_2-\alpha_2\beta_1$. This is a unit because otherwise
$A_1A_2^{\sigma}$ would be divisible by $p$ which is excluded by
our hypothesis. It follows that the determinant is a unit, hence
the loci inside $\CT$ where $j_1$ resp.\ $j_2$ deform intersect
transversally.

Combined with Corollary 10.3 this gives the following statement.
Here we content ourselves with a statement concerning the special
fibre.

\proclaim{Proposition 10.4}
Let $T\in {\roman{Sym}}_2(\Z_{(p)})_{>0}$. Then
$\overline{\CZ}(T,\o)^{\roman{ord}}$ is an artinian scheme of
length $p^{{\roman{ord}}_p({\roman{det}} T)}$ at each of its
points.
\qed
\endproclaim

\proclaim{Corollary 10.5}
Assume that $\xi\in\overline{\CZ}(t_1,\o_1)\times_{\overline{\Cal
M}}\overline{\CZ}(t_2,\o_2)({\Bbb F})$ is an isolated ordinary
point. Then the intersection multiplicity of the divisors
$\overline{\CZ}(t_1,\o_1)^{\roman{ord}}$ and
$\overline{\CZ}(t_2,\o_2)^{\roman{ord}}$ at $\xi$ is equal to
$$e(\xi)=p^{{\roman{ord}}_p({\roman{det}}\, T_{\xi})}\ \ ,$$
where $T_{\xi}\in {\roman{Sym}}_2(\Z_{(p)})_{>0}$ is the fundamental
matrix of $\xi$, comp.\ end of section 6.
\qed
\endproclaim

{\bf Remark.} Let $t\in\Z_{(p)}$ with ${\roman{ord}}_pt$ even and
$>0$. Combining Theorem 6.1 with Corollary 10.3 we see that, outside of the
supersingular locus, the
corresponding cycle $\overline{\CZ}(t,\o)$ has  a $p$-power multiplicity and furthermore
$\overline{\CZ}(t,\o)$ contains some projective lines inside the supersingular locus. If
${\roman{ord}}_pt=0$, then,
outside the
supersingular locus,  $\overline{\CZ}(t,\o)$ is a smooth divisor and
$\overline{\CZ}(t,\o)$ cuts the supersingular
locus in finitely many superspecial points.

We conclude this section with one remark of a global nature.

\proclaim{Proposition 10.6}
Let $T\in{\roman{Sym}}_n(\Z_{(p)})_{>0}$. The special cycle
$\CZ(T,\o)$ is proper over ${\roman{Spec}}\, \Z_{(p)}$, unless

\roster
\item"{a)}" $B=M_2(\kay)$
\item"{b)}" $n=1$
\item"{c)}" $T=t\in {\roman{Nm}}\ \kay^{\times}$.
\endroster
\endproclaim

\demo{Proof}
If $B$ is a division algebra, then ${\Cal M}$ is proper over
${\roman{Spec}}\, \Z_{(p)}$ and hence so is $\CZ(T,\o)$ which maps
by a finite morphism to ${\Cal M}$. Hence we may assume a) and
hence
$$C(V)\otimes \kay\simeq M_4(\kay)\ \ .\tag10.28$$
Therefore, if $(A,\lambda,\iota,\overline{\eta}^p)$ is a point of
${\Cal M}$ the abelian variety $A$ is up to isogeny a product
$A_0^4$ where $A_0$ is an abelian variety of dimension 2 equipped
with an action
$$i_0:\kay\hookrightarrow {\roman{End}}^0(A_0)\ \ .\tag10.29$$
Furthermore the polarization $\lambda$ defines a $\kay$-linear
polarization $\lambda_0$ on $A_0$. A special endomorphism $j$ of
$(A,\lambda,\iota,\overline{\eta}^p)$ induces an endomorphism
$j_0:A_0\to A_0$ up to isogeny with
$$j_0^*=j_0\ \  \text{and}\ \ j_0\circ \iota_0(a)=
\iota_0(a^{\sigma})\circ j_0\ ,\  a\in \kay\ \ .\tag10.30$$ We now wish to check the valuative criterion of
properness. Thus we assume that the point $(A,\lambda,\iota,
\overline{\eta}^p)$ of ${\Cal M}$ is over the generic point of a
complete discrete valuation ring with algebraically closed residue
field. We also assume that $A$ and hence $A_0$ has semi-stable {\it
bad}, i.e., multiplicative reduction. We wish to show that b) and
c) hold. We may write the associated rigid-analytic abelian variety
as a quotient \cite{\raynaud}
$$A_0^{an}=L\otimes {\Bbb G}_m^{an} / i(L')\ \ .\tag10.31$$
The polarization $\lambda_0$ identifies $L'\otimes\Q$ with the dual
of $L\otimes\Q$ and the embedding $i:L'\to L\otimes {\Bbb
G}_m^{an}$ defines a pairing
$$L'\times L^*\longrightarrow {\Bbb
G}_m^{an}\buildrel{\roman{val}}\over\longrightarrow\Z$$ which,
combined with the previous identification, defines a {\it
symmetric} bilinear form
$$(\ ,\ )\ :\  L\otimes \Q\times L\otimes
\Q\longrightarrow\Q\tag10.32$$
which is {\it positive-definite} and satisfies
$$(\iota_0(a)x,y)=(x,\iota_0(a)y)\ \ ,\ \ a\in \kay\ \ .$$
After choosing a generator of $L\otimes \Q$ the form (10.32) is
therefore of the form
$$(x,y)={\roman{tr}}(bxy)\tag10.33$$
for a totally positive element $b$ of $k$.

A special endomorphism $j$ defines an endomorphism of $L$ which
sends $x$ to $a\cdot x^{\sigma}$ for a fixed scalar $a\in \kay$ and
which satisfies
$$(ab)^{\sigma} =ab\ \ ,\tag10.34$$
comp.\ (10.14). But equation (10.34) determines $a$ up to a scalar
in $\Q$, hence there cannot be a non-trivial special endomorphism
of $(A,\lambda,\iota,\overline{\eta}^p)$ anticommuting with $j$.
Hence b) is proved and c) follows from the identity
$$Q(j)={\roman{Nm}}(a)\ \ .\qquad\qed$$
\enddemo

The results given in this section can be regarded as first steps
towards preparing the stage for a theory of Hilbert-Blumenthal
surfaces in characteristic $p$, analogous to the classical theory
of Hirzebruch and Zagier, e.g.\ \cite{\hirzebruch},
\cite{\vandergeer}. Namely, one might ask whether the intersection
numbers of the special cycles $\overline{\CZ}(t,\o)$ (for
$t\in\Z_{(p)>0}$) are related to the Fourier coefficients of some
sort of modular form. For this a better understanding of the
singularities of $\overline{\CZ}(t,\o)$ at supersingular points
seems to be necessary. Similarly, one might investigate whether the
special curves $\overline{\CZ}(t,\o)$ can be used to determine the
position of $\overline{\Cal M}$ (suitably compactified if $B\simeq
M_2(\kay))$ in the classification of algebraic surfaces. For this
it would seem necessary to better understand the irreducible
components of special curves. One contrast to the situation in
characteristic zero is that on $\overline{\Cal M}$ it is easy to
find effective ample divisors (to be on the safe side, assume $B$
to be a division algebra). Indeed the supersingular locus
$\overline{\Cal M}^{ss}$ splits into a sum of two divisors (sum of
the {\it even} resp.\ {\it odd} projective lines ${\Bbb P}_L^1$)
which are exchanged by the Frobenius in ${\roman{Gal}}({\Bbb
F}_{p^2}
/{\Bbb F}_p)$ and each of which is ample. Indeed they span a
2-dimensional submodule of the $\ell$-adic cohomology of
$\overline{\Cal M}\times_{{\roman{Spec}}\, {\Bbb F}_p}
{\roman{Spec}}\, {\Bbb F}$ which is stable under the actions of the
Galois group and of the Hecke algebra prime to $p$; by \cite{\hlr}
this submodule must therefore lie in the subspace cut out by a
one-dimensional automorphic representation of $G({\Bbb A})$ and all
effective divisor classes in such a subspace are ample.


\subheading{\Sec11. The case of a split prime}

In this section we briefly sketch the case when $p$ is split in
$\kay$, and we assume, as in section 2, that $B$ is split at each
prime $\wp$ of $\kay$ dividing $p$. Note that the case
$\kay=\Q\oplus\Q$ is allowed. For $p$ split in $\kay$,  ${\Cal
O}_C\otimes \Z_p\simeq M_4(\Z_p)$ as before, but now ${\Cal
O}_{\smallkay}\otimes \Z_p\simeq \Z_p\oplus \Z_p$.

Fix a point $\xi =(A, \lambda, \iota, \overline{\eta}^p)\in {\Cal
M}({\Bbb F})$. Let $A(p)$ be the $p$-divisible group of $A$. The
action of ${\Cal O}_C\otimes \Z_p\simeq M_4(\Z_p)$ on $A(p)$ allows
us to write as before $A(p)={\Cal A}^4$. On the $p$-divisible group
${\Cal A}$ there is an action of ${\Cal O}_{\kay}\otimes \Z_p\simeq
\Z_p\oplus\Z_p$ which allows us to write further
$${\Cal A}={\Cal A}_1\times {\Cal A}_2\ \ ,\tag11.1$$
where ${\Cal A}_1$ and ${\Cal A}_2$ are $p$-divisible groups of
dimension 1 and height 2. The principal
quasi-polarization of ${\Cal A}$ induced by $\lambda$ comes from a
pair of principal quasi-polarizations $\lambda_1$ and $\lambda_2$
of ${\Cal A}_1$ and ${\Cal A}_2$. By the Serre-Tate theorem,
deforming $\xi$ is equivalent to deforming $({\Cal A}_1,
\lambda_1)$ and $({\Cal A}_2, \lambda_2)$. Since $\lambda_i$
deforms automatically with ${\Cal A}_i$, the deformation space of
$({\Cal A}_i, \lambda_i)$ is the same as that of ${\Cal A}_i$,
i.e.\ of a suitable point on the moduli space of elliptic curves.
In particular, we deduce the following result.

\proclaim{Proposition 11.1}
The supersingular locus ${\Cal M}^{ss}$ consists of a finite number
of isolated points.
\qed
\endproclaim

We next give an estimate on the space $V_{\xi}$ of special
endomorphisms of $\xi=(A,\lambda,\iota, \overline{\eta}^p)$. We
have an inclusion
$$V_{\xi}\subset \{ j:{\Cal A}\to {\Cal A},\ {\roman{deg}}\, j=1,\
j^*=j\}\ \ .\tag11.2$$ Here ${\roman{deg}}\, j$ refers to the
$\Z/2$-grading of ${\Cal A}={\Cal A}_1\times {\Cal A}_2$. Therefore
$j$ corresponds to a pair of homomorphisms of $p$-divisible groups
$$j_1: {\Cal A}_1\to {\Cal A}_2\ \ ,\ \ j_2: {\Cal A}_2\to {\Cal
A}_1\ \ .\tag11.3$$ Using the polarizations $\lambda_1$ and
$\lambda_2$, the condition $j^*=j$ is equivalent to
$$j_2= j_1^*\ \ ,\ \ j_1=j_2^*\ \ .\tag11.4$$
The quadratic form is given by the degree of $j_1$ (or $j_2$),
$$Q(j)\cdot {\roman{id}}_{{\Cal A}_1}= j_1^*j_1\ \ .\tag11.5$$
Decomposing ${\Cal A}_1$ and ${\Cal A}_2$ up to isogeny into simple
$p$-divisible groups we now obtain the following statement.
\proclaim{Proposition 11.2}
We have case by case for the rank of the $\Z_{(p)}$-module
$V_{\xi}$ of special endomorphisms of $\xi=(A, \lambda, \iota,
\overline{\eta}^p)$.

\roster
\item"{(i)}" If ${\Cal A}_1$ and ${\Cal A}_2$ are both
supersingular, then ${\roman{rk}}\, V_{\xi}=4$.
\item"{(ii)}" If ${\Cal A}_1$ and ${\Cal A}_2$ are both ordinary,
then ${\roman{rk}}\, V_{\xi}\leq 2$.
\item"{(iii)}" If one of ${\Cal A}_1, {\Cal A}_2$ is supersingular
and the other ordinary, then $V_{\xi}=(0)$. \qed
\endroster
\endproclaim

We now fix $T\in {\roman{Sym}}_3(\Z_{(p)})_{>0}$ and $\omega\subset
V({\Bbb A}_f^p)^3$ and consider the special cycle ${\Cal
Z}(T,\omega)$. By the previous propositions we know that ${\Cal
Z}(T,\omega)$ is a finite set of points lying over the
supersingular set ${\Cal M}^{ss}$. Again invoking the Serre-Tate
theorem we have an identification for $\xi= (A,\lambda, \iota,
\overline{\eta}^p, {\j})\in {\Cal Z}(T,\omega)$
$$\hat{\Cal Z}(T,\omega)_{\xi}={\roman{Def}}({\Cal A}_1, {\Cal
A}_2;\j_1)={\roman{Def}}({\Cal A}_1, {\Cal A}_2;M_1)\ \ .\tag11.6$$
Here $\j_1$ is the triple consisting of the ``first component''
$j_1:{\Cal A}_1\to {\Cal A}_2$ of the three members $j$ of $\j$ and
$M_1$ is the $\Z_p$-module spanned by them. The length of the
Artinian scheme on the RHS of (11.6) was determined by Gross and
Keating \cite{\grosskeating} and we therefore obtain the following
analogue of Proposition 6.2.

\proclaim{Proposition 11.3}
The length of the local ring ${\Cal O}_{{\Cal Z}(T,\omega),\xi}$
only depends on the $GL_3(\Z_p)$-equivalence class of $T$ and is
equal to $e_p(T)$, given in Proposition~6.2 above.
\qed
\endproclaim

As in section 7 we introduce in the setting of (6.19)-(6.23)
$$\langle {\Cal Z}(T_1, \omega_1),\ldots, {\Cal Z}(T_r,
\omega_r)\rangle_p^{\roman{proper}}=\sum_{\xi} e(\xi)\ \ .\tag11.7$$
Here $\xi$ runs over the isolated points of ${\Cal
Z}(T_1,\omega_1)\times_{\Cal M}\ldots\times_{\Cal M}{\Cal Z}(T_r,
\omega_r)$ which by Propositions 11.1 and 11.2 are the points
$\xi$ where the fundamental matrix $T_{\xi}\in {\roman{Sym}}_3(\Q)$
is non-singular. In the special case when $r=1$, the cycle ${\Cal
Z}(T,\omega)$ lies over the supersingular locus and is finite
(without any further condition on the divisibility by $p$ of $T$).
In this case we put $\langle {\Cal
Z}(T,\omega)\rangle_p=\langle{\Cal
Z}(T,\omega)\rangle_p^{\roman{proper}}$, i.e.
$$\langle{\Cal Z}(T,\omega)\rangle_p =\sum\limits_{\xi\in{\Cal
Z}(T,\omega)}e(\xi)\ \ .\tag11.8$$ Appealing to Proposition 11.3 we
have that (11.8) equals
$$\langle{\Cal Z}(T,\omega)\rangle_p =e_p(T)\cdot \vert{\Cal
Z}(T,\omega)({\Bbb F})\vert\ \ .\tag11.9$$

As in the inert case we
introduce the twisted quadratic space $V'$ over $\Q$ which is
positive-definite, is isomorphic to $V$ at all finite places $\neq
p$ and at $p$ has the same determinant as $V(\Q_p)$ but opposite
Hasse invariant. The quaternion algebra $B'=C^+(V')$ over $\kay$ is
ramified at the two infinite places, is isomorphic to $B$ at all
finite places not over $p$ and is ramified at the two places
$\wp_1$ and $\wp_2$ over $p$. Let $G'$ be the corresponding inner
form of $G$. Then we have a bijection
$${\Cal M}^{ss}({\Bbb F})= G'(\Q)\setminus \bigg( G'(\Q_p)/K'_p\times
G({\Bbb A}_f^p)/K^p\bigg)\ \ .\tag11.10$$ Here $K'_p$ is the unique
maximal compact subgroup of $G'(\Q_p)= B_{\wp_1}^{\prime,\times}\times
B_{\wp_2}^{\prime,\times}$.

After identifying $V'(\Q_p)$ with the ramified quaternion algebra
${\Bbb B}_p$ over $\Q_p$ equipped with its norm form, we have a
natural lattice $V'(\Z_p)$ in $V'(\Q_p)$, namely the ring of
integers in ${\Bbb B}_p$. The usual procedure now yields the
following expression for the cardinality of ${\Cal
Z}(T,\omega)({\Bbb F})$, comp.\ Proposition 7.1.

\proclaim{Proposition 11.4}
Let $K'= K'_p\cdot K^p\subset G'({\Bbb A}_f)$. Let
$$\varphi_f^p= {\roman{char}}(\omega)\ \ ,\ \
\varphi'_p={\roman{char}}(V'(\Z_p)^3)$$
and
$$\varphi'_f=\varphi'_p\cdot \varphi_f^p\in S(V'({\Bbb
A}_f)^3)^{K'}\ \ .$$ Then
$$\vert{\Cal Z}(T,\omega)({\Bbb F})\vert= \vol(K')^{-1}\cdot
{\roman{vol}}(Z'(\Q)\setminus Z'({\Bbb A}_f))
\cdot O_T(\varphi'_f)\
\ ,$$
where the orbital integral on the right is defined as in
Proposition 7.1.
\endproclaim

The relationship to the Eisenstein series is now
established just as in section 7 above. We define $\P_p$
and $\P_p'$ as in (7.11) and complete $\P_p$ into an incoherent
standard section as before. Then we obtain formula (7.13) for the
present case. The next proposition gives the values and derivatives
of the relevant Whittaker functions at $s=0$.

\proclaim{Proposition 11.5} Let $T\in \Sym_3(\Q_p)$ with $\det(T)\ne 0$.
\hfill\break
(i) If $W_{T,p}(e,0,\P_p')\ne 0$, then $T\in \Sym_3(\Z_p)$. \hfill\break
(ii) If $T\in\Sym_3(\Z_p)$ is represented by $V'(\Q_p)$, then
$$W_{T,p}(e,0,\P_p) = \gamma(V'_p)\cdot 2 p^{-4} (p+1)^2.$$
(iii) If $T\in\Sym_3(\Z_p)$ is represented by $V'(\Q_p)$, then
$$W_{T,p}'(e,0,\P_p) = \gamma(V_p)\cdot \log(p)\cdot (1-p^{-2})^2\cdot e_p(T).$$
\endproclaim
\demo{Proof}
We again use the relation (7.15) between the Whittaker functional
and the representation densities and obtain the analogue of (7.17)
$$W_{T,p}'(e,0,\P_p) = -\log(p) \cdot \gamma(V_p)
\frac{\partial\ }{\partial X}\bigg\{ A_{S,T}(X)\bigg\}\bigg\vert_{X=1}.$$
Here $S$ is the matrix for the quadratic form on $V(\Q_p)$, i.e.,
$$S=\diag(1,-1,1,-1)= H_4,$$
a hyperbolic form of dimension $4$. By the same argument which
yields Corollary~10.6 of \cite{\krsiegel}, we obtain
$$\frac{\partial\ }{\partial X}\bigg\{ A_{S,T}(X)\bigg\}\bigg\vert_{X=1}
= -(1-p^{-2})^2\cdot e_p(T)$$
where the factor $(1-p^{-4})$ in Corollary~10.6 of \cite{\krsiegel}
has been changed to $(1-p^{-2})$ since the reduction formula (10.11)
of \cite{\krsiegel} and hence the factor $(1+p^{-2})$ of (10.14) there does
not arise in our present situation. This proves (iii).

Finally, (ii) follows from Proposition~6.10 of \cite{\grosskeating},
cf. also Lemma~7.10 of \cite{\krsiegel}, where we note that
$$S'=\diag(1,-\b,-p,p\b)$$
has determinant $p^2$, and that
$$W_{T,p}(e,0,\P_p') = \gamma(V_p')\cdot |\det(S')|^{3/2}\cdot \a_p(S',T).$$
\qed\enddemo

Using these values in (7.13) we obtain the analogue of Theorem~7.3 in
the case of a split prime.
\proclaim{Theorem 11.5} Let $T\in\Sym_3(\Q)_{>0}$ be
represented by $V(\A_f^p)$, but not by $V(\Q_p)$. Also
assume that $\o$ is locally centrally symmetric. \hfill\break
(i) If $T\notin \Sym_3(\Z_{(p)})$, then $\CZ(T,\o)$ is empty
and $E'_T(h,0,\P)=0$.
\hfill\break
(ii) For any such $T\in \Sym_3(\Z_{(p)})_{>0}$,
the cycle $\CZ(T,\o)$ is either empty or zero dimensional and
$$E'_T(h,0,\P) = -\frac12 \vol(SO(V')(\R)\cdot \pr(K))\cdot W_T^2(h)
\cdot \log(p)\cdot \langle{\Cal Z}(T,\omega)\rangle_p.$$
\endproclaim
\demo{Proof}
Using Proposition~11.4 and (7.13), we have
$$\align
E'_T(h,0,\P)&= \vol(SO(V')(\R)\cdot \pr(K'))\cdot W_T^2(h)\\
\nass
{}&\qquad \times \bigg( -\frac12 \log(p)\cdot (p-1)^2 \cdot e_p(T)\cdot
|\CZ(T,\o)(\Bbb F)|\bigg)\\
\nass
{}&=-\frac12 \vol(SO(V')(\R)\cdot \pr(K))\cdot W_T^2(h)
\cdot \langle{\Cal Z}(T,\omega)\rangle_p,
\endalign
$$
where we have used the fact that
$$\frac{\vol(K')}{\vol(K)} = (p-1)^2.$$
\qed\enddemo


\subheading{\Sec12. A global model and generating series}

In this section we formulate a moduli problem which is represented by a scheme over $\Spec \Z\left[\frac1{N}\right]$
for some specified integer $N$ and whose base change to $\Spec \Z_{(p)}$,  for any $p\nmid N$, coincides with the moduli
space introduced in section 2. We similarly extend the special cycles. We then use the
results of sections 7 and 11 to give a partial identification of the generating series for the arithmetic
degrees of these special cycles.

As in section 1, we start with the quadratic space $(V,Q)$ over $\Q$ of signature $(2,2)$ and
the element $\tau\in B = C^+(V)$ with $\tau=-\tau^\iota$ and $\tau^2<0$. On the left
$C(V)$-module $U=C(V)$, we have the nondegenerate alternating form
$$ <x,y>\ =\ \tr^0(y^\iota \tau x). \tag 12.1$$
We also have the open compact subgroup $K$ of $G(\A_f)$ which, as always, is assumed to
be sufficiently small, i.e., the image $\pr(K)\subset SO(V)(\A_f)$ is neat.

Fix a $\Z$-lattice $\Lambda$ in $V(\Q)$. Let $\OC=C(\Lambda)$ be its Clifford algebra
and let $U_\Z=\OC$ be the corresponding lattice in $U(\Q)$. There exists an integer $N$
such that for any prime $p$ with $p\nmid N$ all of the following conditions are satisfied.
\roster
\item"{(i)}" $p\ne 2$.
\item"{(ii)}" $\Lambda\tt\Z_p$ is a self dual $\Z_p$-lattice in $V(\Q_p)$.
\item"{(iii)}" $\tau\in (\OC\tt\Z_p)^\times$ and the form $<\ ,\ >$ is
perfect on the lattice $U(\Z_p)=U_\Z\tt\Z_p$. In particular, $\OC\tt\Z_p$
is invariant under the involution $x\mapsto x^*=\tau x^\iota \tau^{-1}$.
\item"{(iv)}" $K=K_pK^p$ with $K^p\subset G(\A_f^p)$ and
$$K_p= \GSpin(V)(\Q_p) \cap (\OC\tt\Z_p)^\times.$$
\endroster

We now fix such an integer $N$ and formulate a moduli problem over
$\Spec \Z\left[\frac1{N}\right]$. To do so, we introduce the
category of abelian varieties {\it up to $N$-primary isogeny.} As
in the case of the cateory of abelian varieties up to isogeny prime
to $p$ (as used in section 2), this category is obtained from that
of abelian varieties by formally inverting all isogenies whose
degrees involve only prime factors of $N$. The moduli problem
associates to $S\in \text{Sch}/\Z\left[\frac1{N}\right]$ the set of
isomorphism classes of $4$-tuples $(A,\lambda,\iota,
\bar\eta_N)$ where
\roster
\item"{(i)}" $A$ is an abelian scheme over $S$ up to $N$-primary isogeny.
\item"{(ii)}" $\lambda:A\rightarrow \hat A$ is a $\Z\left[\frac1{N}\right]^\times$-class
of principal polarizations on $A$. 
\item"{(iii)}" $\iota:\OC\tt\OK \lra \End(A)\tt\Z\left[\frac1{N}\right]$ is a homomorphism such that
$$\iota(c\tt a)^* = \iota(c^*\tt a),$$
for the Rosati involution of $\End^0(A)$ determined by $\lambda$
and the involution $*$ of $C(V)$.
\item"{(iv)}" $\bar\eta_N$ is a $K_N$-equivalence class of $\OC\tt\OK$-linear isomorphisms
$$\eta_N:V_N(A)\isoarrow U(\A_N)$$
which preserves the symplectic forms up to a scalar in $\A_n^\times$.
\endroster

Here we have used the notation $\A_N = \prod_{p\mid  N} \Q_p$ and
$V_N(A) =\prod_{p\mid  N} V_p(A)$ and have written $K$ as
$$K = K_N\cdot \prod_{p\mid N} K_p,$$
cf. (iv) in the conditions on $N$ above. As before we impose the determinant condition (1.1).

\proclaim{Theorem 12.1} The moduli problem defined above is representable
by a quasi-projective scheme $\Cal M$ over $\Z\left[\frac1{N}\right]$ whose base change
to $\Spec \Z_{(p)}$ is the moduli scheme introduced in section 2 for any $p\nmid N$.
In particular, its generic fiber is the moduli scheme $M$ introduced in section 1.
\endproclaim

We shall content ourselves with showing how to establish a bijection
$\Cal M(\Spec k) \simeq M(\Spec k)$ for an algebraically closed field
$k$ of characteristic $0$. Let $(A,\lambda,\iota, \bar{\eta})\in  M(k)$ and
fix $\eta\in\bar{\eta}$. We obtain a $\hat{\Z}$-lattice in $\hat{V}(A)$
defined by
$$T_\eta(A) := \eta^{-1}(U(\hat \Z)),\tag 12.2$$
where $U(\hat \Z) = U_\Z\tt\hat \Z$. Then, by our choice of $N$,
for any $p\nmid N$, the lattice $T_\eta(A)\tt\Z_p$ in $V_p(A)$ is the
unique $\OC\tt\OK$-stable lattice which is self-dual for the form
$<\ ,\ >_\eta\ = \eta^*(<\ ,\ >)$ induced by $\eta$. In particular, for all $p\nmid N$, $T_\eta(A)\tt\Z_p$ is
independent of the choice of $\eta$.  
Let $B$ be the abelian variety in the isogeny class of $A$ with
$$\hat{T}(B) = T_\eta(A). \tag 12.3$$
Then $B$ is equipped with an action of $\OC\tt\OK$ and,
by the above, it is unique up to $N$-primary isogeny. 
Furthermore, there exists a polarization $\lambda_0\in \lambda$ and a
trivialization of the roots of unity of $k$ such that, for the pairing
$<\ ,\ >_{\l_0}$ on $\hat V(A)$ thus determined,
$$<\ ,\ >_{\lambda_0} \ = \ <\ ,\ >_\eta. \tag12.4$$
Then $\lambda_0$ defines a polarization of $B$ which is principal in the
category of abelian varieties varieties up to $N$-primary isogeny.
Finally, $\bar{\eta}_N$ is obtained by simply ignoring the
components of $\bar{\eta}$ prime to $N$. Then
$(B,\lambda_0,\iota,\bar{\eta}_N)\in \Cal M(k)$ and this procedure
defines the map $M(k)\rightarrow \Cal M(k)$. The map in the opposite direction
is constructed in a similar way.

We next turn to special cycles. The definition of special
endomorphism, Definition~1.2, carries over in the obvious way to
$(A,\lambda,\iota, \bar{\eta}_N)\in \Cal M(S)$. If $S$ is
connected, the special endomorphisms form a quadratic space over
$\Z\left[\frac1{N}\right]$. For the remainder of this section, we
fix a compact open subset $\o_N\subset V(\A_N)^3$, stable under the
action of $K_N$. For any $T\in
\Sym_3(\Z\left[\frac1{N}\right])_{>0}$, we define $\Cal Z(T) =\Cal
Z(T,\o_N)$ as the scheme representing the functor of $5$-tuples
$(A,\lambda,\iota,\bar\eta_N;\j)$, where $\xi=
(A,\lambda,\iota,\bar\eta_N)\in \Cal M(S)$, and $\j=(j_1,j_2,j_3)$
is a triple of special endomorphisms with $Q(\j)=T$ and where
$\eta\circ \j\circ \eta^{-1}\in \o_N$ for one, and hence for all,
$\eta\in\bar\eta_N$. We note that the generic fiber of $\Cal
Z(T,\o_N)$ is the special cycle associated in section 1 to $T$ and
$$\o=\o_N\times\prod_{p\nmid N} \o_p,\tag12.5$$
where $\o_p =(\Lambda\tt\Z_p)^3$.
Similarly, the base change of $\Cal Z(T,\o_N)$ to $\Spec \Z_{(p)}$ for $p\nmid N$ is
the special cycle associated in section 2 to $T$ and
$$\o^p = \o_N\times\prod_{\scr\ell\nmid N\atop \ell\ne p}\o_\ell,\tag12.6$$
where $\o_\ell=(\Lambda\tt\Z_\ell)^3$.

We next introduce the Eisenstein series associated to $\o_N$. Let
$\ph_f\in S(V(\A_f)^3)$ be the characteristic function of the set $\o$ defined in (12.5)
and let $\P_f$ be the standard
section of the induced representation with $\P_f(0)=\l_f(\ph_f)$,
where $\l_f:S(V(\A_f)^3) \rightarrow I_f(0,\chi)$.
We again complete $\P_f$ into an incoherent section
$\P = \P_\infty\cdot \P_f$ with $\P_\infty =\P_\infty^2$ as usual and we form the
corresponding Eisenstein series $E(g,s,\P)$.

For $\tau=u+iv\in \frak H_3$, the Siegel half space of genus $3$, we let
$$h_\tau = \pmatrix 1&u\\{}&1\endpmatrix\pmatrix a&{}\\{}&{}^ta^{-1}\endpmatrix\in Sp_6(\R)\tag12.7$$
where $a\in GL_3(\R)$ with $\det(a)>0$ and $v=a{}^ta$. Note that $h_\tau(i)=\tau\in \frak H_3$. Then
$$\phi(\tau):= \det(v)^{-1}\cdot E'(h_\tau,0,\P)\tag12.8$$
is a (non-holomorphic) Siegel modular form of weight $2$ with respect to
the arithmetic subgroup $\Gamma_H=Sp_6(\Q)\cap K_H$ where $K_H\subset Sp_6(\A_f)$ is
an open compact subgroup fixing $\P_f$ in the induced representation.

It will be convenient to introduce the following definition.

{\bf Definition 12.2.} An element $T\in \Sym_3(\Z\left[\frac1{N}\right])_{>0}$
is {\it regular} if (i) $\text{Diff}(T) =\{p\}$, i.e., if $T$ is represented by  $V(\Q_\ell)$
for all primes $\ell\ne p$ but is not represented by $V(\Q_p)$,
(ii) $p\nmid N$, and (iii) if $p$ is inert in $\kay$, then $p\nmid T$.

If $T$ is regular with associated prime $p$, we can introduce the degree
$$\degh( \Cal Z(T,\o_N) ) = \log(p)\cdot <\Cal Z(T,\o_N)\tt\Z_{(p)}>_p,\tag12.9$$
as in sections 7 and 11.

The next result provides a partial analogue for the case of Hilbert-Blumenthal surfaces
of the main results of \cite{\annals}, and \cite{\tiny}.

\proclaim{Theorem 12.3} The modular form $\phi(\tau)$ has a Fourier expansion
of the form
$$\align
C^{-1}\cdot \phi(\tau) &= \sum_{T\ \text{\rm regular}} \degh(\Cal Z(T,\o_N)) \cdot q^T
+\sum_{p<\infty} \sum_{\matrix\scr T\in \Sym_3(\Q)_{>0}\\\scr \Diff(T) =\{p\}\\\scr T\ \text{not regular}\endmatrix }
a_T\cdot q^T \\
\nass
\nass
\nass
{}&\qquad+\sum_{\matrix \scr T\in \Sym_3(\Q)\\ \scr\sig(T)=(2,1)\endmatrix } a_T(v)\cdot q^T
+ \sum_{\matrix \scr T\in \Sym_3(\Q)\\ \scr\det(T)=0\endmatrix } a_T(v)\cdot q^T ,\endalign
$$
where
$$C=-\frac12\cdot \vol(SO(V')(\R)\pr(K)),$$
and
where $q^T = e^{2\pi i \tr(T\tau)}$.
\endproclaim

\demo{Proof} First recall that the only nonsingular $T\in \Sym_3(\Q)$ which
contribute to the Fourier expansion of $\phi(\tau)$ are those $T$
for which $|\Diff(T)|=1$, \cite{\annals}, Theorem~6.1. Furthermore,
if $\Diff(T)=\{p\}$ for a finite prime $p$, then the associated
term in the Fourier expansion is holomorphic and hence has the form
$a_T\cdot q^T$. This follows from (iii) of Theorem~6.1 of
\cite{\annals}, since, in our present situation, the quadratic
space $V^{(p)}$ is positive definite and hence the theta integral
whose $T$-th Fourier coefficient contributes the $\tau$-dependence
to the $T$-th Fourier coefficient of $\phi$ is holomorphic. Suppose
that, in addition, $T$ is regular, so that $p\nmid N$ and, if $p$
is inert in $\kay$, then $p\nmid T$. Then Theorems~7.3 and~11.5
provide the claimed expression for the Fourier coefficient. Note
that
$$q^T = \det(v)^{-1}\cdot W_T^2(h_\tau).\tag12.10$$
\qed\enddemo

{\bf Remark 12.4.} Obviously one would like to extend this result to include an
arithmetic interpretation of all of the Fourier coefficients of $\phi(\tau)$.
For example, consider the terms in the second sum for $T$'s with $\Diff(T)=\{p\}$
with $p\nmid N$  but with $p$ inert in $\kay$ and $p\mid T$. For such $T$'s one would like
to identify $a_T$ with some kind of arithmetic degree of the {\it one dimensional}
cycle $\Cal Z(T,\o^p)$, where $\o^p$ is as in (12.6), whose structure is described
in section 8. For $p$ dividing $N$, one has to deal with problems of bad reduction.


\subheading{\Sec13. On the hereditary nature of special cycles}

In this series of papers \cite{\annals}, \cite{\krsiegel},
\cite{\krdrin}, \cite{\tiny} we have investigated special cycles on
Shimura varieties associated to orthogonal groups of signature
$(r,2)$ for $r=3,2,1,0$. In this section we make some remarks on
the relations among these cycles in the various cases.

Analytically speaking, this relation is fairly obvious. Quite
generally, let $(V,q)$ be a quadratic space of signature $(r,2)$
over $Q$. Let $V'\subset V$ be a subspace of signature $(r-n,2)$ of
the form $V'=<x>^{\perp}$ for some $x\in V({\Bbb Q})^n$. Denoting
by $G$ resp.\ $G'$ the groups of spinorial similitudes and by
${\Cal D}$ resp.\ ${\Cal D'}$ the spaces of oriented negative
2-planes of $V(\R)$ resp.\ $V'(\R)$ we obtain injections
$$G'\hookrightarrow G\ \ ,\ \ {\Cal D}\hookrightarrow {\Cal D}\ \
.\tag13.1$$
For $K'=K\cap G'({\Bbb A}_f)$ we also obtain a (generically 1-1) morphism of
Shimura varieties
$$Sh(G',{\Cal D})_{K'}\hookrightarrow Sh(G,{\Cal D})_K\ \
.\tag13.2$$
Under this map, a special cycle for $Sh(G',{\Cal D'})_{K'}$
defined as in \cite{\duke} is mapped into a certain disjoint sum of
connected components of some special cycles on $Sh(G,{\Cal D})_K$,
comp.\ Proposition 1.4 above. In this sense the special cycles are
hereditary for the various Shimura varieties, cf.~also section 9 of \cite{\duke}.

For $r=0,\ldots,3$ we have a moduli-theoretic definition of the
corresponding Shimura varieties. In each of these cases the special
cycles were defined by imposing {\it special endomorphisms} on the
abelian varieties with additional structure parametrized by the
model in question. However, in each of the cases the notion of
special endomorphisms depended on the moduli problem and the
hereditary nature of these notions is hardly transparent.

Our aim in this section is more modest. We only wish to compare the
infinitesimal deformation theory of a supersingular point of the
moduli problem with values in ${\Bbb F}$ and of the special cycles
passing through this point, for various $r$. By the theorem of
Serre-Tate this leads case by case to the following deformation
problems.

{\bf r=3} {\it (good reduction case):} ${\roman{Def}}({\Cal A},
\lambda)$, where ${\Cal A}$ is a supersingular formal group of
dimension 2 and height 4 and where $\lambda$ is a principal
quasi-polarization of ${\Cal A}$. A {\it special endomorphism} is
an element $j\in{\roman{End}}({\Cal A})$ such that $j^{\ast}=j$ and
${\roman{tr}}^0(j)=0$.

{\bf r=2} {\it (good reduction case):} ${\roman{Def}}({\Cal A},
\lambda,\iota)$, where ${\Cal A}$ is a supersingular formal group
of dimension 2 and height 4 and where $\iota:{\Cal
O}\to{\roman{End}}({\Cal A})$ is an action of ${\Cal O}=\Z_{p^2}$
or ${\Cal O}=\Z_p\oplus\Z_p$ satisfying the {\it special condition}
and where $\lambda$ is a principal quasi-polarization which is
$\iota$-linear. A {\it special endomorphism} is an element
$j\in{\roman{End}}({\Cal A})$ such that $j^*=j$ and
$j\circ\iota(a)= \iota(\overline a)\circ j$, $a\in{\Cal O}$. In
this case, automatically, ${\roman{tr}}^0(j)=0$.

{\bf r=1} {\it (good reduction case):} ${\roman{Def}}({\Cal A})$,
where ${\Cal A}$ is a supersingular formal group of dimension 1 and
height 2. A {\it special endomorphism} is an element
$j\in{\roman{End}}({\Cal A})$ such that ${\roman{tr}}^0(j)=0$.

{\bf r=1} {\it ($p$-adic uniformization case):}
${\roman{Def}}({\Cal A}, \iota)$, where $({\Cal A}, \iota)$ is a
s.f.\ $O_{B_p}$-module (here $B_p$ is the quaternion division
algebra over $\Q_p$). A {\it special endomorphism} is an element
$j\in {\roman{End}}({\Cal A}, \iota)$ such that
 ${\roman{tr}}^0(j)=0$.

{\bf r=0} {\it (good reduction case):} ${\roman{Def}}({\Cal
A},\iota)$ where ${\Cal A}$ is a supersingular formal group of
dimension 1 and height 2 and where $\iota: {\Cal O}\to
{\roman{End}}({\Cal A})$ satisfies the special condition. Here
${\Cal O}$ is the ring of integers in a quadratic field extension
of $\Q_p$. A {\it special endomorphism} is an element $j\in
{\roman{End}}({\Cal A})$ such that $j\circ \iota(a)=
\iota(\overline a)\circ j$, $a\in {\Cal O}$.

To illustrate the hereditary nature of these definitions we now
start with an object of the deformation problem for a given $r$ and
an $n$-tuple $\j=(j_1,\ldots, j_n)$ of special endomorphisms and
construct under certain hypotheses an object of the deformation
problem for $r-1$ and an $(n-1)$-tuple $\j'=(j'_1,\ldots, j'_n)$ of
special endomorphisms of the deformation problem for $r-1$.

{\bf r=3.} Suppose that we are given $({\Cal A},\lambda)$ and a collection $\j$
of special endomorphisms with $q(j_1)= 1$, or
$\Delta$. Then $j_1$ defines a {\it special} action $\iota:{\Cal
O}\to{\roman{End}}({\Cal A})$, where ${\Cal O}=\Z_p\oplus \Z_p$ or
${\Cal O}=\Z_{p^2}$, and $\j'=(j_2,\ldots, j_n)$ is an
$(n-1)$-tuple of special endomorphisms of $({\Cal A},\lambda,\iota)$.

{\bf r=2.} Suppose that we are given $({\Cal A}, \lambda,\iota)$ where we assume ${\Cal
O}=\Z_{p^2}$. A collection $\j$ of special endomorphisms with $q(j_1)=1$ (resp. $q(j_1)=\Delta$)
determines idempotents via (6.7) (resp. (6.13)) and hence an isomorphism
$$({\Cal A},\lambda)\simeq ({\Cal A}_1, \lambda_1)\times ({\Cal A}_2,
\lambda_2).$$
We also have a natural isomorphism from ${\Cal A}_1$ to
${\Cal A}_2$, namely $\underline{\delta}$ (resp.\ $\psi''_1$)  in the
notation of (6.5) (resp. (6.11)). Composing this natural isomorphism
with $j_2,\ldots, j_n$ we obtain
$$j'_2,\ldots, j'_n:{\Cal A}_1\longrightarrow {\Cal A}_1\ \ .$$
Since $j_i^{'*}= -j'_i$, it follows that ${\roman{tr}}^0(j'_i)=0$ and hence
that $(j'_2,\ldots, j'_n)$ is an $(n-1)$-tuple of special
endomorphisms of ${\Cal A}_1$.

{\bf r=1.} Let ${\Cal A}$. Let $\j$ with $q(j_1)=\Delta$ or $p$ or
$\Delta p$. Then $j_1$ defines a special action $\iota: {\Cal
O}\to{\roman{End}}({\Cal A})$, where ${\Cal O}=\Z_p[j_1]$.
Furthermore $(j_2,\ldots, j_n)$ is an $(n-1)$-tuple of special
endomorphisms of $({\Cal A},
\iota)$.

We remark that the hereditary nature of these deformation problems
is a weak support of our conjecture on the singularities of the
intersections of special cycles in the case of isolated
intersections, comp.\ \cite{\krsiegel}, Conjecture 6.3. In our present
context of arithmetic Hirzebruch-Zagier cycles, the conjecture
states that for an isolated intersection point $\xi$ of the special
cycles $\CZ(T_1,\o_1),\ldots, \CZ(T_r, \o_r)$ we have
$$({\Cal O}_{\CZ(T_1,\o_1)}\otimes^{\Bbb L}\ldots \otimes^{\Bbb
L}{\Cal O}_{\CZ(T_r,\o_r)})_{\xi}= ({\Cal O}_{\CZ(T_1,\o_1)}
\otimes \ldots\otimes {\Cal O}_{\CZ(T_r,\o_r)})_{\xi}\ \ .\tag13.3$$
Hence the entity $e(\xi)$ in Proposition 6.2 would be the
intersection multiplicity in the sense of Serre.

We close this section with a curiosity. Let us start with an object
$({\Cal A},\lambda,\iota)$ of the deformation problem for $r=2$
where ${\Cal O}=\Z_{p^2}$. Let $\j$ with $j_1^2=p$ (the case where
$j_1^2=\Delta p$ is analogous). Then $\Z_{p^2}[j_1]= O_{B'_p}$ and
$\iota$ and $j_1$ define the structure $({\Cal A},\iota')$ of a
s.f.\ $O_{B'_p}$-module. However, for $i\geq 2$, the special
endomorphism $j_i$ {\it anticommutes} with the $O_{B'_p}$-action
(in the sense of the main involution of $O_{B'_p}$). Hence $j_i$ is
{\it not} a special endomorphism in the sense of {\bf r=1} {\it
($p$-adic uniformization case)} of the s.f.\ $O_{B'_p}$-module
$({\Cal A}, \iota')$. We have no explanation for this seeming
perplexity, except to point out that this remark seems to exclude a
simple method of creating bad singularities in the special cycles
in the good reduction cases considered above. Note that, in the
case {\bf r=1} {\it ($p$-adic uniformization case)},  the special
cycles can have bad singularities, e.g.\ have embedded components;
comp.\ \cite{\krdrin}.


\hfuzz=20pt
\widestnumber\key{666}

\vskip ,2in
\subheading{References}

\ref\key{\eichler}
\by M. Eichler
\book Quadratische Formen und orthogonale Gruppen
\publ Springer-Verlag
\publaddr Berlin
\yr 1974
\endref

\ref\key{\vandergeer}
\by G.\ van der Geer
\book Hilbert modular surfaces
\publ  Springer Verlag
\publaddr Berlin
\yr 1988
\endref

\ref\key{\grosskeating}
\by B. H. Gross and K. Keating
\paper  On the intersection of modular correspondences
\jour Invent.\ math.
\oldvol 112
\yr 1993
\pages 225--245
\endref

\ref\key{\hlr}
\by G.\ Harder, R.\ Langlands, M.\ Rapoport
\paper Algebraische Zyklen auf Hilbert-Blu\-men\-thal-Fl\"achen
\jour J.\ reine angew.\ Math.
\oldvol 366
\yr 1986
\pages 53--120
\endref

\ref\key{\hirzebruch}
\by F.\ Hirzebruch and G.\ van der Geer
\paper Lectures on Hilbert modular surfaces
\jour Sem.\ de Math.\ Sup.\ {\bf 77}, Presses Univ.\ de
Montr\'eal, Montreal 1981
\endref

\ref\key{\kaiser}
\by C. Kaiser
\paper Ein getwistetes fundamentales Lemma f\"ur die $\text{\rm GSp}_4$
\jour Dissertation, Bonn
\yr 1997
\endref

\ref\key{\kitaokatwo}
\by Y. Kitaoka
\paper A note on local representation densities of quadratic forms
\jour Nagoya Math. J.
\oldvol 92
\yr 1983
\pages 145--152
\endref

\ref
\key{\kitaoka}
\by Y. Kitaoka
\paper Fourier coefficients of Eisenstein series of degree 3
\jour Proc. of Japan Acad.
\oldvol 60
\yr 1984
\pages 259--261
\endref

\ref\key{\kottwitz}
\by R. Kottwitz
\paper Points on some Shimura varieties over finite fields
\jour JAMS
\oldvol 5
\yr 1992
\pages 373--444
\endref

\ref\key{\duke}   
\by S. Kudla
\paper Algebraic cycles on Shimura varieties of orthogonal type
\jour Duke Math. J.
\oldvol 86
\yr 1997
\pages 39--78
\endref

\ref\key{\annals}   
\by S. Kudla
\paper Central derivatives of Eisenstein series and height pairings
\jour Ann.\ of Math.
\oldvol 146
\yr 1997
\pages 545--646
\endref

\ref\key{\krsiegel}
\by S.\ Kudla and M.\ Rapoport
\paper Cycles on Siegel 3-folds and derivatives of Eisenstein series
\jour preprint K\"oln
\yr 1997
\endref

\ref\key{\krdrin}
\by S.\ Kudla and M.\ Rapoport
\paper Heights on Shimura curves and $p$-adic uniformization
\jour preprint
\yr 1998
\endref

\ref\key{\tiny}
\by S.\ Kudla, M.\ Rapoport, T.\ Yang
\paper On the derivative of an Eisenstein series of weight one
\jour to appear
\endref

\ref\key{\langlands}
\by R.P.\ Langlands
\book Base change for $GL(2)$
\bookinfo Ann.\ of Math.\ Studies 
\oldvol 96
\publ Princeton University Press
\publaddr Princeton, NJ
\yr 1979
\endref

\ref\key{\milne}
\by J. Milne
\paper Points on Shimura varieties mod p
\inbook Proc. Symp. Pure Math.
\oldvol 33, {\it Part 2}
\yr 1979
\pages 165--184
\publ AMS
\publaddr Providence, RI
\endref

\ref\key{\raynaud}
\by M.\ Raynaud
\paper Vari\'et\'es ab\'eliennes et g\'eom\'etrie rigide
\inbook Actes du Congr\`es International des Math\'ematiciens Nice
\pages 473-477
\publ Gauthier-Villars
\publaddr Paris
\yr 1971
\endref

\ref\key{\stamm}
\by H. Stamm
\paper On the reduction of the Hilbert-Blumenthal-moduli scheme with
$\Gamma_0(p)$-level structure
\jour Forum Math.
\oldvol 9
\yr 1997
\pages 405--455
\endref

\ref\key{\yang}
\by T.\ Yang
\paper An explicit formula for local densities of quadratic forms
\jour J. Number Theory
\oldvol 72
\pages 309--356
\yr 1998
\endref

\medskip\medskip
\noindent
\line{Stephen S.\ Kudla\hfill Michael Rapoport}
\line{Department of Mathematics\hfill Mathematisches Institut}
\line{University of Maryland\hfill der Universit\"at zu K\"oln}
\line{College Park, MD 20742\hfill Weyertal 86-90}
\line{{}\hfill D -- 50931 K\"oln}
\line{USA\hfill Germany}

\bye